\documentclass[11pt,a4paper]{article}

\usepackage[english]{babel}
\usepackage{amsmath}
\usepackage{amssymb}
\usepackage{graphicx}        
\usepackage{epstopdf}                             

\usepackage{amsfonts}
\usepackage{color}
\usepackage{verbatim}
\usepackage{caption}
\usepackage{subcaption}

\usepackage{array}

\usepackage{epsfig}
\usepackage{float}

\usepackage{algorithm}
\usepackage[noend]{algpseudocode}
\usepackage{url}

\usepackage{theorem}

\usepackage{algpseudocode} 
\usepackage{multirow}
\usepackage{amstext}
\usepackage{alltt}
\usepackage{verbatim}
\usepackage{array}
\usepackage{longtable}

\usepackage{fullpage}

\newtheorem{rem}{Remark}[section]

\newcommand{\EndProofMarker}{$\Box$}

\def\XXint#1#2#3{{\setbox0=\hbox{$#1{#2#3}{\int}$ }
\vcenter{\hbox{$#2#3$ }}\kern-.6\wd0}}

\graphicspath{{img/}}

\begin{document}

\title{Learning patient--specific parameters for a diffuse interface glioblastoma model
from neuroimaging data}  
\author{Abramo Agosti$^\ddag$, Pasquale Ciarletta$^\ddag$, Harald Garcke$^\sharp$, Michael Hinze$^*$}

\date{}

\maketitle 

\begin{center}
{\small $^\ddag$ MOX--Dipartimento di Matematica, 
        Politecnico di Milano \\
        via Bonardi 9, 20133 Milano, Italy\\
{\tt abramo.agosti@polimi.it}}        
\end{center}
\begin{center}
{\small $^\sharp$ Fakult\"at f\"ur Mathematik,
        Universit\"at Regensburg \\
        93040 Regensburg, Germany\\
{\tt harald.garcke@ur.de}}     
\end{center}
\begin{center}
{\small $^*$ Mathematisches Institut,
        Universit\"at Koblenz-Landau \\
        D-56070 Koblenz, Germany\\
{\tt hinze@uni-koblenz.de}}        
\end{center}
\maketitle 

  \abstract Parameters in mathematical models for glioblastoma
  multiforme (GBM) tumour growth are highly patient specific. Here we
  aim to estimate parameters in a Cahn--Hilliard type diffuse interface
  model in an optimised way using model order reduction (MOR) based on
  proper orthogonal decomposition (POD). Based on snapshots derived
  from finite element simulations for the full order model (FOM) we
  use POD for dimension reduction and solve the parameter estimation
  for the reduced order model (ROM). Neuroimaging data are used to
  define the highly inhomogeneous diffusion tensors as well as to
  define a target functional in a patient specific manner. The reduced
  order model heavily relies on the discrete empirical interpolation
  method (DEIM) which has to be appropriately adapted in order to deal
  with the highly nonlinear and degenerate parabolic PDEs. A feature
  of the approach is that we iterate between full order solves with
  new parameters to compute a POD basis function and sensitivity based
  parameter estimation for the ROM problems. The algorithm is applied
  using neuroimaging data for two clinical test cases and we can
  demonstrate that the reduced order approach drastically decreases
  the computational effort.



\section{Introduction}

Glioblastoma multiforme (GBM) is a malignant primary brain tumour
characterised by high infiltration into the parenchyma and wide
phenotypic heterogeneity \cite{leece}.  These characteristic features
of GBM provoke recurrence and marked resistance to adjuvant therapy,
resulting into poor prognosis and very low survival rates
\cite{ostrom}.  Thus, the emergent development of precision medicine
in neuro-oncology mainly concern the patient-specific optimisation of
the clinical treatment of GBM, with the aim to guide the decision
making of medical doctors for improving the quality of life of each
patient \cite{mirnezami2012preparing}.

In this context, mathematical models have proved useful as in-silico
benchmarks to improve the prognostic prediction and to tailor
personalised strategies in clinical practice
\cite{66,jackson2014mathematical}.  Most existing mathematical
approaches to neuro-oncology are based on reaction-diffusion partial
differential systems or agent-based models, that mimic the chemical
exchanges driving the tumour growth and transport properties of the
tumour cells as well as the response to adjuvant therapy
\cite{hatzikirou2005mathematical,jackson2015patient,prados2015toward,
  alfonso2017biology}.  Recent developments also enable the
possibility to integrate neuroimaging data in the virtual
reconstruction of the patient's brain, gaining insight on the effect
of the brain micro--structure on the invasive pathway
\cite{jbabdi2005simulation,65,painter2013mathematical,lipkova,MR3382724,MR3795322,MR3609861}. Despite the
great progress in assessing accurate mathematical predictions of the
prognostic clinical outcomes, the complexity underlying the physical
and biological cues driving GBM invasion make it particularly
difficult to quantify the accuracy of a given class of models in
reproducing the observable clinical events \cite{hawkins2013bayesian}.

In this work we propose a new strategy for optimising the parameter
estimation of a nonlinear degenerate diffuse interface model which
describes the GBM evolution integrating neuroimaging data, recently
proposed and analysed in \cite{Agosti1,Agosti2}. This partial
differential model consists of a Cahn--Hilliard equation with a
single--well potential of Lennard-Jones type, a non-conserved order
parameter and a degenerate mobility, that couples the growth of the
tumour phase with a reaction--diffusion equation for the  oxygen
concentration in the brain, including the effects of the standard
Stupp protocol of adjuvant therapy. In particular, it accounts for
the augmented tumour motility along white matter fibers tracts, which
is a typical hallmark of GBM, through the definition of heterogeneous
diffusion and chemotactic coefficients directly extracted from
Magnetic Resonance (MRI) and Diffusion Tensor (DTI) imaging
data. Because of the peculiar non-convexity and nonlinearity of the
chemical potential driving the local cell-cell interactions,  the
finite element approximation of the discrete model  has a high
computational cost, since it requires sophisticated numerical
techniques to select the physical solution representing the expanding
GBM boundary whilst avoiding numerical instabilities \cite{Agosti1}.\\

A proof-of-concept of the predictive ability of this model in clinical
practice has been presented in \cite{Agosti2}, showing by a manual
tuning of the model parameters how the numerical simulations on a
growing GBM tumour inside the virtual brain reconstructed by
segmentation of neuroimaging data could accurately fit the observed
invasion patterns observed at key clinical stages after surgical
removal and during adjuvant therapy. Here we propose a robust
automated procedure to optimise the parameter estimation by minimising
the $L^2$-distance between the indicator functions of the tumour
distribution sets in numerical simulations and the corresponding
clinical data of the GBM mass at a key time identified by our medical
collaborators in the clinical protocol. The constrained optimisation
problem is formulated at the numerical level as a Mathematical Problem
with Equilibrium Constraints (MPEC) \cite{luo}. Due to the high
computational cost of solving the Full Order Model (FOM), the
numerical solution of the MPEC will be given using model order
reduction. In particular, we propose an iterative algorithm extending
the one proposed in \cite{hinze2}, which computes a snapshot based
POD-ROM with the help of simulations 
at the FOM level and estimates parameter for the Reduced Order
Model (ROM) level through sensitivity analysis.  The main challenge
with respect to existing approach is the derivation of a ROM dealing
with a singular single-well potential and a degenerate anisotropic
mobility. The goal is to derive a robust iterative algorithm which
converges to an optimal state, that explores new regions in the
parameter space by changing the ROM basis and minimises the
optimisation functional at the ROM level avoiding the ROM solution to
violate the physical constraints satisfied by the full order
solution. For this scope, the bottleneck is the definition of an
effective order reduction of the degenerate and nonlinear terms of the
diffuse interface model.

The paper is organised as follows. In Section 2 we summarise the
diffuse interface model of GBM invasion and we derive the
corresponding FOM and ROM discretised problems. In Section 3, we
introduce a novel optimisation algorithm for parameter estimation. In
Section 4, we apply this algorithm using neuroimaging data
corresponding to two clinical test cases: the growth of a primary GBM
and the recurrence pattern after surgical resections.  The accuracy
and the computational gain of the proposed numerical procedure are
finally discussed in Section 5, together with few concluding remarks.

\section{Mathematical model}
In this section, we first summarise the diffuse interface model
employed for the patient-specific description of GBM evolution,
followed by the presentation of its FOM and ROM discretisations.
\subsection{The diffuse interface model}
The patient-specific GBM evolution is described using the diffuse
interface model proposed in \cite{Agosti1}. The model considers the
brain tissue as a saturated mixture composed by a tumour phase with a
volume fraction $\phi=\phi({\bf x},t)$ that expands at the expense of
another phase made of cells and fluids, so that $0 \le \phi \le
1$. The mass exchanges are regulated by the oxygen concentration
$n=n({\bf x},t)$ that is produced by the vascular network and consumed
by the tumour cells. This multi-phase framework has proved to give a
more realistic representation of the mechano-biological features
underlying the tumour growth processes
\cite{preziosi2009multiphase,wise2008three,42,garcke}.

Within the domain  $\Omega$ representing the brain, the mathematical model is given by  the following coupled PDEs:

\begin{equation}
\label{eqn:1}
\begin{cases}
  \frac{\partial \phi}{\partial t}= \nabla \cdot \biggl(
  \frac{\phi(1-\phi)^2}{M}\mathbf{T}\nabla \Sigma(\phi)\biggr)  +
  \Gamma_\phi (\phi,n) -\nabla \cdot \bigl(\chi_n\phi
  (1-\phi)^2\mathbf{T}\nabla n\bigr) &\;\text{in} \; \Omega\times(0,T), \\
  \frac{\partial {n}}{\partial t}=  \nabla \cdot (\mathbf{D}\nabla
  {n})+\Gamma_n (\phi,n) &\;\text{in} \; \Omega\times (0,T),\\
\end{cases}
\end{equation}
where $\mathbf{D}$ and $\mathbf{T}$ are the diffusion tensors of
oxygen and the tensor of preferential mobility, that  can be extracted
from neuroimaging data as in \cite{Agosti2}, $M$ is a friction
parameter that penalises the relative velocity between the phases, and
$\chi_n$  is a chemotactic coefficient, that  is considered to be $4$
times higher in the White Matter (WM) than its value in the Grey Matter
(GM) and in the Cerebrospinal Fluid (CSF). We set $\chi_n=k_n\chi$ where
$k_n$ is a chemotactic parameter which we need to estimate and
$\chi=4$ in the White Matter (WM) and $\chi=1$ in the Grey Matter (GM)
and the Cerebrospinal Fluid (CSF). 

The first equation in (\ref{eqn:1}) is a Cahn--Hilliard (CH) type
equation with degenerate mobility and non-conserved order parameter;
the chemical potential is defined by
\begin{equation}\label{eqn:sigma}\Sigma(\phi)= E(1-\phi_e)\psi_1'(\phi)+E\psi_2'(\phi)-\gamma^2\Delta\phi,  \\ 
\end{equation}
where $E$ is the Young modulus of the healthy tissue, $\gamma$ a
characteristic short-range interaction coefficient, and $\phi_e$ is
the homeostatic value of the volume fraction. We note that, due to the
non-smoothness in space of the chemotactic coefficient $\chi_n$, the
chemotactic term is inserted in the first equation of system
\eqref{eqn:1}, due to stability issues of the associated numerical
approximation, while it should be more natural to consider it as a
micro--force term associated to a coupling energy between cells and
nutrient directly in the form of the chemical potential $\Sigma$ as in
\cite{garcke}. The local interaction potential is given by:
\[
\psi_1(\phi):=-\log(1-\phi), \quad \psi_2(\phi):=-\frac{\phi^3}{3}-(1-\phi_e)\biggl(\frac{\phi^2}{2}+\phi\biggr),
\]
and it is split into a convex  $\psi_1$ and a non-convex $\psi_2$ term for future convenience. Such a functional form of the Lennard-Jones type describes attraction at low volume fraction and repulsion at beyond the homeostatic threshold $\phi_e$, as proposed in \cite{byrne2003modelling}. A simple functional form is given for the growth term $\Gamma_n$ for the oxygen : 
\[ \Gamma_n (\phi,n) =-\delta_n \phi {n} +S_n(1-{n})(1-\phi)\]
where $S_n$ is the production rate from the vascular network and $\delta_n$ gives the characteristic decay time. Similarly, the tumour growth rate is assumed in the form:
\[ \Gamma_\phi (\phi,n) = \nu\phi  (n-\delta ) (1-\phi) - k_{T}(t) \phi\]
where $\nu$ is the production rate mediated by the local oxygen concentration, and $k_T$ is a decay rate that accounts for apoptosis and/or adjuvant therapy.
The latter contribution is defined by the clinical Stupp protocol consisting  of radiotherapy and chemotherapy as in \cite{104}, reading:
\begin{equation}\label{eqn:2a}
k_{T}(t) \phi=k_{R}(t) \phi + k_C(t) \phi.
\end{equation}
The functions $k_R(t)$ and $k_C(t)$ are the temporal profiles of the radio- and chemo-therapy schedules, respectively: 
\begin{equation}
k_R(t)=
\begin{cases}
R_{eff}  & t_i \le t \le t_{i+1} \\
0 & \mbox{otherwise}
\end{cases}
\qquad
k_C(t)=
\begin{cases}
k_{C1}  & s_0 \le t \le s_1 \\
k_{C2}  & s_2 \le t \le s_3 \\
k_{C3}  & s_j \le t \le s_{j+1} \\
0 & \mbox{otherwise}
\end{cases}
\end{equation}
where the coefficients $k_{C1}$, $k_{C2}$ and $k_{C3}$ are extracted by clinical data and reflect the increasing drug dosage at different chemotherapy cycles, while $t_i,t_{i+1}$ are the days at the beginning and at the end of radiotherapy administration periods and $s_j,s_{j+1}$ are the days at the beginning and at the end of chemotherapy administration cycles.
The radiotherapy death rate  $R_{eff}$ is modelled via the linear-quadratic (LQ) model \cite{101} as
\begin{equation} \label{eq:Reff}
R_{eff}=\alpha m d + \beta m d^2 \, ,
\end{equation}
where $d$ [Gy] is the dose of radiation for every fraction, $m$ is the number of fractions per day (here, $m = 1$ day$^{-1}$), $\alpha$ [Gy$^{-1}$] and $\beta$ [Gy$^{-2}$] are two tissue-dependent parameters for cell kill \cite{102,110}. \\
\noindent The model (\ref{eqn:1}) is complemented by the following initial and boundary conditions:

\begin{equation}
\label{eqn:2bc}
\begin{cases}
\nabla \phi \cdot \boldsymbol{\nu}=\nabla \Sigma \cdot \boldsymbol{\nu}=\nabla n \cdot \boldsymbol{\nu}=0 &\;\text{on} \; \partial \Omega\times (0,T),\\
\phi(.,0)=\phi_0, \; n(.,0)=1 &\;\text{in} \; \Omega,
\end{cases}
\end{equation}
where $\boldsymbol{\nu}$ is the outer unit normal vector to
$\partial \Omega$ and $\phi_0$ is the initial distribution of tumour
concentration, that will be given by a neuroimaging datum at the
beginning of the clinical follow-up.
A list of reference biological ranges for the values of the parameters in \eqref{eqn:1} and the corresponding source is reported in Table \ref{tab:1}.\\

\begin{table}[!h]
\hspace*{+1.2cm}\resizebox{0.845\textwidth}{!}{\begin{minipage}{\textwidth}
\begin{tabular}{@{} l l l l @{}}
\hline
     & \textbf{Parameter description} & \textbf{Range of values} & \textbf{Ref.} \\ 
     \hline
$M$  & Tumour inter-phase friction & 1377.86\textendash 5032.2 (Pa day)/mm$^2$ & \cite{85}  \\
$\nu$ & Tumour cells proliferation rate &  0.012\textendash 0.5 day$^{-1}$ & \cite{65,gonzalez}  \\
$k_n$ & Chemotactic parameter & 0.007\textendash90.72 mm$^2$/(day) & \cite{Agosti1,Agosti2,ford} \\     
$S_n$ & Oxygen supply rate & $10^{4}$ day$^{-1}$ & \cite{42} \\
$\delta_n$ & Oxygen consumption rate & 8640 day$^{-1}$ &  \cite{frieboes,gonzalez}\\ 
  $\gamma$ & Diffuse interface thickness & 0.29\textendash 0.78  Pa$^{1/2}$mm &                                                        \cite{Agosti1} \\
$E$ & Brain Young modulus & 689.8  Pa\textendash 698.2 Pa & \cite{87}  \\
$\delta$ & Hypoxia threshold & 0.1\textendash 0.33 & \cite{bedogni,mascheroni} \\
$\phi_e$ & Equilibrium cell volume fraction & 0.389\textendash 0.8 & \cite{88}  \\
\hline
$m$ & Radiation fractions per day & 1 day$^{-1}$ & \cite{105}  \\
$d$ & Radiation dose & 2 Gy & \cite{105}  \\
$\alpha$ & Linear coefficient for RT induced cell kill &   0.027 Gy$^{-1}$ & \cite{102,104} \\
$\beta$  & Quadratic coefficient for RT induced cell kill  & 0.0027 Gy$^{-2}$ &  \cite{104,112,113}\\  
$k_{C1}$ & Concomitant chemotherapy death rate & 0.00735 day$^{-1}$ & \cite{104} \\
$k_{C2}$ & First cycle of adjuvant CHT death rate & 0.0147 day$^{-1}$ & \cite{104} \\
$k_{C3}$ & Remain. cycles of adjuvant CHT death rate & 0.0196 day$^{-1}$ & \cite{104} \\
\hline
\end{tabular}
\end{minipage} }
\caption[Model parameters]{Values or ranges of values for parameters used in the model.
}
\label{tab:1}
\end{table}

\noindent
As first highlighted in \cite{bernis}, the presence of compactly
supported weak solutions of the degenerate Cahn--Hilliard equation in
\eqref{eqn:1} is linked to the non-uniqueness of the solution, with
the existence of physical solutions with moving free boundary at a
finite speed and unphysical ones with fixed support in time. Due to
the degeneracy of the mobility and the logarithmic singularity of the
potential, we remark that a weak solution of \eqref{eqn:1},
\eqref{eqn:2bc} satisfies the positivity and separation constraints
\[
0\leq \phi <1 \quad \text{a.e} \; \text{in} \; \Omega\times (0,T).
\]  
For future convenience, we finally define the following set of model parameters
\[
\mathcal{P}\equiv\{L:=1/M,\nu,k_n,S_n,\delta_n,\gamma^2,E,\delta,c_e:=1-\phi_e\},
\]
 whose optimal value will be searched for in the parameters optimisation problem introduced in the following sections.
 
\subsection{FOM discretisation}
In order to discretise the initial boundary value problem given by (\ref{eqn:1}), (\ref{eqn:2bc}) we start from the MRI and DTI data collected at the initial time $t=0$ to generate the tetrahedral mesh for the discretised geometry of the brain $\mathcal{T}_h$, the map which specifies the location of the brain tissues, the initial tumour distribution and the tensors $\mathbf{D}$ and $\mathbf{T}$. \\
\noindent We solve the following \textbf{Initialisation} problem:\\ 
\textbf{Problem 1}\\
\emph{
\textbf{Initialisation:}
Given MRI$(t=0)$ and DTI$(t=0)$, determine the set
\begin{equation}
\label{eqn:initialization}
\begin{cases}
\mathcal{T}_h,\\
\text{map(WM,GM,CSF)},\\
\phi_0,\\
\mathbf{D},\mathbf{T}
\end{cases}
=\textbf{Initialization}\text{(MRI(t=0), DTI(t=0))}.
\end{equation}
}

The \textbf{Initialisation} procedure is defined as follows.  Starting
from a segmentation of the brain tissues and of the initial tumour
distribution, we extract the external brain surface and generate the
computational surface and the tetrahedral internal mesh, conveniently
refined around the tumour region. We also generate the labelled mesh
$\text{map}(WM,GM,CSF)$ which maps each cell of the mesh to an integer
value identifying the cerebral tissue the cell barycentre belongs to,
together with the map for the characteristic function of the tumour
distribution.  The initial tumour distribution $\phi_0$ is then
obtained as the characteristic function of the tumour extension
multiplied by the equilibrium value $\bar{\phi}$. The value
$\bar{\phi}$ is obtained by solving $\Gamma_\phi(\bar{\phi},\bar{n}) =
0 = \Gamma_n(\bar{\phi},\bar{n})$ with $k_T(0)=0$,
i.e. $\bar{n}=\delta$ and 
\[\bar{\phi}=\frac{S_n(1-\delta)}{S_n+\delta(\delta_n-S_n)}.\]
This means that we are considering that the tumour phase is in an
homogeneous equilibrium value in the tumour support at the initial
temporal stage, before the application of the therapy. We finally
generate the meshes containing the values of the independent
components of the
tensors $\mathbf{D}$ and $\mathbf{T}$ by analysing the log--signal associated to the DTI data. The extrapolation processes performed in \textbf{Problem 1} are described in detail in \cite{Agosti1,Agosti2}.\\
The FOM discretisation of system \eqref{eqn:1} is obtained through the
finite element method, by following the approach introduced in
\cite{Agosti1,BBG}. Let us introduce the following finite element spaces:
\begin{align}
\notag & Q_{h} := \{\chi \in C(\bar{\Omega}):\chi |_{K}\in \mathbb{P}^{1}(K) \; \forall K\in \mathcal{T}_{h}\}\subset H^{1}(\Omega),\\
\notag & Q_{h}^+ := \{\chi \in Q_{h}: \chi \geq 0\; \rm in \, \Omega\}
\end{align}
where $\mathbb{P}_{1}(K)$ indicates the space of polynomials of total order one on $K$.
We set $\Delta t = T/N$ for a $N \in \mathbb{N}$ and $t_{n}=n\Delta t$, {$n=0,...,N$}.
Starting from initial data $\phi_{h}^{0}=\pi^{h}\phi_{0}$, $n_{h}^{0}=\pi^{h}n_{0}$, where $\pi^h$ is the standard Lagrangian interpolation operator, and from a set of parameters $\mathcal{P}_k$, we consider the following \textbf{FOM} discretisation. (Note that if no mesh vertex is on the boundary of the tumour support, the Lagrangian projection $\phi_h^{0}=\pi^h(\phi_0)$ is well defined. Otherwise, we take a smoothing of $\phi_0$ such that $\phi_0\in C(\bar{\Omega})$.)\\
\textbf{Problem 2}\\
\emph{
\textbf{FOM:} \\ 
For $n=1,\dots,N$, given $(\phi_h^{n-1}, n_h^{n-1}) \in Q_h^+\times Q_h^+$, we define
\begin{equation}
\label{eqn:fom}
(\phi_h^{n},    \Sigma_h^{n},n_h^{n})=\textbf{F}_k^n(\textbf{Initialization},\mathcal{P}_k)=\textbf{F}_k^n(\mathcal{T}_h,\text{map}(WM,GM,CSF),\phi_0,\mathbf{D},\mathbf{T},\mathcal{P}_k),
\end{equation}
where $(\phi_h^{n},    \Sigma_h^{n},n_h^{n}) \in Q_h^+ \times Q_h \times Q_h^+$ satisfies, $\forall (v_h,w_h,q_h) \in Q_h \times     Q_h^+ \times Q_h$,
    \begin{equation}
    \label{eqn:3}
    \begin{cases}
      \displaystyle \biggl ( \frac{\phi_h^{n} -\phi_h^{n-1}}{\Delta t}
      , v_h \biggr )^h =& -L_k ( \phi_h^{n-1} (1-\phi_h^{n-1})^2
      \mathbf{T} \nabla \Sigma_h^{n} , \nabla v_h ) + \nu_k (
      \phi_h^{n-1} (n_h^{n}-\delta_k)(1-\phi_h^{n-1}) , v_h )^h \\&
      +k_{nk} ( \chi \phi_h^{n-1} (1-\phi_h^{n-1})^2 \mathbf{T} \nabla
      n_h^{n} , \nabla v_h )-(k_T\phi_h^{n-1},v_h)^h,
      \\
      \displaystyle \gamma^2_k(\nabla
      \phi_h^n,\nabla(w_h-\phi_h^n))+&(E_kc_{ek}\psi'_1(\phi_h^n),w_h-\phi_h^n)^h
      \ge \biggr(\Sigma_h^n+E_k\psi'_2(\phi_h^{n-1}),w_h-\phi_h^n\biggl)^h, \\
      \displaystyle \biggl ( \frac{n_h^{n} -n_h^{n-1}}{\Delta t} , q_h
      \biggr )^h =& - ( \mathbf{D} \nabla n_h^{n} , \nabla q_h )
      +S_{nk}( (1-n_h^{n})(1-\phi_h^{n-1}) , q_h )^h -\delta_{nk} (
      \phi_h^{n-1} n_h^{n} , q_h )^h.
    \end{cases}
    \end{equation}
    }
    
    The finite element approximation \eqref{eqn:3} has the form of a
    discrete variational inequality, where the positivity of the
    discrete solution is enforced as a constraint. Indeed, $\phi_h^n$
    is projected onto the space with positive values $Q_h^+$. This
    projection aims at recovering the analytical properties of the
    continuous solution \cite{Agosti1,BBG}.
\begin{rem}
\label{rem:equality}
For a solution of \eqref{eqn:3}, we obtain that 
$\phi_h^n(\mathbf{x}_j)\ge 0$ and if $\phi^n_h(\mathbf{x}_j)>0$ holds the
equality
\begin{equation}\label{eq:sigma}
  \gamma^2_k(\nabla \phi_h^n,\nabla
  \chi_j)+(E_kc_{ek}\psi'_1(\phi_h^n)+E_k\psi'_2(\phi_h^{n-1})-\Sigma_h^n,\chi_j)^h=0
\end{equation}  
is true, with $\chi_j$ the basis function associated to the node $j$.
To prove this we choose in the case $\phi^n_h(\mathbf{x}_j)>0$ in
\eqref{eqn:3} the test functions $w_h=\phi^n_h\pm \frac12
\phi^n_h(\mathbf{x}_j)\chi_j$ and obtain \eqref{eq:sigma} as
$\phi^n_h(\mathbf{x}_j)>0$. This shows that $\Sigma^n_h(\mathbf{x}_j)$
is uniquely defined if $\phi^n_h(\mathbf{x}_j)>0$.
\end{rem}
The lumped mass approximation of the $L^2$-scalar product is
introduced in \eqref{eqn:3} in order for the discrete solution to be
able to track compactly supported solutions with a moving free
boundary from the unphysical ones with fixed support.
\\
Note that the convex part of the cellular potential is treated
implicitly in time, whereas the concave part is treated
explicitly. The well posedness of system \eqref{eqn:3} can be shown
expressing its equations as the KKT conditions of a convex
minimisation problem with regular inequality constraints. It is solved
by a null--space gradient projection algorithm.

\subsection{ROM discretisation}
We use Proper Orthogonal Decomposition (POD) reduced order modelling to
obtain a ROM of the full order system \eqref{eqn:3}. We follow
\cite{hinze1}, which formulates a ROM of a Cahn--Hilliard equation with constant
mobility and advection by projecting the equations for $\phi$ and
$\Sigma$ onto the spaces spanned by the POD basis obtained from the
time snapshots matrices associated to $\phi$ and $\Sigma$,
respectively, and which approximates the nonlinear convex term in the
double--well potential using the Discrete Empirical Interpolation
Method (DEIM). In the present case, this formalism is extended to take
into account for the degeneracy of the mobility and the physical
constraints expected to be satisfied by the ROM solution, as well as
for the presence of anisotropy and the presence of the dynamics of the
nutrient. In particular, our strategy is to use DEIM interpolation to
approximate the nonlinear terms $\psi_1'$ and $\psi_1''$, (the latter
term is needed to solve the nonlinear ROM system by a Newton method), and to project the equations for $\phi$, $\Sigma$ and $n$ in \eqref{eqn:3} onto the space spanned by the POD basis associated to $\phi$, $\Sigma$ and $n$, respectively. The mobility and chemotactic terms, which contain polynomial nonlinearities in the variable $\phi$, are written as higher order tensors acting on the reduced order space associated to $\phi$.\\
We start by calculating the POD basis obtained from the snapshots matrices associated to $\phi$, $\Sigma$, $n$, $\psi_1'(\phi)$ and $\psi_1''(\phi)$. {In what follows we do not distinguish between a finite element function $f_h$ and its nodal vector and it will be clear from the context which object is meant. For a given snapshot matrix $F=\left[f_h^0,\dots,f_h^N\right]$ the POD basis elements $\xi_s^{f}$ of the POD basis $\left\{\xi_l^f\right\}_{l=1,\dots,N_{\text{POD}}}$ are obtained as follows: 
\begin{itemize}
\item prescribe the required information content to be covered by the POD basis as $ic\in (0,1]$;
\item compute the trace $tr(F^tF)$ of the correlation matrix $F^tF=(f_h^m,f_h^l)_{ml}\in M(N+1,\mathbb{R})$, where $(\cdot,\cdot)$ denotes the chosen inner product;
\item set $N_f^{\text{POD}}:= \min\left\{m,\left(\sum\limits_{i\le m} \lambda_i\right)/tr(F^tF) \ge ic\right\}$;
\item (successively) compute the eigensystem $\{v^i,\lambda_i\}_{i=1,\dots,N_f^{\text{POD}}}$ of $F^tF$;
\item set $\xi_s^{f}:=\frac{1}{\sqrt{\lambda_s}}\sum\limits_{j}v^s_j f^j_h$ $(1\le s \le N_f^{\text{POD}})$.
\end{itemize}
With this procedure we set up the POD bases for our reduced order
models, whose construction is summarised in \\ \textbf{Problem 3}:\\
\emph{
\textbf{POD:}\\
With the given time snapshots form the corresponding snapshot matrices in correspondence of the parameter set $\mathcal{P}_k$ according to $\mathbf{F}_1:=\left[\phi_h^0,\dots,\phi_h^N\right]$, $\mathbf{F}_2:=\left[\Sigma_h^0,\dots,\Sigma_h^N\right]$, $\mathbf{F}_3:=\left[n_h^0,\dots,n_h^N\right]$, $\psi_1'(\mathbf{F}_1):=\left[\psi_1'(\phi_h^0),\dots,\psi_1'(\phi_h^N)\right]$, and $\psi_1''(\mathbf{F}_1):=\left[\psi_1''(\phi_h^0),\dots,\psi_1''(\phi_h^N)\right]$. Then compute the POD systems corresponding to the correlation matrices matrices $(\mathbf{F}_1)^T\mathbf{F}_1$, $(\mathbf{F}_2)^T\mathbf{F}_2$, $(\mathbf{F}_3)^T\mathbf{F}_3$, $(\psi_1'(\mathbf{F}_1))^T\psi_1'(\mathbf{F}_1)$, $(\psi_1''(\mathbf{F}_1))^T\psi_1''(\mathbf{F}_1)$ and set
\[
N_{\text{POD}}:=\max\{N_{\phi}^{\text{POD}}, N_{\Sigma}^{\text{POD}},  N_{n}^{\text{POD}},N_{\psi_1'}^{\text{POD}}, N_{\psi_1''}^{\text{POD}}\}.
\]
If the respective bases are computed successively complete them by adding the remaining $N_{\text{POD}}-N_{\theta}^{\text{POD}}$ basis elements to the basis $\left\{\xi_l^{\theta}\right\}_{l=1,\dots,N_{\theta}^{\text{POD}}}$, where $\theta \in \{\phi,\Sigma,n,\psi_1',\psi_2''\}$. Finally assemble the respective bases in the array
\begin{align}
\label{eqn:pod}
\mathbf{P}_k:=\left(\mathbf{P}_{k1},\mathbf{P}_{k2},\mathbf{P}_{k3},\mathbf{P}_{k4},\mathbf{P}_{k5}\right),
\end{align}	
where
\begin{multline*}
\mathbf{P}_{k1}:=\left\{\xi_l^{\phi}\right\}_{l=1,\dots,N_{\text{POD}}},\mathbf{P}_{k2}:=\left\{\xi_l^{\Sigma}\right\}_{l=1,\dots,N_{\text{POD}}},\mathbf{P}_{k3}:=\left\{\xi_l^{n}\right\}_{l=1,\dots,N_{\text{POD}}},\\
\mathbf{P}_{k4}:=\left\{\xi_l^{\psi_1'}\right\}_{l=1,\dots,N_{\text{POD}}},\text{
and }\mathbf{P}_{k5}:=\left\{\xi_l^{\psi_1''}\right\}_{l=1,\dots,N_{\text{POD}}}.
\end{multline*}
}
We note that POD basis elements are finite element functions. In the numerical examples we will specify $N_{\phi}^{\text{POD}}$, $N_{\Sigma}^{\text{POD}}$,
$N_{n}^{\text{POD}}$, $N_{\psi_1'}^{\text{POD}}$, $N_{\psi_1''}^{\text{POD}}$ such that the required information contents of the POD bases satisfy $ic=0.9999$, i.e. all the POD bases contain at least 99.99\% of the snapshot information.

}
Then we make the ansatz
\begin{equation}
 \label{eqn:4}
 \phi_h^n=\sum_{i=1}^{N_{\text{POD}}}\alpha_{ik}^n\xi_{i}^{\phi}, \quad \Sigma_h^n=\sum_{i=1}^{N_{\text{POD}}}\beta_{ik}^n\xi_{i}^{\Sigma}, \quad n_h^n=\sum_{i=1}^{N_{\text{POD}}}\eta_{ik}^n\xi_{i}^{n}.
\end{equation}
We moreover approximate the singular nonlinear terms $\psi_1'(\phi_h^n)$ and $\psi_1''(\phi_h^n)$ by a greedy algorithm using DEIM interpolation \cite{deim}, i.e. by computing the nonlinearities only on the nodes of the mesh which give the greatest interpolation contribution for each of their POD basis elements,  
\begin{align}
\label{eqn:deim1}
& \psi_1'\biggl(\sum_{i=1}^{N_{\text{POD}}}\alpha_{ik}^n\xi_{i}^{\phi}\biggr)=\sum_{i=1}^{N_{\text{POD}}}(P_{2,ij}^TU_{\psi_1',js})^{-1}\psi_1'(P_{2,sl}^T\Phi_{lm}\alpha_{mk}^n)\xi_i^{\psi_1'},\\
\label{eqn:deim2}
& \psi_1''\biggl(\sum_{i=1}^{N_{\text{POD}}}\alpha_{ik}^n\xi_{i}^{\phi}\biggr)=\sum_{i=1}^{N_{\text{POD}}}(P_{2,ij}^TU_{\psi_1'',js})^{-1}\psi_1''(P_{2,sl}^T\Phi_{lm}\alpha_{mk}^n)\xi_i^{\psi_1''},
\end{align}
where $i,j,s,l,m=1,\dots,N_{\text{POD}}$, $U_{\psi_1'}:=\left[\xi_{1}^{\psi_1'}, \dots, \xi_{N_{\text{POD}}}^{\psi_1'}\right]$, $U_{\psi_1''}:=\left[\xi_{1}^{\psi_1''}, \dots, \xi_{N_{\text{POD}}}^{\psi_1''}\right]$, $\Phi:=\left[\xi_{1}^{\phi}, \dots, \xi_{N_{\text{POD}}}^{\phi}\right]$ and $P_2$ is the row selection operator of the DEIM interpolation nodes for the term $\psi_1''(\cdot)$. The DEIM algorithm is reported in the Appendix. We observe that we are computing both \eqref{eqn:deim1} and \eqref{eqn:deim2} on the same interpolation nodes (associated to the $\psi_1''(\cdot)$ term only). This is needed to practically enforce the separation property $\sum_{i=1}^{N_{\text{POD}}}\alpha_{ik}^{n,p}\xi_{i}^{\phi}<1$, driven by the singularities in the terms $\psi_1'(\cdot)$ and $\psi_1''(\cdot)$, throughout the steps of the Newton algorithm.\\
Substituting \eqref{eqn:4} and \eqref{eqn:deim1} into \eqref{eqn:3}, choosing $v_h\equiv \xi_m^{\phi}$, $w_h\equiv \xi_m^{\Sigma}$, $q_h\equiv \xi_m^{n}$, we - after replacing the original variational inequality by an equation as in the original equation \eqref{eqn:sigma} - obtain the following \textbf{ROM} system
\begin{equation}
\label{eqn:6}
\begin{cases}
\displaystyle V_{1,mi}\frac{\alpha_{ik}^n-\alpha_{ik}^{n-1}}{\Delta t}&=-L_k\bigl(\alpha_{ik}^{n-1}\alpha_{jk}^{n-1}\alpha_{sk}^{n-1} V_{2,ijsml}\beta_{lk}^{n} -2\alpha_{ik}^{n-1}\alpha_{jk}^{n-1}V_{3,ijms}\beta_{sk}^{n}+\alpha_{ik}^{n-1}V_{4,imj}\beta_{jk}^{n}\bigr)+\\
\displaystyle  &k_{nk}\bigl(\alpha_{ik}^{n-1}\alpha_{jk}^{n-1}\alpha_{sk}^{n-1} V_{8,ijsml}\eta_{lk}^{n} -2\alpha_{ik}^{n-1}\alpha_{jk}^{n-1}V_{9,ijms}\eta_{sk}^{n}+\alpha_{ik}^{n-1}V_{10,imj}\eta_{jk}^{n}\bigr)+\\
\displaystyle &\nu_k (\alpha_{ik}^{n-1} V_{5,ims}-\alpha_{ik}^{n-1} \alpha_{jk}^{n-1} V_{6,ijms})\eta_{sk}^n -(\nu_k \delta_k + K_T)V_{1,mi}\alpha_{ik}^{n-1}+\nu_k \delta_k \alpha_{ik}^{n-1} V_{7,ims}\alpha_{sk}^{n-1},\\ \\
\displaystyle U_{1,mi}\beta_{ik}^n&=\gamma_k^2U_{6,mi}\alpha_{ik}^n+E_kc_{ek}U_{2,mi}(P_{2,ij}^TU_{\psi_1',js})^{-1}\psi_1'(P_{2,sl}^T\Phi_{lm}\alpha_{mk}^n)-E_k\alpha_{ik}^{n-1}U_{3,imj}\alpha_{jk}^{n-1}-\\
\displaystyle & E_kc_{ek}U_{4,mi}\alpha_{ik}^{n-1}-E_kc_{ek}U_{5,m},\\ \\
\displaystyle  W_{1,mi}\frac{\eta_{ik}^n-\eta_{ik}^{n-1}}{\Delta t}&=-W_{2,mi}\eta_{ik}^n+S_{nk} \bigl(\alpha_{ik}^{n-1}W_{3,imj}\eta_{jk}^n+W_{4,m}-W_{5,mi}\alpha_{ik}^{n-1}-W_{1,mi}\eta_{ik}^{n-1}\bigr)-\\
\displaystyle & \delta_{nk}\alpha_{ik}^{n-1}W_{3,imj}\gamma_{jk}^{n}.
\end{cases}
\end{equation}
Here, $i,j,s,l,m=1,\dots,N_{\text{POD}}$, and the initial conditions are given by $\alpha_{ik}^0=(\phi_h^0,\xi_i^{\phi})^h$ and $\eta_{ik}^0=(n_h^0,\xi_i^{n})^h$. The second order and higher order tensors $V_1,V_2,\dots,V_{10}$, $U_1,U_{2},\dots,U_6$, $W_1,\dots,W_5$ in \eqref{eqn:6} and $U_{7}$ (needed for the Newton problem) are defined by the following \textbf{Assemble} problem in terms of $\mathbf P_k$ defined in \eqref{eqn:pod}. 
\\
\textbf{Problem 4}\\
\emph{
\textbf{Assemble:}
\begin{equation}
\label{eqn:assemble}
\mathbf{A}_k(\mathbf{P}_{k1},\mathbf{P}_{k2},\mathbf{P}_{k3},\mathbf{P}_{k4},\mathbf{P}_{k5}):= (V_1,\dots ,V_{10}, U_1,U_{2},\dots ,U_6, W_1,\dots ,W_5, U_{7}),
\end{equation}
where
\[
V_{1,ji}:=(\xi_i^{\phi}, \xi_j^{\phi})^h, \quad U_{1,ji}:=(\xi_i^{\Sigma}, \xi_j^{\Sigma})^h, \quad U_{2,ji}:=(\psi_{1,i}', \xi_j^{\Sigma})^h,\] 
\[
U_{4,ji}:=(\xi_i^{\phi}, \xi_j^{\Sigma})^h, \quad U_{5,i}:=(1, \xi_i^{\Sigma})^h, \quad U_{6,ji}:=(\nabla \xi_i^{\phi}, \nabla \xi_j^{\Sigma}), 
\]
\[
W_{1,ji}:=(\xi_i^{n}, \xi_j^{n})^h, \quad W_{2,ji}:=(\mathbf{D}\nabla \xi_i^{n}, \nabla \xi_j^{n}), \quad
W_{4,i}:=(1, \xi_i^{n})^h, \quad W_{5,ji}:=(\xi_i^{\phi}, \xi_j^{n})^h,
\]
for the second order tensors corresponding to the bilinear forms in \eqref{eqn:3}, with $i,j=1,\dots,N_{\text{POD}}$.\\
For the higher order tensors of the polynomial nonlinear forms we for $i,j,k,l,m=1,\dots,N_{\text{POD}}$ find 
\[
V_{2,ijkml}:=(\xi_i^{\phi}\xi_j^{\phi}\xi_k^{\phi}\mathbf{T}\nabla \xi_l^{\Sigma},\nabla \xi_m^{\phi}), \quad  V_{3,ijlk}:=(\xi_i^{\phi}\xi_j^{\phi}\mathbf{T}\nabla \xi_k^{\Sigma},\nabla \xi_l^{\phi}), \quad V_{4,ikj}:=(\xi_i^{\phi}\mathbf{T}\nabla \xi_j^{\Sigma},\nabla \xi_k^{\phi}),\quad\] 
for the mobility term;
\[
V_{8,ijkml}:=(\chi\xi_i^{\phi}\xi_j^{\phi}\xi_k^{\phi}\mathbf{T}\nabla \xi_l^{n},\nabla \xi_m^{\phi}), \quad  V_{9,ijlk}:=(\chi\xi_i^{\phi}\xi_j^{\phi}\mathbf{T}\nabla \xi_k^{n},\nabla \xi_l^{\phi}), \quad V_{10,ikj}:=(\chi\xi_i^{\phi}\mathbf{T}\nabla \xi_j^{n},\nabla \xi_k^{\phi}),\quad\] 
for the chemotactic term;
\[
V_{5,ikj}:=(\xi_i^{\phi}\xi_j^{n}, \xi_k^{\phi})^h, \quad V_{6,ijlk}:=(\xi_i^{\phi}\xi_j^{\phi}\xi_k^{n}, \xi_l^{\phi})^h, \quad V_{7,ijk}:=(\xi_i^{\phi}\xi_j^{\phi}, \xi_k^{\phi})^h,\] 
for the source term for $\phi$, 
\[
U_{7,ikj}:=(\psi_{1,i}''\xi_j^{\phi}, \xi_k^{\Sigma})^h, \quad U_{3,ikj}:=(\xi_i^{\phi}\xi_j^{\phi}, \xi_k^{\Sigma})^h,\] 
for the terms containing $\psi_1''$ and $\psi_2'$, and
\[
W_{3,ikj}:=(\xi_i^{\phi}\xi_j^{n}, \xi_k^{n})^h
\]
for the source term for $n$.
} 
 
We highlight that it is of utmost importance to assemble the higher order
tensors $V_2$ , $V_3$ , $V_4$ and $V_8$ , $V_9$ , $V_{10}$, which project the mobility and chemotactic
terms onto the ROM space, instead of approximating them through tensor
interpolation. The former approach avoids to loose informations about the
anisotropy  of the tensor of preferential directions \textbf{T} and the heterogeneity of the chemotactic function $\chi$.\\

Due to the nonlinearity in the term $\psi_1'$, we solve \eqref{eqn:6} by means of the Newton method, defining the \textbf{ROM Newton} problem: \\
\textbf{Problem 5}\\ 
\emph{
\textbf{ROM Newton:}
\begin{equation}
\label{eqn:romnewton}
(\alpha_{ik}^n,\beta_{ik}^n,\eta_{ik}^n)_{i=1,\dots,N_{\text{POD}};n=0,\dots,N}=\mathbf{RN}_k^n(\mathbf{A}_k,\mathcal{P}_k,\phi_h^0,n_h^0).
\end{equation}
}

The Newton algorithm, which defines the function $\mathbf{RN}_k^n$ is reported in the Appendix.\\
\noindent We finally search for the solutions of the nine \textbf{ROM linearised} systems obtained from varying $\mathcal{P}_k=\mathcal{P}_k+\delta \mathcal{P}$ in \eqref{eqn:6}, defining\\ 
\textbf{Problem 6}\\
\emph{
\textbf{ROM linearised:}
\begin{align}
\label{eqn:romlinearized}
& \notag \nabla_{\mathcal{P}_k}\vec{\alpha}_k=\biggl[\frac{\partial \vec{\alpha}_k}{\partial L_k}\, \frac{\partial \vec{\alpha}_k}{\partial \nu_k}\, \frac{\partial \vec{\alpha}_k}{\partial k_{nk}}\, \frac{\partial \vec{\alpha}_k}{\partial S_{nk}}\, \frac{\partial \vec{\alpha}_k}{\partial \delta_{nk}}\, \frac{\partial \vec{\alpha}_k}{\partial \gamma_k^2}\, \frac{\partial \vec{\alpha}_k}{\partial E_k}\, \frac{\partial \vec{\alpha}_k}{\partial \delta_k}\, \frac{\partial \vec{\alpha}_k}{\partial c_{ek}}\biggr]^t=\\
& \mathbf{RL}_k(\mathbf{A}_k,\mathcal{P}_k,(\mathbf{RN}_k^n)_{n=0,\dots,N}).
\end{align}  
}

The construction of the functions $\mathbf{RL}_k$ is reported in the Appendix.


\section{Optimisation algorithm}
{In the following we  propose an algorithm to learn the parameters of our model from clinical neuroimaging data. For this purpose we set up a minimisation  problem for the model parameters which we solve iteratively with the help of reduced order models aligned with the parameter sets associated to the respective iteration. We refer to \cite{AH01} where this concept was
proposed for optimal flow control using POD surrogate models.

The set of model parameters is given by
\[
\mathcal{P}=\{L,\nu,k_n,S_n,\delta_n,\gamma^2,E,\delta, c_e\}.
\]
To formulate the parameter learning problem as an optimisation problem we introduce the functional
}
\begin{equation}
\label{eqn:2}
J(\phi(\mathcal{P}),\mathcal{P})=\frac{1}{2||H(\phi_{\text{data}}(T))||_{L^2(\Omega)}^2}||H_{\phi_e}(\phi(T))-H(\phi_{\text{data}}(T))||_{L^2(\Omega)}^2+\frac{\eta}{2}\sum_{m=1}^{|\mathcal{P}|}\biggl(\frac{\mathcal{P}_{m}-\mathcal{P}_{\text{exp,m}}}{\mathcal{P}_{\text{exp,m}}}\biggr)^2,
\end{equation}
where $\phi$ is a solution of the system \eqref{eqn:1} supplemented with
the initial and boundary conditions \eqref{eqn:2a}. The set 
$\mathcal{P}_{\text{exp}}$ contains estimates of expected values for the
parameters, $\eta$ is a regularisation parameter and $T$ is the time
at which we compare the simulations and the data tumor extensions. A
regularised Heaviside function $H_{\phi_e}$ with slope $2/\phi_e$ is
used to approximate the characteristic function of the tumor extension
\[
H_{\phi_e}(\phi):=\begin{cases}
1, \quad \text{if} \; \phi\geq \phi_e/2,\\
2\phi/\phi_e \quad \text{if} \; \phi\geq 0 \; \text{and} \; \phi\leq \phi_e/2,\\
0 , \quad \text{if} \; \phi\leq 0.
\end{cases}
\]
This means that we are considering $\phi=\phi_e/2$ as the equation for
the hypersurface defining the boundary of the tumour extension, since
we can assume that the tumour profile given as a solution of the
degenerate Cahn--Hilliard equation (without growth) with the single well cellular potential has a kink--like form between the two equilibria $\phi=0$ and $\phi=\phi_e$. Finally, $H(\phi_{\text{data}}(T))$ is the characteristic function of the tumour extension from data, computed through the \textbf{Target} problem: \\ 
\textbf{Problem 7}\\
\emph{
\textbf{Target:}\\
Given MRI($t=0$), determine
\begin{equation}
\label{eqn:target}
H(\phi_{\text{data}}(T))=\mathbf{Target}(\text{MRI}(t=T)),
\end{equation}
where, given a segmentation map of the tumour extension from MRI($t=T$), the function $H(\phi_{\text{data}}(T))$ takes the value $1$ on the tumour map and the value $0$ outside.\\
}
The functional $J$ measures the $L^2(\Omega)$ distance between the characteristic functions of the tumor extensions from simulations and data. Note that, due to the fact that we cannot easily obtain informations about the tumor cell densities from the MRI images, we are going to consider the distance between the characteristic functions of the tumor extensions from simulations and data, and not the distance between the tumor distributions.\\
The value of the functional $J$ calculated at the FOM level is
\begin{equation}
\label{eqn:2fom}
J(\mathbf{F}_{1k},\mathcal{P}_k)=\frac{1}{2||\pi_h H(\phi_{\text{data}}(T))||_{h}^2}||\pi_h[H_{\phi_e}(\phi_h^N)- H(\phi_{\text{data}}(T))]||_{h}^2+\frac{\eta}{2}\sum_{m=1}^{|\mathcal{P}|}\biggl(\frac{\mathcal{P}_{k,m}-\mathcal{P}_{\text{exp,m}}}{\mathcal{P}_{\text{exp,m}}}\biggr)^2,
\end{equation}
where $||f_h||_h^2=(f_h,f_h)_h=\{f_h\}^TM_h\{f_h\}$ is the lumped $L^2(\Omega)$ norm, with $M_h$ the lumped mass matrix.
We define the following \\
\textbf{Problem 8}\\
\emph{
\textbf{FOM optimisation problem:}
\begin{equation}
\label{eqn:optfom}
\min_{\mathcal{P}_k\in \mathcal{P}_{\text{bio}}} J(\mathbf{F}_{1k},\mathcal{P}_k), \quad \mathbf{F}_k \; \text{solution of} \, \eqref{eqn:3}.
\end{equation}
}

Here, $\mathcal{P}_{\text{bio}}$ is the set of biological ranges for each parameter in the set $\mathcal{P}$, as given in Table \ref{tab:1}. 
\\
In order to solve \eqref{eqn:optfom}, we will write an iterative algorithm which, given the FOM solution $\mathbf{F}_k$ corresponding to a parameter set $\mathcal{P}_k$ at iteration $k$, computes the associated \textbf{ROM} and \textbf{ROM linearised} solutions from \textbf{Problem 5} and \textbf{Problem 6}, and minimises the functional $J$ at the ROM level through sensitivity analysis, updating the parameter set and initiating a new iteration $k+1$ until convergence. The value of the functional $J$ calculated at the ROM level is
\begin{equation}
\label{eqn:2rom}
J(\mathbf{RN}_{1k},\mathcal{P}_k)=\frac{1}{2||\pi_h H(\phi_{\text{data}}(T))||_h^2}\biggl|\biggl|\pi_h H_{\phi_e}\biggl(\sum_{i=1}^{N_{\text{POD}}}\alpha_{ik}^N\xi_{i}^{\phi}\biggr)- \pi_h H(\phi_{\text{data}}(T))\biggr|\biggr|_h^2+\frac{\eta}{2}\sum_{m=1}^{|{\mathcal{P}}|}\biggl(\frac{\mathcal{P}_{k,m}-\mathcal{P}_{\text{exp,m}}}{\mathcal{P}_{\text{exp,m}}}\biggr)^2.
\end{equation}
In order to calculate a minimum for \eqref{eqn:2rom}, for $\mathcal{P}_k\in \mathcal{P}_{\text{bio}}$ and with $\mathbf{RN}$ the solution of \eqref{eqn:6}, we use sensitivity analysis to define a projected gradient algorithm which updates the parameters set along descent directions of the functional $J$.
We define the weighted gradient directions, for each component $m=1,\dots,9$ in the parameter set $\mathcal{P}_k$,
\begin{align}
\label{eqn:wg}
&\nabla_{\mathcal{P}_k,w} J(\mathbf{RN}_{1k},\mathcal{P}_{k})|_m:= \\
& \notag \biggl(J_{\alpha}(\mathbf{RN}_{1k},\mathcal{P}_{k})^T (\mathbf{RL}_{1k}^N\text{diag}[d\mathcal{P}_0])\biggr)_m^T+\text{diag}[d\mathcal{P}_0]J_{\mathcal{P}_k}(\mathbf{RN}_{1k},\mathcal{P}_{k})|_m=\frac{1}{||\pi_h H(\phi_{\text{data}}(T))||_h^2}\times\\
& \notag \biggl[\bigl\{\pi_h H_{\phi_e}(\Phi \vec{\alpha}_k^N)- \pi_h H(\phi_{\text{data}}(T))\bigr\}^TM_h\biggl(\text{diag}\biggl\{\pi_h\frac{\partial H_{\phi_e}}{\partial \Phi \vec{\alpha}_k^N}(\Phi \vec{\alpha}_k^N)\biggr\}\Phi \biggl(\frac{\partial \vec{\alpha}_k^N}{\partial \mathcal{P}_{k,m}}d\mathcal{P}_{0,m}\biggr)\biggr)+\\
&\notag \bigl\{ \pi_h H_{\phi_e}(\Phi \vec{\alpha}_k^N)-\pi_h H(\phi_{\text{data}}(T))\bigr\}^T
M_h\biggl\{\pi_h\frac{\partial H_{\phi_e}}{\partial c_{ek}}(\Phi \vec{\alpha}_k^N)dc_{ek,0}\biggr\}\biggr|_{m=9}\biggr]+\eta\frac{(\mathcal{P}_{k,m}-\mathcal{P}_{\text{exp},m})}{\mathcal{P}_{\text{exp},m}^2}d\mathcal{P}_{0,m},
\end{align}
where
\[
d\mathcal{P}_0:=10^{-n_w}\mathcal{P}_0, \; n_w\in \mathbb{N}_+,
\]
is the vector of weights, which defines a weighted Euclidean scalar product in the parameter space. 
 $M_h$ denotes the lumped scalar {mass} matrix. We note that the weighted gradient \eqref{eqn:wg} can be obtained in an equivalent way using the standard Euclidean {inner} product in parameter space and solving the linearised problems \eqref{eqn:romlinearized} at $\mathcal{P}_k$ for a definite variation $\mathcal{P}_k+d\mathcal{P}_0$.\\
Then, we define the {projection} function 
\[
\mathcal{P}_k(\lambda):=\max(\mathcal{P}_{\text{bio,min}},\min(\mathcal{P}_{\text{bio,max}},\mathcal{P}_k-\lambda \nabla_{\mathcal{P}_k,w} J(\mathbf{RN}_{1k},\mathcal{P}_{k}))),
\]
which updates the values of the parameters $\mathcal{P}_k$ along the weighted gradient directions with a learning rate $\lambda$ and projects them onto the feasible set.\\
Following \cite{kelley}, we define the\\
\textbf{Problem 9}\\
\emph{
\textbf{Projected Weighted Gradient Algorithm:}
Given \textbf{RN}, \textbf{RL}, $\mathcal{P}_k$, we update
\\ \\
\textbf{a}: Given $\beta<1$ and $\lambda=\beta$, compute $\mathcal{P}_k(\lambda)$.
\\
\textbf{b}: Find the least integer $m$ such that $\lambda_k=\beta^m$ and
\[
J(\mathbf{RN}(\mathcal{P}_k(\lambda_k)),\mathcal{P}_k(\lambda_k))-J(\mathbf{RN}_{1k},\mathcal{P}_k)\leq \frac{-10^{-4}}{\lambda_k}|\mathcal{P}_k(\lambda_k)-\mathcal{P}_k|^2;
\]
\\
\textbf{c}: \begin{equation}
\label{eqn:pgl}
\mathcal{P}_{k+1}=\mathbf{PWG}(\mathbf{RN}_k,\mathbf{RL}_k,\mathcal{P}_k):=\mathcal{P}_k(\lambda_k).
\end{equation}
\\
}

We finally formulate the following \textbf{Optimization algorithm}:
\begin{algorithm}[h!]
\caption{Optimization Algorithm}\label{alg:opt}
\begin{algorithmic}
\Require MRI(t=0), DTI(t=0), MRI(t=T), $\mathcal{P}_0,\mathcal{P}_{\text{bio}},\mathcal{P}_{\text{av}}$;
\State \textbf{Initialisation}(MRI(t=0), DTI(t=0)) (Problem \eqref{eqn:initialization}); 
\State \textbf{Target}(MRI(t=T)) (Problem \eqref{eqn:target});
\For{$k\geq 0$}  
    \State \textbf{Step 1--\textbf{FOM}}: $\mathbf{F}_k$(\textbf{Initialisation},$\mathcal{P}_k$) (Problem \eqref{eqn:fom});
    \State Compute $J(\mathbf{F}_{1k},\mathcal{P}_k);$
    \If {$k\geq 1$ \textbf{and} $J((\mathbf{F}_{1k},\mathcal{P}_k)\geq J((\mathbf{F}_{1k-1},\mathcal{P}_{k-1})$}
    \State {$\mathcal{P}_{\text{opt}}\leftarrow\mathcal{P}_{k-1}$;}
    \State break;
\ElsIf {$k\geq 1$ \textbf{and} $|J(\mathbf{F}_{1k},\mathcal{P}_k)-J(\mathbf{F}_{1k-1},\mathcal{P}_{k-1})|\leq \text{tol}_F|J(\mathbf{F}_{11},\mathcal{P}_1)-J(\mathbf{F}_{10},\mathcal{P}_0)|$}
    \State {$\mathcal{P}_{\text{opt}}\leftarrow\mathcal{P}_{k}$;}
    \State break;
\EndIf 
\State \textbf{Step 2--\textbf{POD}:} $\mathbf{P}_k(\mathbf{F}_k)$ (problem \eqref{eqn:pod});
\State \textbf{Step 3--Assemble the \textbf{ROM} systems:} $\mathbf{A}_k(\mathbf{P}_k)$ (problem \eqref{eqn:assemble});
 \State\textbf{Step 4--ROM Optimization:} 
\For{$l\geq 0$}  
\State {$\mathcal{P}_{l}\leftarrow\mathcal{P}_{k}$},
\State {\textbf{Step A}: $\mathbf{RN}_l(\mathbf{A}_l,\mathcal{P}_l,\phi_h^0,n_h^0)$; $\mathbf{RL}_l(\mathbf{A}_l,\mathcal{P}_l,\mathbf{RN}_l)$} (problems \eqref{eqn:romnewton},\eqref{eqn:romlinearized}); 
\State {\textbf{Step B} Compute $J(\mathbf{RN}_{1l},\mathcal{P}_l)$};
\State \textbf{Step C:} $\mathcal{P}_{l+1} = \mathbf{PWG}(\mathbf{RN}_l,\mathbf{RL}_l,\mathcal{P}_l)$ (problem \eqref{eqn:pgl});
\If {$\max_{i=1,\dots,|\mathcal{P}|} \biggl((\mathcal{P}_{i,l+1}-\mathcal{P}_{i,l})/\mathcal{P}_{i,l}\biggr)\leq \text{tol}_{Ra}$ \textbf{and} \\ $|J(\mathbf{RN}_{1l+1},\mathcal{P}_{l+1})-J(\mathbf{RN}_{1l},\mathcal{P}_{l})|\leq \text{tol}_{Rb}|J(\mathbf{RN}_{11},\mathcal{P}_1)-J(\mathbf{RN}_{10},\mathcal{P}_0)|$\textbf{and} \\
$|\mathcal{P}_{l+1}(1)-\mathcal{P}_{l+1}|\leq \text{tol}_{Pa}|\mathcal{P}_0|+\text{tol}_{Pr}|\mathcal{P}_0(1)-\mathcal{P}_0|$}
    \State {$\mathcal{P}_{k+1}\leftarrow\mathcal{P}_{l+1};$}
    \State break.
\EndIf 
\EndFor  
\EndFor
\end{algorithmic}
\end{algorithm}
\newpage 
The Algorithm \ref{alg:opt} stops at an iteration $k$ when the functional $J$, calculated at the \textbf{FOM} level, decreases, with respect to its value at $k-1$, by an amount which is a sufficiently small fraction of the initial decrease between iteration $1$ and $0$, which means that a local minimum is being approached. In the same way, the projected gradient iterations stop at an iteration $l$ when the following termination criteria are simultaneously satisfied: the functional $J$, calculated at the \textbf{ROM} level, decreases of a sufficiently small amount with respect to the initial decrease attained at the first iteration, the parameters are changing by a relative small amount along the descent directions and the Euclidean norm of $\mathcal{P}_l(1)-\mathcal{P}_l$, which is a measure of stationarity related to the magnitude of the weighted gradient of $J$, is sufficiently small with respect to the norm of $\mathcal{P}_0$ and the initial norm of $\mathcal{P}_0(1)-\mathcal{P}_0$. After \textbf{Step 4}, we go back to \textbf{Step 1} and calculate the new POD basis associated to the FOM solutions obtained with the new set of parameters $\mathcal{P}_{l+1}$. We thus dynamically span the space of parameters by solving the optimization algorithm over ROM systems associated to different POD bases for each set $\mathcal{P}_k$ at each step $k$. This is indeed an alternative way to consider parameters variability in the MOR of evolution equations with respect to the local reduced basis method used e.g. in \cite{manzoni}, which would request the static computation by k-means clustering of different local POD basis from FOM solutions performed for different sets of parameters, performing the minimisation problem on a ROM level by choosing properly the local basis along the flow of projected gradient parameter updates.\\
\\
We finally observe that the weights in \eqref{eqn:wg} are needed to precondition the ill--conditioned gradient projected algorithm \eqref{eqn:pgl}. Indeed, due to the large differences in the order of magnitude of the parameters in the set $\mathcal{P}$ (see Table \ref{tab:1}), the solutions of the linearised systems \eqref{eqn:romlinearized} differ accordingly (higher linearised solutions for smaller parameters) by order of magnitudes and some parameters may undergo large variations during the first step of the Armijo procedure in \textbf{Problem 9}. These variations could not be fully represented by the information contained in the starting POD basis, which is based on FOM solutions which satisfy the physical constraints, thus causing the Newton algorithm to compute $\mathbf{RN}(\mathcal{P}_k(\beta)))$ to diverge and the solution to violate the physical constraints. The latter fact happens when the negative values associated to the higher order POD basis are amplified during the algorithm.
An alternative way to proceed would be to regularise the ROM system in order to eliminate its instability when the ROM solutions take values in an unphysical range and to penalise unphysical solutions, by choosing a convex potential containing a  smooth penalisation of negative values. This method was found to be unnecessary when weights are introduced in \eqref{eqn:wg}.
\section{Results}
In this section we apply the proposed algorithm to optimise the parameter estimation from the neuroimaging data of  two test cases provided by a clinical study conducted  at the Istituto Neurologico Besta in Milan. Test case $1$ is a clinical follow-up of a primary tumour subjected to adjuvant therapy, optimising the model parameters by following-up the tumour growth until surgical removal. Test case $2$ concerns the recurrence pattern of a GBM after surgery until the start of radiotherapy. While in the first case the tumour mass keeps a round shape, in the  latter it  grows in an irregular manner infiltrating  the peritumoral brain tissue after surgery. 

Since the biological range of the model parameters  to be estimated can vary  as summarised in Table \ref{tab:1}, we use as our initial guess for the optimisation algorithm the manually tuned parameters  in \cite{Agosti1,Agosti2}, reading
\begin{align}
\label{eqn:5}
 \notag \mathcal{P}_0\equiv &\{(1/5000),(0.08),(2.0),(10^4),(8640),(0.1225),(694),(0.3),(0.611)\},\\
 \notag \mathcal{P}_{\text{bio}}\equiv &\{[1/5032.2,1/1377.86],[0.012,0.5],[0.007,90.72],[10^3,10^5], [10^3,10^5],\\
& \notag[0.0841,0.6084],[106.66,1533.3],[0.1,0.33], [0.2,0.611]\},\\
 \mathcal{P}_{\text{exp}}\equiv &\{(1/3991.06),(0.06),(2.0),(10^4),(8640),(0.1225),(694),(0.3),(0.611)\},
\end{align}
with units $mm^2/(\text{Pa}\,\text{day})$, $\text{day}^{-1}$, $mm^2/\text{day}$, $\text{day}^{-1}$, $\text{day}^{-1}$, $\text{Pa}\,mm^2$ and $\text{Pa}$ for the first seven parameters respectively, whereas $\delta$ and $c_e$ are dimensionless.
Moreover, we set $\eta=10^{-4}$, $\text{tol}_F=\text{tol}_{Ra}=\text{tol}_{Rb}=\text{tol}_{Pr}=10^{-3}$, $\text{tol}_{Pa}=10^{-6}$, $n_w=1$.

\subsection{Test case $1$: clinical follow-up of a primary tumour}
\label{sec:2}
\begin{figure}[!ht]

\includegraphics[width=0.7\linewidth]{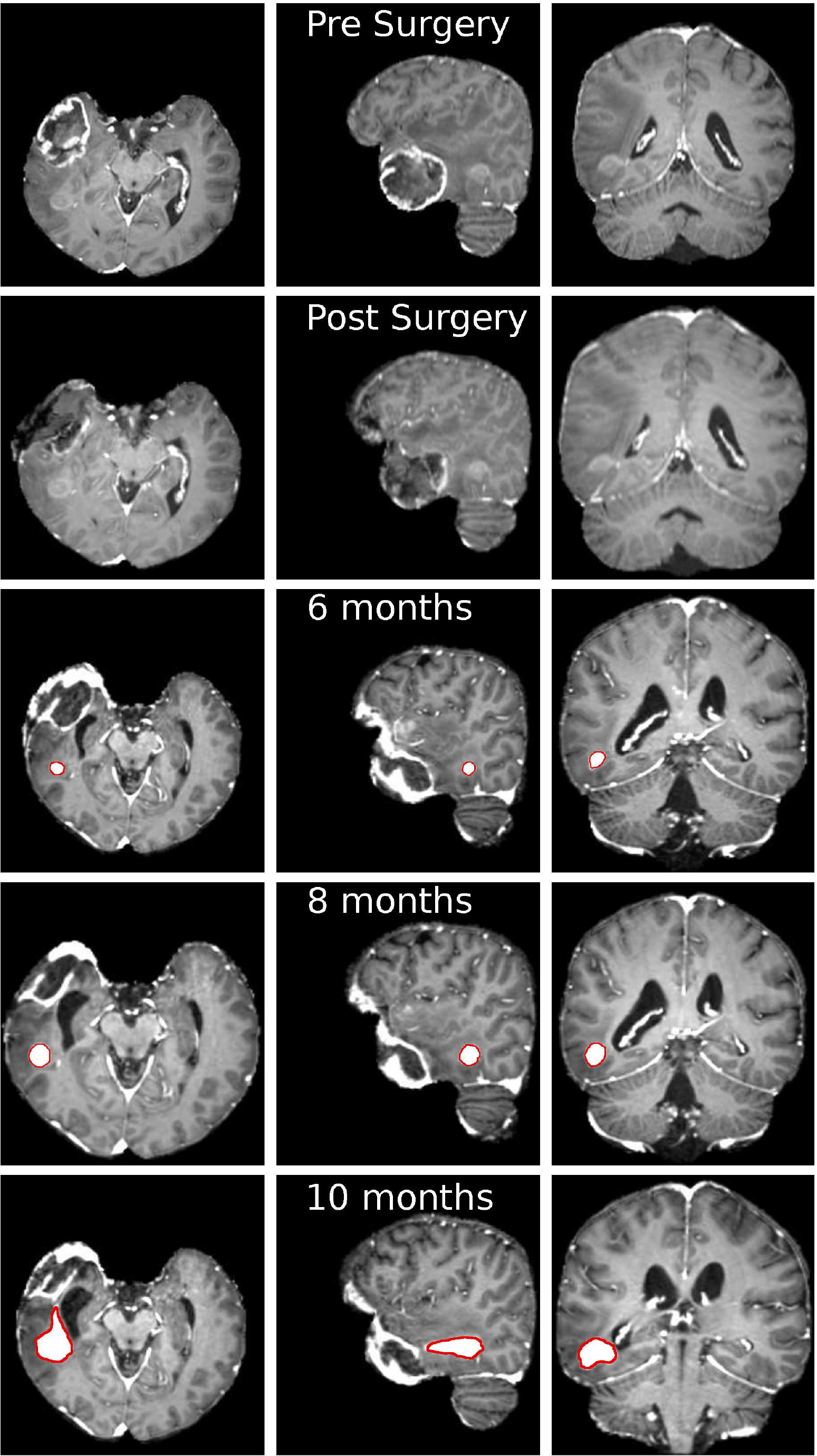}
\centering
\caption{ Axial (first column), sagittal (second column) and coronal (third column) slices of the T1-weighted MRI at different temporal stages. First row: before surgery; second raw: after surgery; third row: 6 months after surgery; fourth row: 8 months after surgery; fifth row: 10 months after surgery. It is possible to appreciate the progressive volumetric increase of the posterior temporal mass. The segmented boundary of the tumour is highlighted in red color.}
\label{fig:1}
\end{figure}

We first apply the proposed Optimization Algorithm for a test case which investigates the clinical follow-up of a primary tumour.\\
A patient diagnosed with multiple GBM lesions underwent  a surgical removal, which left one posterior temporal mass untouched. The patient started radiotherapy with concomitant chemotherapy (Temozolomide)  25 days after surgery following the Stupp protocol; the pre-Radiotherapy MRI confirmed the presence of the primary GBM mass in the posterior temporal area. The patient completed the standard radiation protocol, and performed post-radiotherapy MRI immediately after the treatment and then every two months. These MRI scans depict the progression of the posterior temporal lesion.
Further 6 MRI scans were taken on following the schedule of the clinical protocol until a post-radiotherapy stage at 8 months after surgery.
Our numerical simulations investigate  the follow-up of the growing posterior lesion from $6$ months after surgery (initial time $t=0$ days in the simulations)  to $8$ months after surgery (final time $t=T:=2$ months). At $t=T$ we compare data and simulations, searching for the optimal set of parameters $\mathcal{P}_{\text{opt}}$ which locally minimises the functional \eqref{eqn:2}, obtained by solving \textbf{Algorithm 1}. Between $t=0$ and $t=2$ months the patient underwent two chemotherapy cycles, so that
\[
K_T(t)=\begin{cases}
K_{C_3} \quad 0\leq t \leq 8, \quad 33\leq t \leq 38 \; \text{[days]},\\
0 \quad \text{otherwise}.
\end{cases}
\] 
In Figure \ref{fig:1} we show the axial, sagittal and coronal slices of the T1-weighted MRI at different temporal stages.
The clinicians observed that at $t=4$ months ($10$ months after surgery), the posterior lesion joins  the infiltrated mass re--grown from the peritumoral area affected by  surgery. Therefore, we choose to set the optimization problem only between $t=0$ and $t=2$ months, when the posterior lesion can be considered as a single  tumour mass expanding in a healthy brain tissue.
\subsubsection{Initialisation}
\begin{figure}[!ht]
\centering
\includegraphics[width=0.7\linewidth]{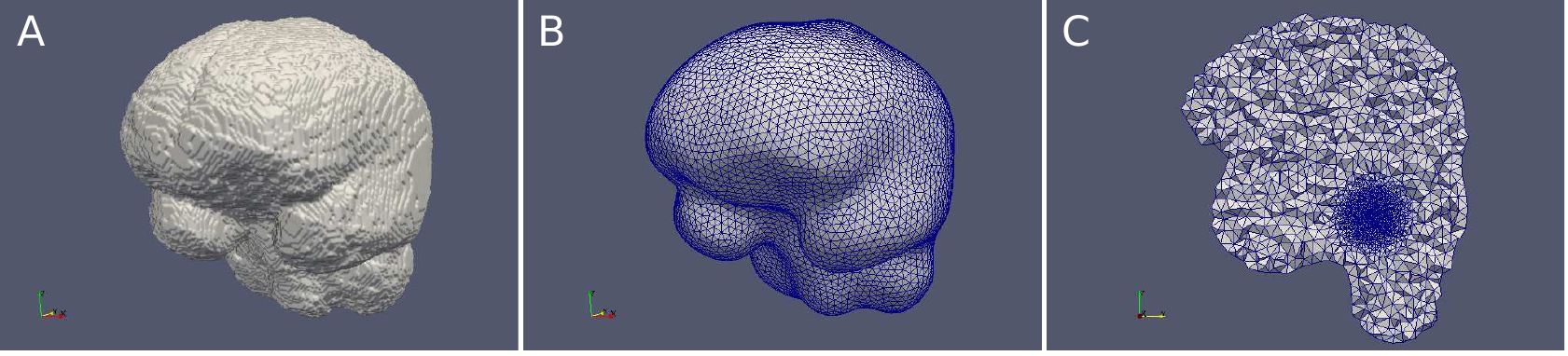}
\caption{ (A) External brain surface extracted from the medical images; (B) Smoothed and re--meshed external surface; (C) Tetrahedral mesh generated within
the external surface, conveniently refined in the peritumoral area.}
\label{fig:2}
\end{figure}

\begin{figure}[!ht]
\includegraphics[width=0.7 \linewidth]{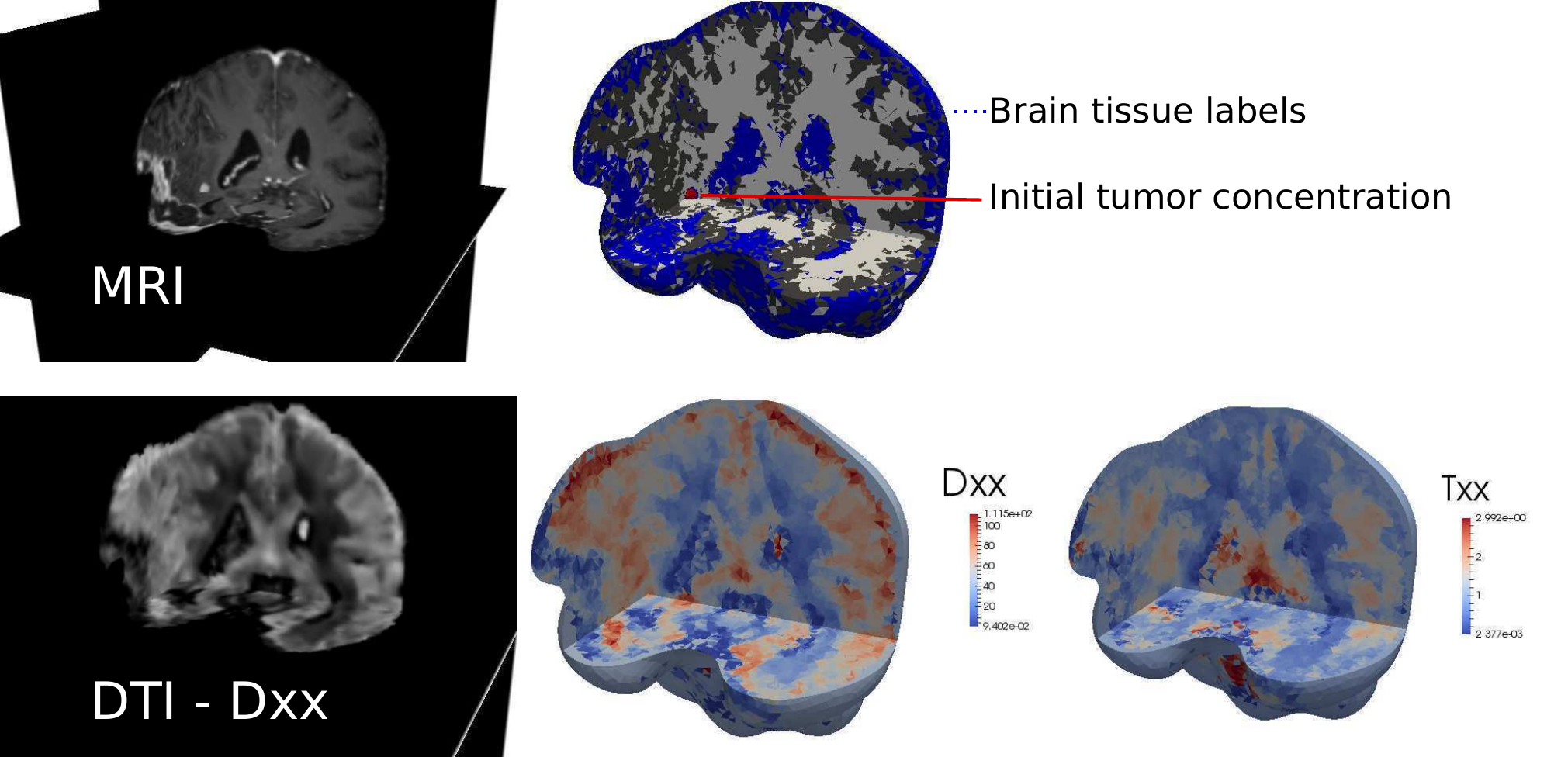}
\centering
\caption{ 3D virtual reconstructions of the MRI (top) and DTI (bottom) data, with the corresponding
computational meshes containing the labels of the brain tissues and the $xx$ component of the tensors \textbf{D} and \textbf{T}. White matter, grey matter and CSF are highlighted in white, grey
and blue colors respectively. The initial tumour distribution is also highlighted in brown color, as segmented from the T1-weighted MRI.
}
\label{fig:3}
\end{figure}

In Figures \ref{fig:2} and \ref{fig:3} we represent the results of the \textbf{initialisation} step of \textbf{Algorithm 1}, which defines the domain $\Omega$ (Figure \ref{fig:2}), the map(WM,GM,CSF), the initial condition $\phi_h^0$ and the tensors \textbf{D} and \textbf{T} (Figure \ref{fig:3}), extracted from the segmentation of the MR images and from the study of the log-signal associated to the DTI  at $t=0$. 

The number of elements and the number of nodes of the mesh $\mathcal{T}_h$ are $301867$ and $50713$ respectively. Moreover, we choose $\Delta T=0.1225$ (days), and $N=490$.\\
The characteristic function of the tumour extension $\phi_{\text{data}}(T)$ has been extracted  from the segmentation of the MR images at $t=2$ month.
\subsubsection{Step 1}
In Figure \ref{fig:4} we report the values of the functional $J(\phi_h^N(\mathcal{P}_k),\mathcal{P}_k)$, calculated in \textbf{step 1} of \textbf{Algorithm 1}, and of the set of parameters $\mathcal{P}_k$, for different values of $k$. We also plot the iso--surfaces $\phi_{\text{data}}(T)=0$ from the MRI data and $\phi_h^N(\mathcal{P}_k)=\phi_e/2$ from the FOM simulations, reporting the value of the Jaccard index, defined as the intersection over union ratio between the two volumes  enclosed within these two surfaces.
{\footnotesize
\centering
\begin{table}[h!]
\resizebox{0.9\textwidth}{!}{\begin{minipage}{\textwidth}
\begin{tabular}{ c c c c c c c c c c c }
 \hline
 \parbox[t][][t]{2.2cm}{\centering \textbf{Iteration}\\\textbf{k=0}} & \parbox[t][][t]{1.6cm}{\centering $\mathbf{J}(\mathcal{P}_0)$\\0.32328} &\parbox[t][][t]{0.9cm}{\centering $\mathbf{L}_0$\\0.0002} & \parbox[t][][t]{0.9cm}{\centering $\boldsymbol{\nu}_0$\\0.08} &  \parbox[t][][t]{0.9cm}{\centering $\mathbf{k}_{n0}$\\2} & \parbox[t][][t]{0.9cm}{\centering $\mathbf{S}_{n0}$\\10000}  & \parbox[t][][t]{0.9cm}{\centering $\boldsymbol{\delta}_{n0}$\\8640} & \parbox[t][][t]{0.9cm}{\centering $\boldsymbol{\gamma}_0^2$\\0.1225} &  \parbox[t][][t]{0.9cm}{\centering $\mathbf{E}_{0}$\\694} &  \parbox[t][][t]{0.9cm}{\centering $\boldsymbol{\delta}_{0}$\\0.3} &  \parbox[t][][t]{0.9cm}{\centering $\textbf{c}_{e0}$\\0.611}\\  \hline \\
 \textbf{MRI} &  &  & &  \textbf{FOM} &   &  & \textbf{Comparison} &  & & \\ \\
\end{tabular}
\begin{tabular}{l}
\raisebox{-0.5\height}{\includegraphics[scale=0.92]{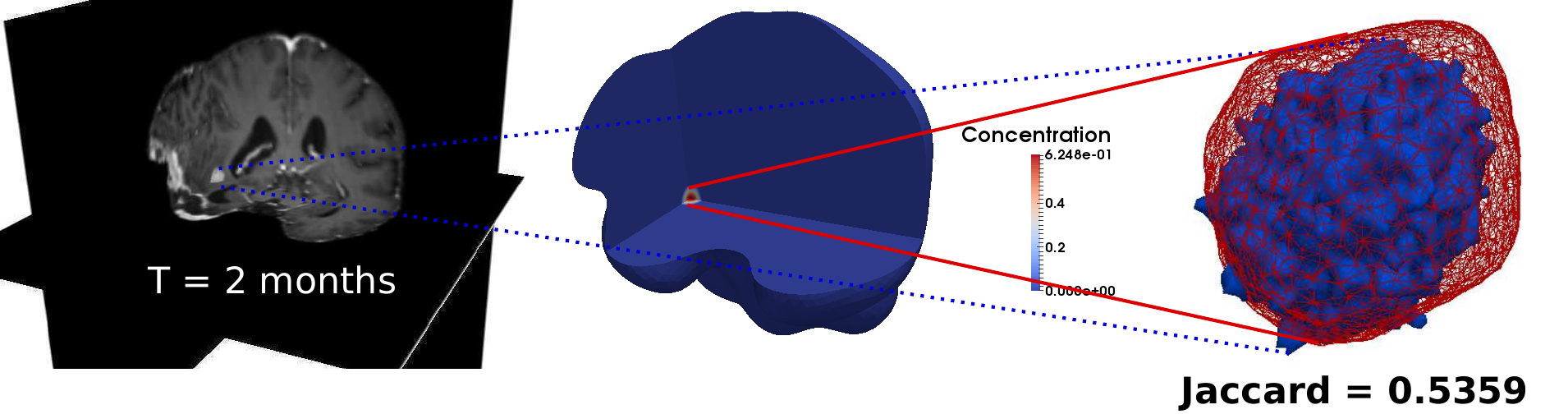}}  \\
  \hline
  \end{tabular}
  \begin{tabular}{ c c c c c c c c c c c   }
 \parbox[t][][t]{2.2cm}{\centering \textbf{Iteration}\\\textbf{k=1}} & \parbox[t][][t]{1.2cm}{\centering $\mathbf{J}(\mathcal{P}_1)$\\0.08001} &\parbox[t][][t]{1.2cm}{\centering $\mathbf{L}_1$\\0.0002} & \parbox[t][][t]{1.cm}{\centering $\boldsymbol{\nu}_1$\\0.02352} &  \parbox[t][][t]{1.cm}{\centering $\mathbf{k}_{n1}$\\1.9769} & \parbox[t][][t]{1.cm}{\centering $\mathbf{S}_{n1}$\\9999.99}  & \parbox[t][][t]{1.cm}{\centering $\boldsymbol{\delta}_{n1}$\\8640.00} & \parbox[t][][t]{0.9cm}{\centering $\boldsymbol{\gamma}_1^2$\\0.1225} &  \parbox[t][][t]{0.9cm}{\centering $\mathbf{E}_{1}$\\693.99}&  \parbox[t][][t]{0.9cm}{\centering $\boldsymbol{\delta}_{0}$\\0.3237} &  \parbox[t][][t]{0.9cm}{\centering $\textbf{c}_{e0}$\\0.5755} \\  \hline \\
 \textbf{MRI} &  &  & & \textbf{FOM}    &  &  & \textbf{Comparison} &  & &   \\ \\
\end{tabular}
\begin{tabular}{l}
\raisebox{-0.5\height}{\includegraphics[scale=0.92]{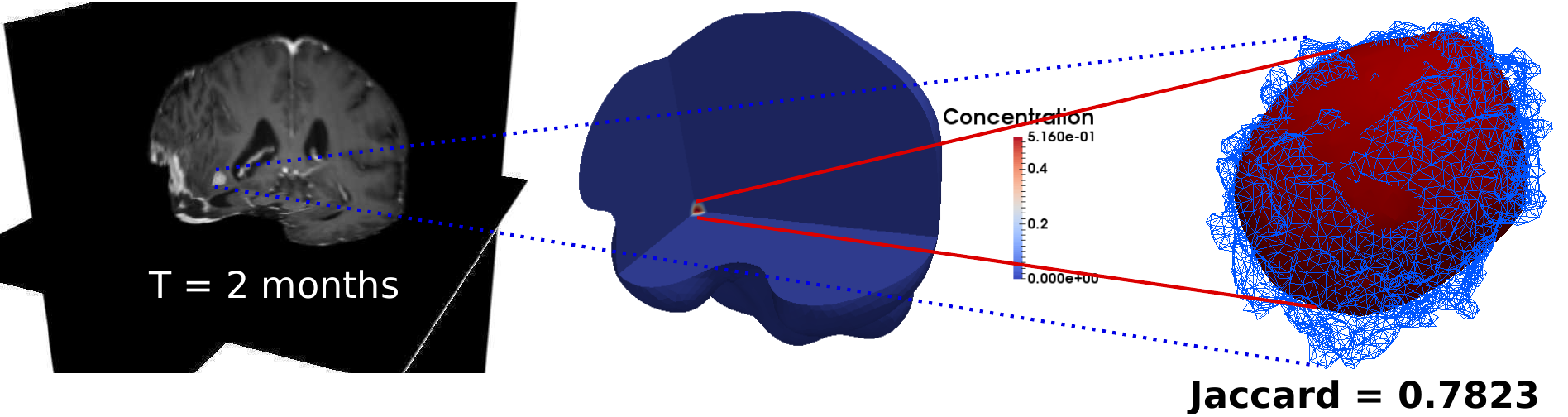}}  \\
  \hline
  \end{tabular}
  \begin{tabular}{  c c c c c c c c c  c c  }
 \parbox[t][][t]{2.2cm}{\centering \textbf{Iteration}\\\textbf{k=2}} & \parbox[t][][t]{1.2cm}{\centering $\mathbf{J}(\mathcal{P}_2)$\\0.07167} &\parbox[t][][t]{1.2cm}{\centering $\mathbf{L}_2$\\0.0002} & \parbox[t][][t]{1.cm}{\centering $\boldsymbol{\nu}_2$\\0.0213} &  \parbox[t][][t]{1.cm}{\centering $\mathbf{k}_{n2}$\\1.9842} & \parbox[t][][t]{1.cm}{\centering $\mathbf{S}_{n2}$\\10000.04}  & \parbox[t][][t]{1.cm}{\centering $\boldsymbol{\delta}_{n2}$\\8639.95} & \parbox[t][][t]{0.9cm}{\centering $\boldsymbol{\gamma}_2^2$\\0.1225} &  \parbox[t][][t]{0.9cm}{\centering $\mathbf{E}_{2}$\\693.99} &  \parbox[t][][t]{0.9cm}{\centering $\boldsymbol{\delta}_{0}$\\0.3271} &  \parbox[t][][t]{0.9cm}{\centering $\textbf{c}_{e0}$\\0.611}\\  \hline \\
 \textbf{MRI} &  &  & &\textbf{FOM}    &  &  & \textbf{Comparison} &  & & \\ \\
\end{tabular}
\begin{tabular}{l}
\raisebox{-0.5\height}{\includegraphics[scale=0.92]{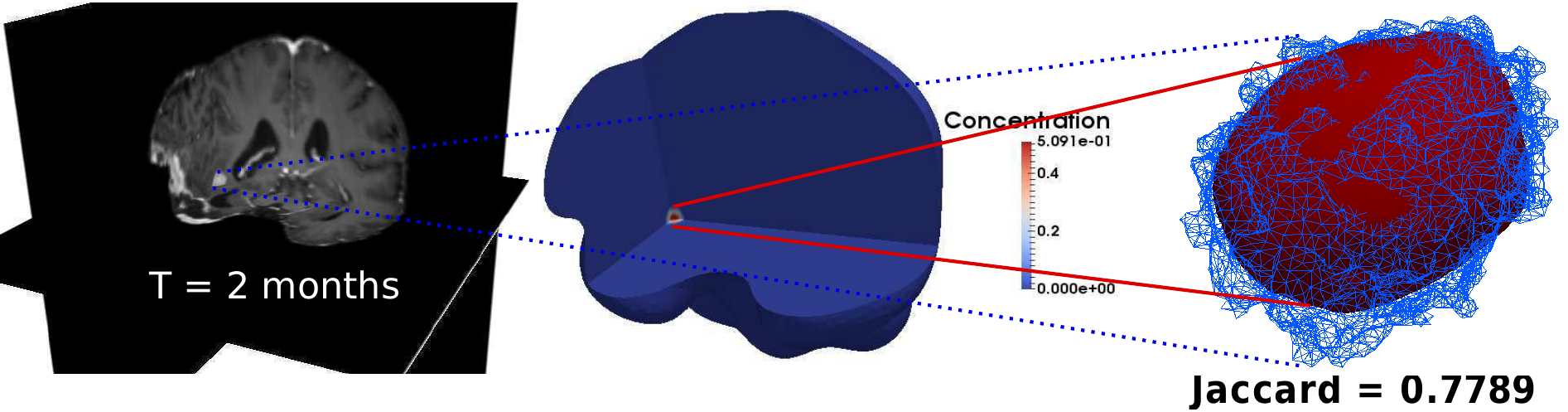}}  \\
  \hline
  \end{tabular}
  \end{minipage} }
  \captionof{figure}[]{Values of $J(\phi_h^N(\mathcal{P}_k),\mathcal{P}_k)$ and of $\mathcal{P}_k$ for different iteration steps $k$ of \textbf{Algorithm 1}, with a comparison between the isosurfaces $\phi_{\text{data}}(T)=0$ (highlighted in blue color) and $\phi_h^N(\mathcal{P}_k)=\phi_e/2$ (highlighted in red color).}
  \label{fig:4}
\end{table}
}
\newpage
\noindent The \textbf{Algorithm 1} stops at $k=2$, since
\[J(\phi_h^N(\mathcal{P}_{k+1}),\mathcal{P}_{k+1})=J(\phi_h^N(\mathcal{P}_{k}),\mathcal{P}_{k})\] for $k=2$. Indeed, the ROM optimization algorithm for $k=2$ makes no advances in the parameter space (see Figure \ref{fig:7}). 
We thus identify 
\begin{align}
\label{eqn:opttest1}
\notag &\mathcal{P}_{\text{opt}}\equiv \mathcal{P}_2=\\
&\notag \{L=0.0002,\nu=0.0213,k_n=1.9842,S_n=10000.04,\delta_n=8639.95,\gamma^2=0.1225,\\
& E=693.99,\delta=0.3271,c_e=0.611\},
\end{align} 
with the corresponding units.
We observe that during the optimization process at 
different $k$ levels the overall overlapping between the tumour extensions from FOM simulations and from data is increasing, which is shown by the corresponding increase in the value of the Jaccard index.
\subsubsection{Step 2}
In Table \ref{tab:2} we report, for each step $k$ of \textbf{Algorithm 1}, the values of the cumulated fractions of $\text{tr}\mathbf{F}_1^T\mathbf{F}_1,\text{tr}\mathbf{F}_2^T\mathbf{F}_2,\text{tr}\mathbf{F}_3^T\mathbf{F}_3,\text{tr}(\psi_1'(\mathbf{F}_1))^T\psi_1'(\mathbf{F}_1),\text{tr}(\psi_1''(\mathbf{F}_1))^T\psi_1''(\mathbf{F}_1)$ associated to the eigenvalues of the corresponding matrices, arranging them starting from the eigenvalue with the highest magnitude and following a decreasing order.
\begin{table}[h!]
\caption{Values of the cumulated fraction of $\text{tr}\mathbf{F}_1^T\mathbf{F}_1$,$\text{tr}\mathbf{F}_2^T\mathbf{F}_2$,$\text{tr}\mathbf{F}_3^T\mathbf{F}_3$,$\text{tr}(\psi_1'(\mathbf{F}_1))^T\psi_1'(\mathbf{F}_1)$,$\text{tr}(\psi_1''(\mathbf{F}_1))^T\psi_1''(\mathbf{F}_1)$ for the first eigenvalues with the highest magnitude.}
\label{tab:2}    
\begin{center}
\begin{tabular}{lllllll}
\hline
\parbox[t][][t]{2.2cm}{\centering \textbf{Iteration}\\\textbf{k=0}} & Eigenvalue & \% $\text{tr}\mathbf{F}_1^T\mathbf{F}_1$ & \% $\text{tr}\mathbf{F}_2^T\mathbf{F}_2$ & \% $\text{tr}\mathbf{F}_3^T\mathbf{F}_3$ & \% $\text{tr}(\psi_1')^T\psi_1'$ & \% $\text{tr}(\psi_1'')^T\psi_1''$\\
\hline
& First & $95.1394$ & $99.9195$ & $99.9695$ & $99.9098$ & $99.3874$ \\ 
& Second & $99.4332$ & $99.9868$ & $99.9964$ & $99.9778$ & $99.8069$\\
& Third & $99.8773$ & $99.9970$ & $99.9991$ & $99.9950$ & $99.9483$\\ 
& Fourth & $99.9695$ & $99.9986$ & $99.9996$ & $99.9983$ & $99.9882$\\ 
& Fifth & $99.9912$ & $99.9997$ & $99.9999$ & $99.9997$ & $99.9980$\\
 \hline
 \parbox[t][][t]{2.2cm}{\centering \textbf{Iteration}\\\textbf{k=1}} & Eigenvalue & \% $\text{tr}\mathbf{F}_1^T\mathbf{F}_1$ & \% $\text{tr}\mathbf{F}_2^T\mathbf{F}_2$ & \% $\text{tr}\mathbf{F}_3^T\mathbf{F}_3$ & \% $\text{tr}(\psi_1')^T\psi_1'$ & \% $\text{tr}(\psi_1'')^T\psi_1''$\\
\hline
& First & $97.6975$ & $99.9483$ & $99.9891$ & $99.9133$ & $99.0374$ \\ 
& Second & $99.6790$ & $99.9950$ & $99.9987$ & $99.9932$ & $99.9278$\\
& Third & $99.9524$ & $99.9989$ & $99.9997$ & $99.9985$ & $99.9849$\\ 
& Fourth & $99.9912$ & $99.9996$ & $99.9998$ & $99.9996$ & $99.9971$\\ 
 \hline
 \parbox[t][][t]{2.2cm}{\centering \textbf{Iteration}\\\textbf{k=2}} & Eigenvalue & \% $\text{tr}\mathbf{F}_1^T\mathbf{F}_1$ & \% $\text{tr}\mathbf{F}_2^T\mathbf{F}_2$ & \% $\text{tr}\mathbf{F}_3^T\mathbf{F}_3$ & \% $\text{tr}(\psi_1')^T\psi_1'$ & \% $\text{tr}(\psi_1'')^T\psi_1''$\\
\hline
& First & $97.7632$ & $99.9500$ & $99.9890$ & $99.9105$ & $99.0029$ \\ 
& Second & $99.6881$ & $99.9953$ & $99.9988$ & $99.9933$ & $99.9280$\\
& Third & $99.9542$ & $99.9990$ & $99.9997$ & $99.9986$ & $99.9855$\\ 
& Fourth & $99.9917$ & $99.9997$ & $99.9999$ & $99.9996$ & $99.9973$\\ 
 \hline
\end{tabular}
\end{center}
\end{table}\\
We thus have that $N_{\text{POD}}=N_{\phi}^{\text{POD}}=5$ for $k=0$ and $N_{\text{POD}}=N_{\phi}^{\text{POD}}=4$ for $k=1,2$. 
In Figure \ref{fig:5} we show the basis elements $\xi_{i}^{\phi}$, corresponding to the highest eigenvalues needed to explain the $99.99\%$ variance of the data, for $k=0,1$, superposed with the initial condition and final distribution of cell concentration (highlighted by a distribution of green and red points respectively).
We observe that $\xi_{1}^{\phi}$ and $\xi_{2}^{\phi}$ are distributed over the bulk of the final state $\phi_h^N$ and the initial condition $\phi_h^0$ respectively, whereas $\xi_{3}^{\phi}$, $\xi_{4}^{\phi}$ and $\xi_{5}^{\phi}$ are oscillating functions over the set where the tumour is expanding during its temporal evolution, and thus contain the information about the tumour boundary and its expansion. \newpage
\noindent We observe that the number of basis functions needed to explain the $99.99\%$ variance of the data is small, which depends on the fact that the region of tumour expansion small and there are no significant topological changes in the evolution dynamics of the FOM solution, which also spreads in a smooth manner. In order to deal with a sufficiently smooth FOM dynamics we needed to choose a mesh $\mathcal{T}_h$ sufficiently refined in the region of tumour evolution. This turned out to be necessary in order to deal with low dimensional higher order tensors in \textbf{Problem 4} \eqref{eqn:assemble} and to deal with ROM systems which are solvable with low computational resources and in highly reduced computational times.\\

\begin{figure}[!h]
\centering 
\includegraphics[width=0.9 \linewidth]{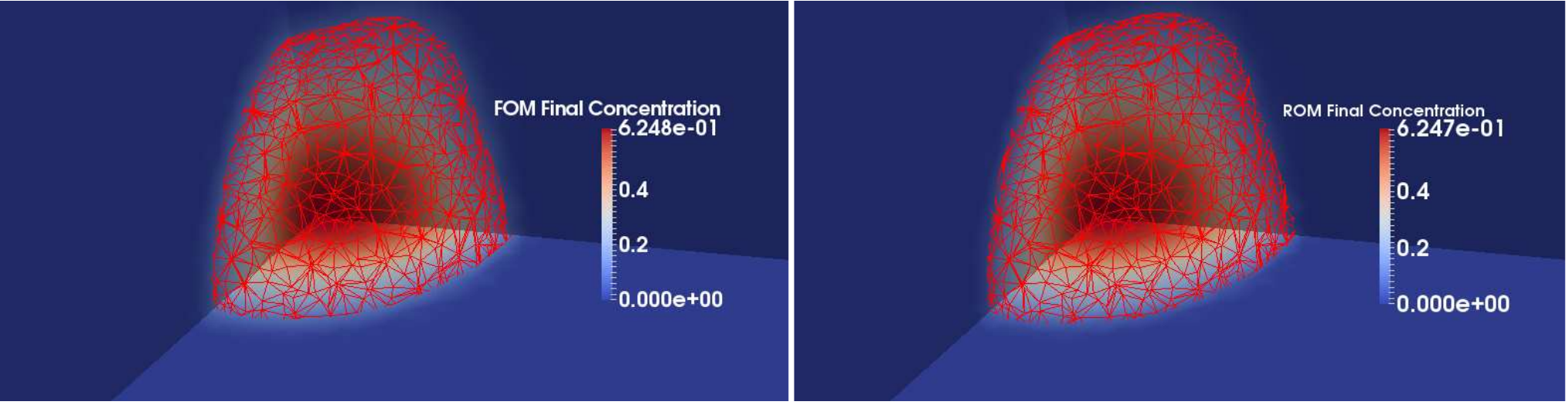}
\caption{Comparison between the final state $\phi_h^N$, solution of the FOM simulation with parameter set $\mathcal{P}_0$ and the corresponding final state $\sum_{i=1}^{N_{\text{POD}}}\alpha_{i0}^N\xi_{i}^{\phi}$, solution of the ROM system \eqref{eqn:6}. The iso--surfaces $\phi_h^N=\phi_e/2$ and $\sum_{i=1}^{N_{\text{POD}}}\alpha_{i0}^N\xi_{i}^{\phi}=\phi_e/2$ are highlighted in red colors.
}
\label{fig:6}
\end{figure}

\begin{table}[h!]
\resizebox{0.9\textwidth}{!}{\begin{minipage}{\textwidth}
\begin{tabular}{ c }
 \hline
 \parbox[t][][t]{17.4cm}{\centering \textbf{Iteration k=0}} \\  \hline \\
\end{tabular}\\
\begin{tabular}{l}
\raisebox{-0.5\height}{\centering \includegraphics[scale=1.3]{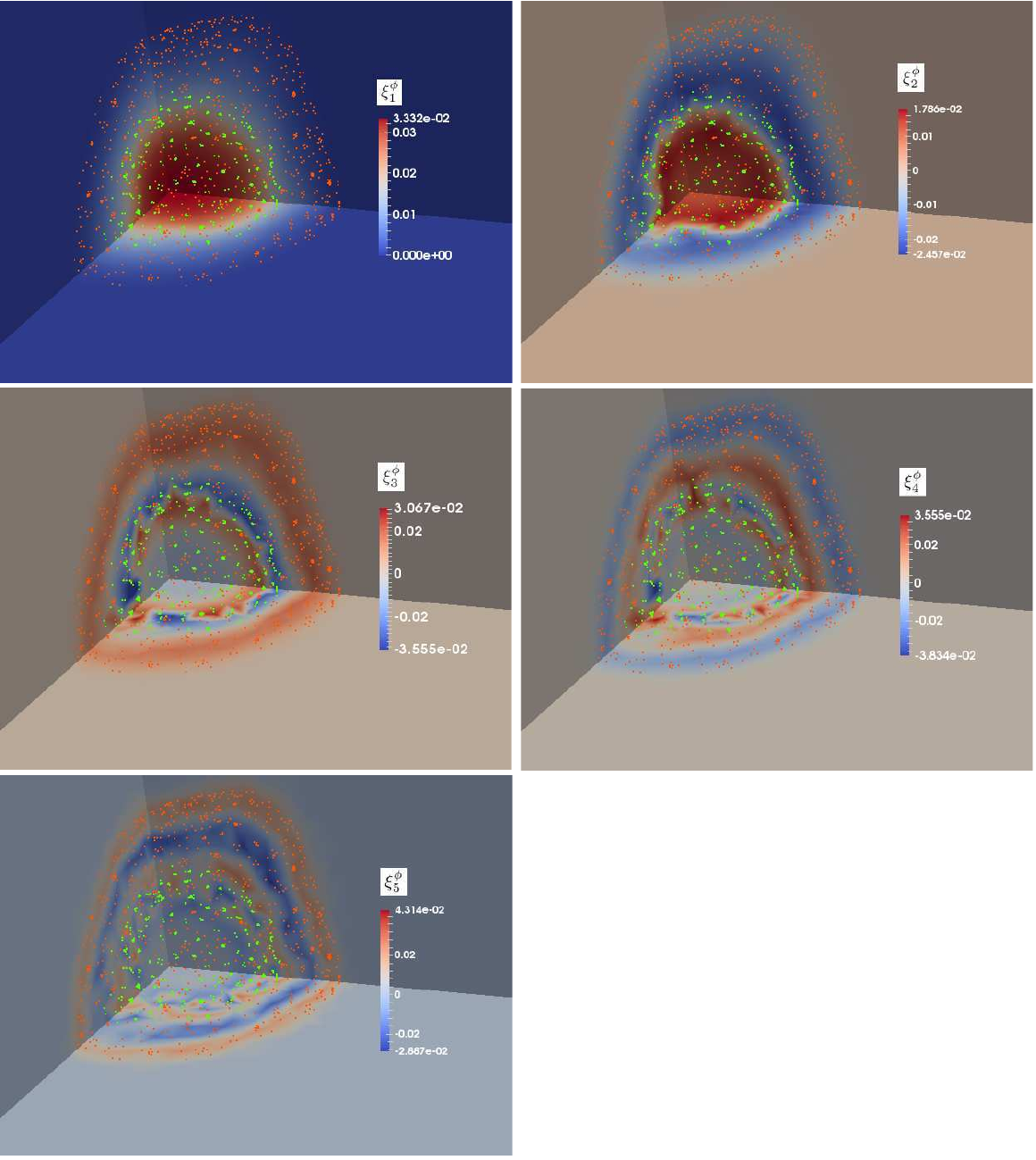}}  \\
  \hline
  \end{tabular}\\
  \end{minipage} }%
  \end{table}%
  \begin{table}[h!]
\resizebox{0.9\textwidth}{!}{\begin{minipage}{\textwidth}
  \begin{tabular}{ c }
  \hline
 \parbox[t][][t]{17.6cm}{\centering \textbf{Iteration k=1}} \\  \hline \\
\end{tabular}\\
\begin{tabular}{l}
\raisebox{-0.5\height}{\centering \includegraphics[scale=1.3]{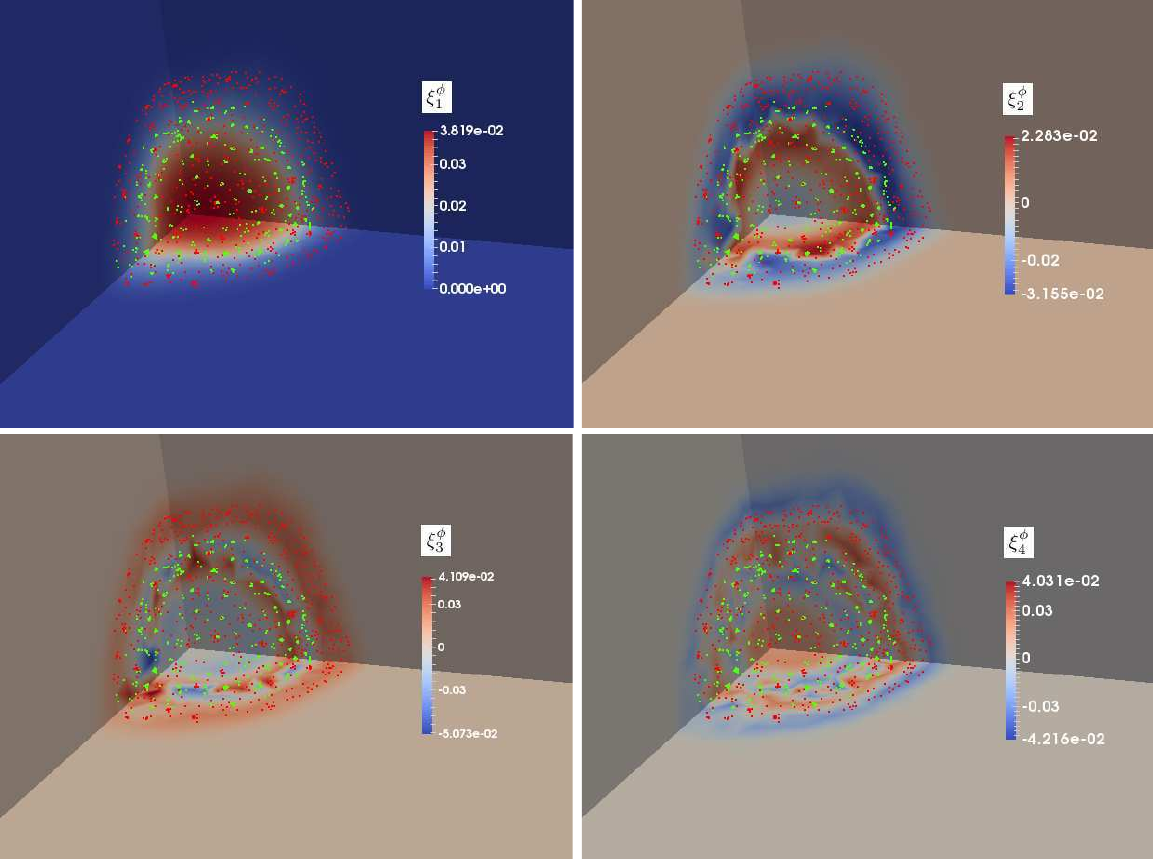}}  \\
  \hline
  \end{tabular}
  \end{minipage} }
  \captionof{figure}[]{Plot of the basis elements $\xi_{i}^{\phi}$ corresponding to the highest eigenvalues needed to explain the $99.99\%$ variance of the data, for $k=0,1$. Green and red points are
distributed over the initial condition and final distribution of cell concentration, respectively.}
  \label{fig:5}
\end{table}
\newpage
We finally show in Figure \ref{fig:6} a comparison between the final state $\phi_h^N$ calculated from the FOM simulation through \textbf{Algorithm 2} with parameter set $\mathcal{P}_0$ and the corresponding final state $\sum_{i=1}^{N_{\text{POD}}}\alpha_{i0}^N\xi_{i}^{\phi}$ obtained as a solution of the ROM system \eqref{eqn:6} through \textbf{Algorithm 4}. 

We observe that the ROM solution is approximating the FOM solution with a very high fidelity.

\subsubsection{Steps 3 and 4}
In Figure \ref{fig:7} we report the values of the functional $J(\vec{\alpha}_l,\mathcal{P}_l)$, of the normalised set of parameters 
\[
\mathcal{P}_l/\mathcal{P}_{\text{exp}}=\{L_l/L_{\text{exp}},\nu_l/\nu_{\text{exp}},k_{nl}/k_{n\text{exp}},S_{nl}/S_{n\text{exp}},\delta_{nl}/\delta_{n\text{exp}},\gamma_l^2/\gamma_{\text{exp}}^2,E_l/E_{\text{exp}},\delta_c/\delta_{c\text{exp}},c_e/c_{e\text{exp}}\},
\]
and of $|\mathcal{P}_l(1)-\mathcal{P}_l|$, computed in \textbf{Steps 3} and \textbf{4}  of \textbf{Algorithm 1}, for $k=0,1,2$. We also plot the iso--surfaces $\phi_{\text{data}}(T)=0$ from the MRI data and $\sum_{i=1}^{N_{\text{POD}}}\alpha_{i\bar{l}}^N\xi_{i}^{\phi}=\phi_e/2$ from the ROM simulations, where $\bar{l}$ is the number of the last iteration of \textbf{Step 4}, reporting the value of the Jaccard index between the two volumes enclosed by these surfaces.\\
For $k=0$ the ROM optimization process in \textbf{Steps 3} and \textbf{4} goes through $11$ steps before matching the termination conditions.
We also observe that only the model parameters $L,\nu,\delta,c_e$ change significantly from their starting values during the optimization process, being the system quite insensitive to changes of the remaining parameters $k_n,S_n,\delta_n,\gamma^2,E$. The proliferation rate $\nu$ is the most sensitive parameter whose variations lead to the functional minimisation. This is in accordance to the sensitivity analysis found in the literature for tumour growth models based on Cahn--Hilliard--Darcy--Forchheimer--Brinkman equations with logistic growth \cite{fritz}.\\
For $k=1$ the ROM optimization process in \textbf{Steps 3} and \textbf{4} goes through $45$ steps. Finally, for $k=2$ the ROM optimization process is making no progress.

\begin{table}[p]
\resizebox{0.9\textwidth}{!}{\begin{minipage}{\textwidth}
\begin{tabular}{ c c c }
 \hline
 \parbox[t][][t]{5.2cm}{\centering \textbf{Iteration k=0}} & &\\ \hline \\ 
 \parbox[t][][t]{5.6cm}{\centering $\mathbf{J}(\vec{\alpha}_l,\mathcal{P}_l)$} &\parbox[t][][t]{6.6cm}{\centering $\mathcal{P}_l/\mathcal{P}_{\text{exp}}$}& \\  
  \raisebox{-0.5\height}{\includegraphics[scale=0.23]{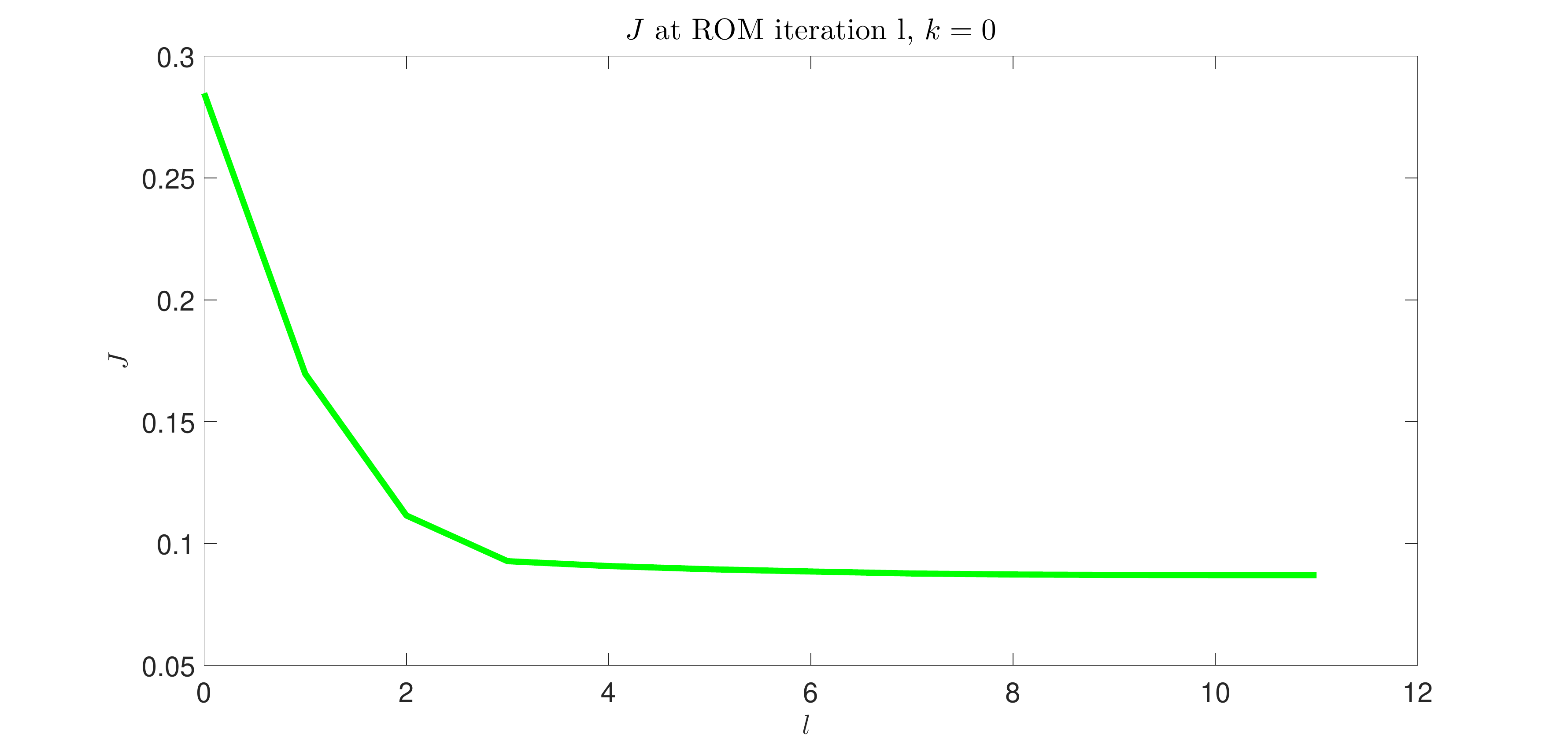}}&\raisebox{-0.5\height}{\includegraphics[scale=0.23]{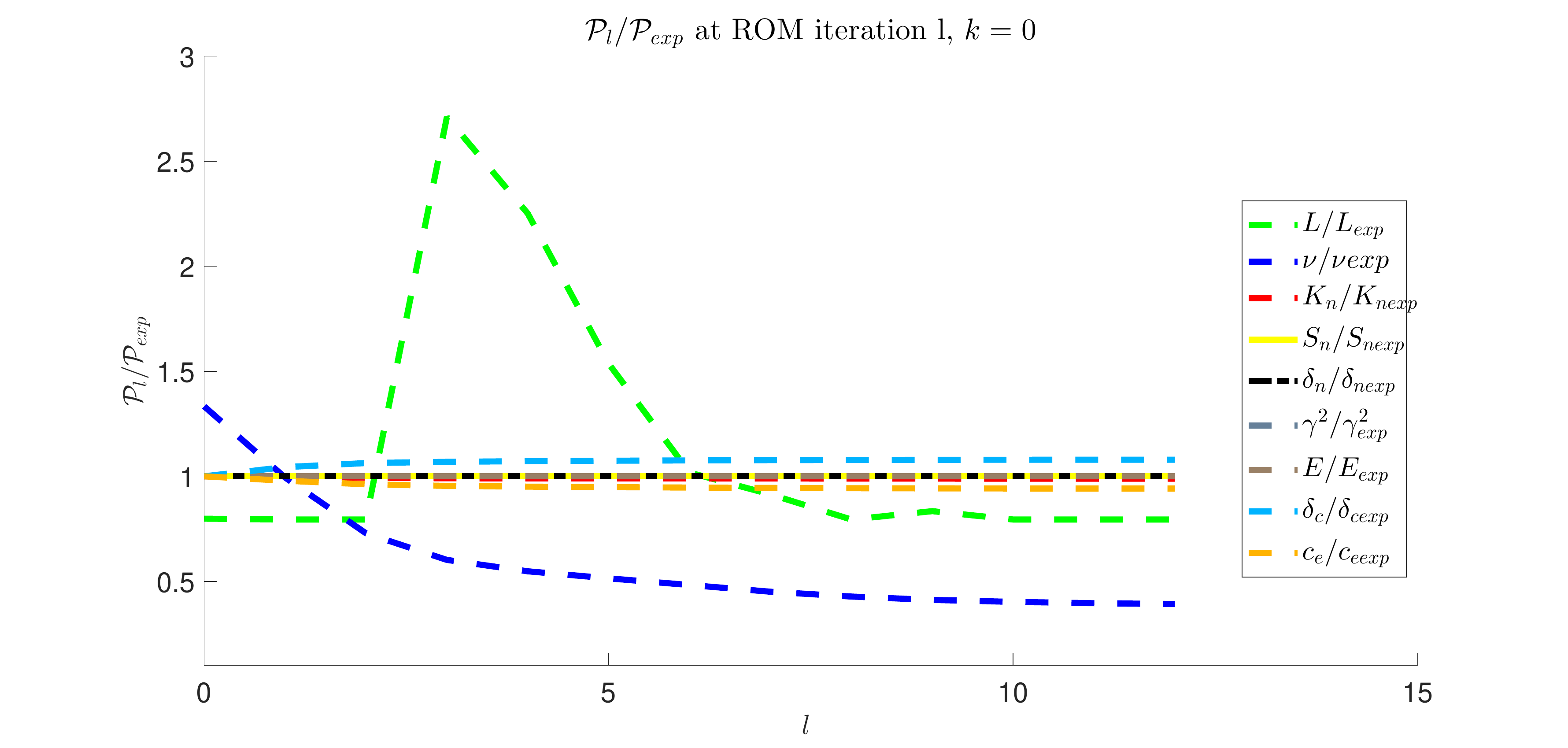}}&\\ \\
\parbox[t][][t]{5.6cm}{\centering $|\mathcal{P}_l(1)-\mathcal{P}_l|$} & & \\  
  \raisebox{-0.5\height}{\includegraphics[scale=0.23]{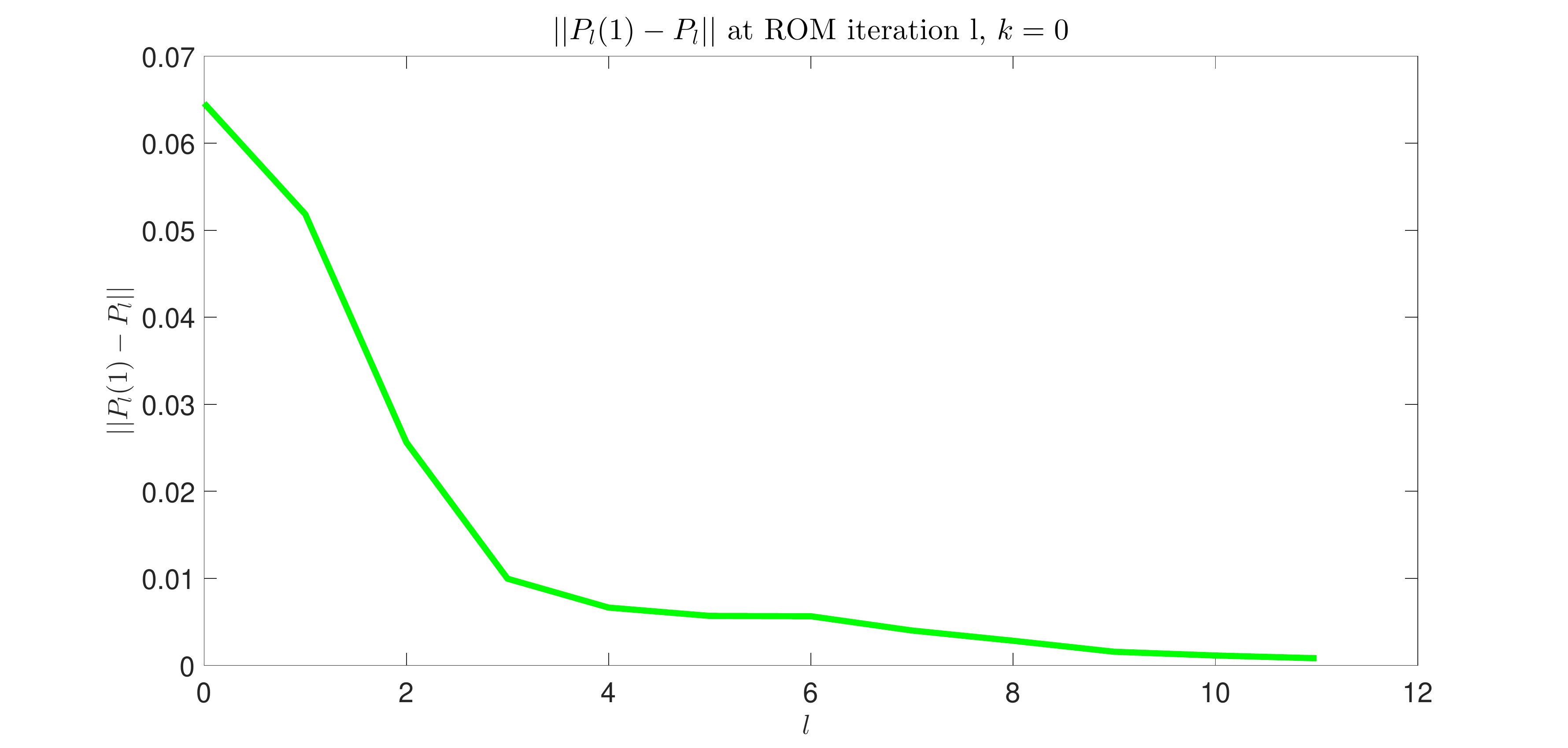}}& &\\ \\  
 \end{tabular}
 \begin{tabular}{c c c}
 \parbox[t][][t]{3.2cm}{\centering\textbf{MRI}} &   \parbox[t][][t]{7.2cm}{\centering\textbf{ROM ($\mathcal{P}_{\bar{l}=11})$}}  &  \parbox[t][][t]{4.2cm}{\centering\textbf{Comparison}}  \\ \\
\end{tabular}\\
\begin{tabular}{l}
\raisebox{-0.5\height}{\includegraphics[scale=0.8]{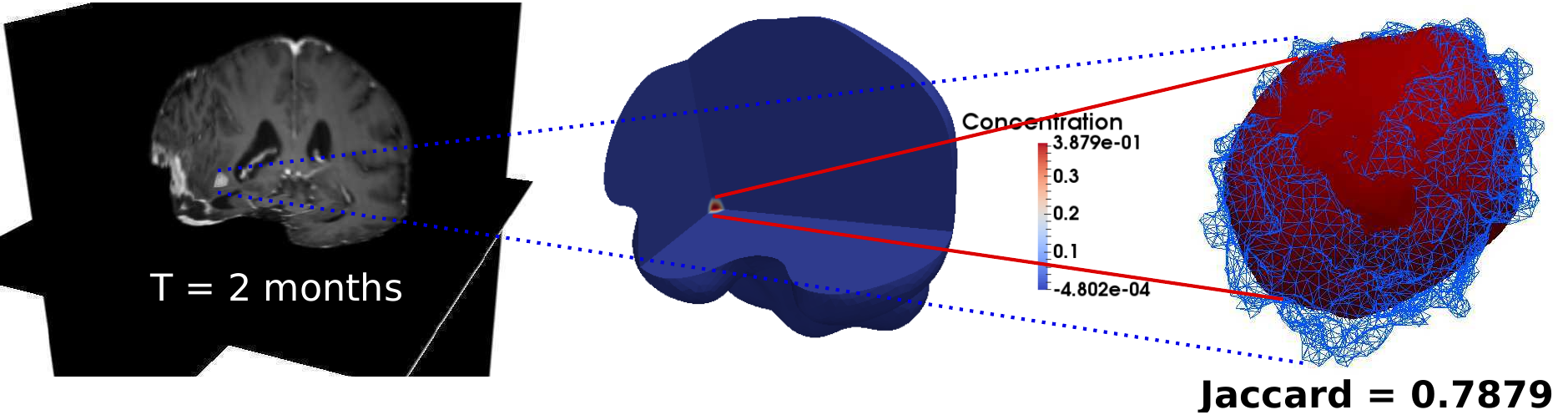}}  \\
  \end{tabular}
  \begin{tabular}{ c c c }
 \hline
 \parbox[t][][t]{5.2cm}{\centering \textbf{Iteration k=1}} & &\\ \hline \\ 
 \parbox[t][][t]{5.6cm}{\centering $\mathbf{J}(\vec{\alpha}_l,\mathcal{P}_l)$} &\parbox[t][][t]{6.6cm}{\centering $\mathcal{P}_l/\mathcal{P}_0$} &\\  \
  \raisebox{-0.5\height}{\includegraphics[scale=0.23]{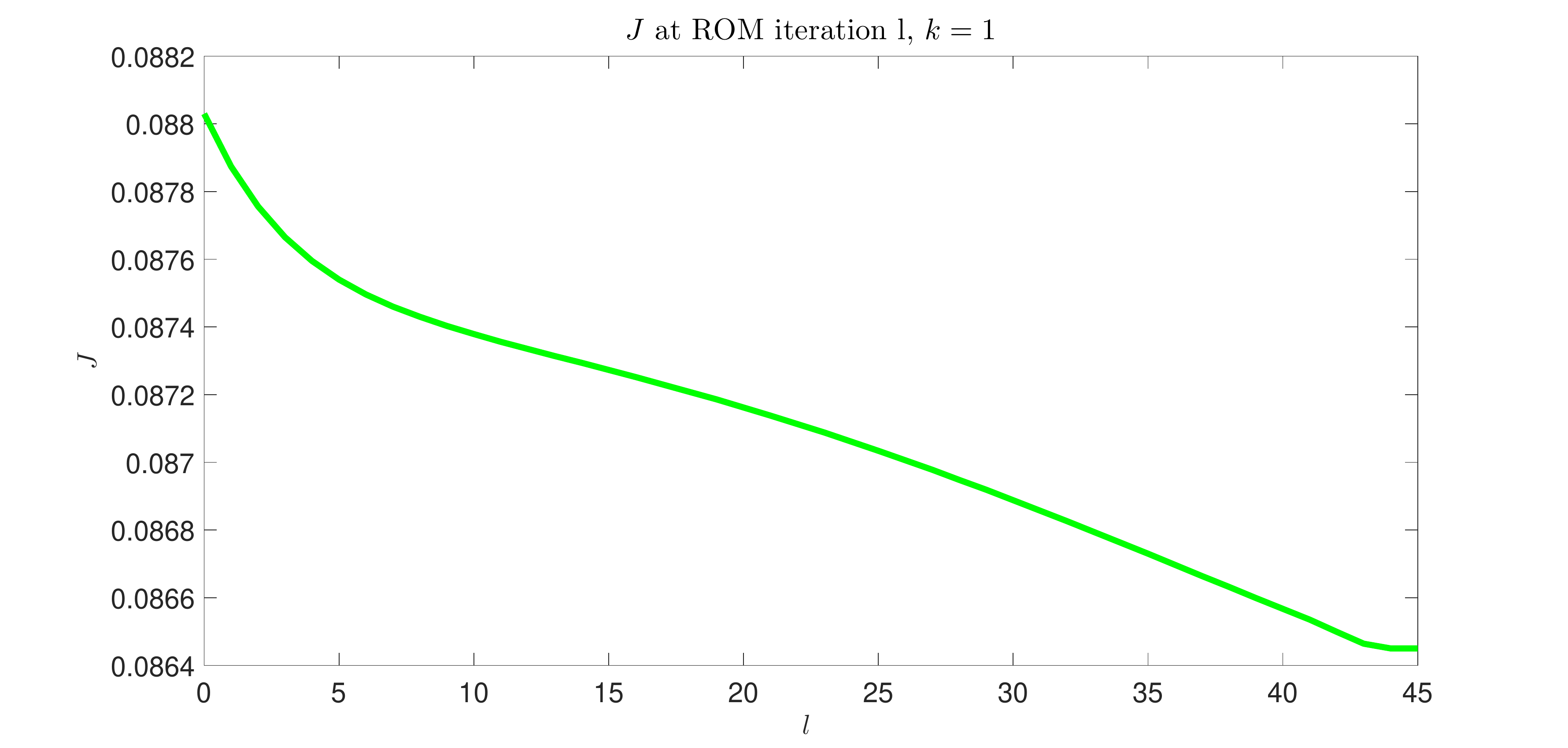}}&\raisebox{-0.5\height}{\includegraphics[scale=0.23]{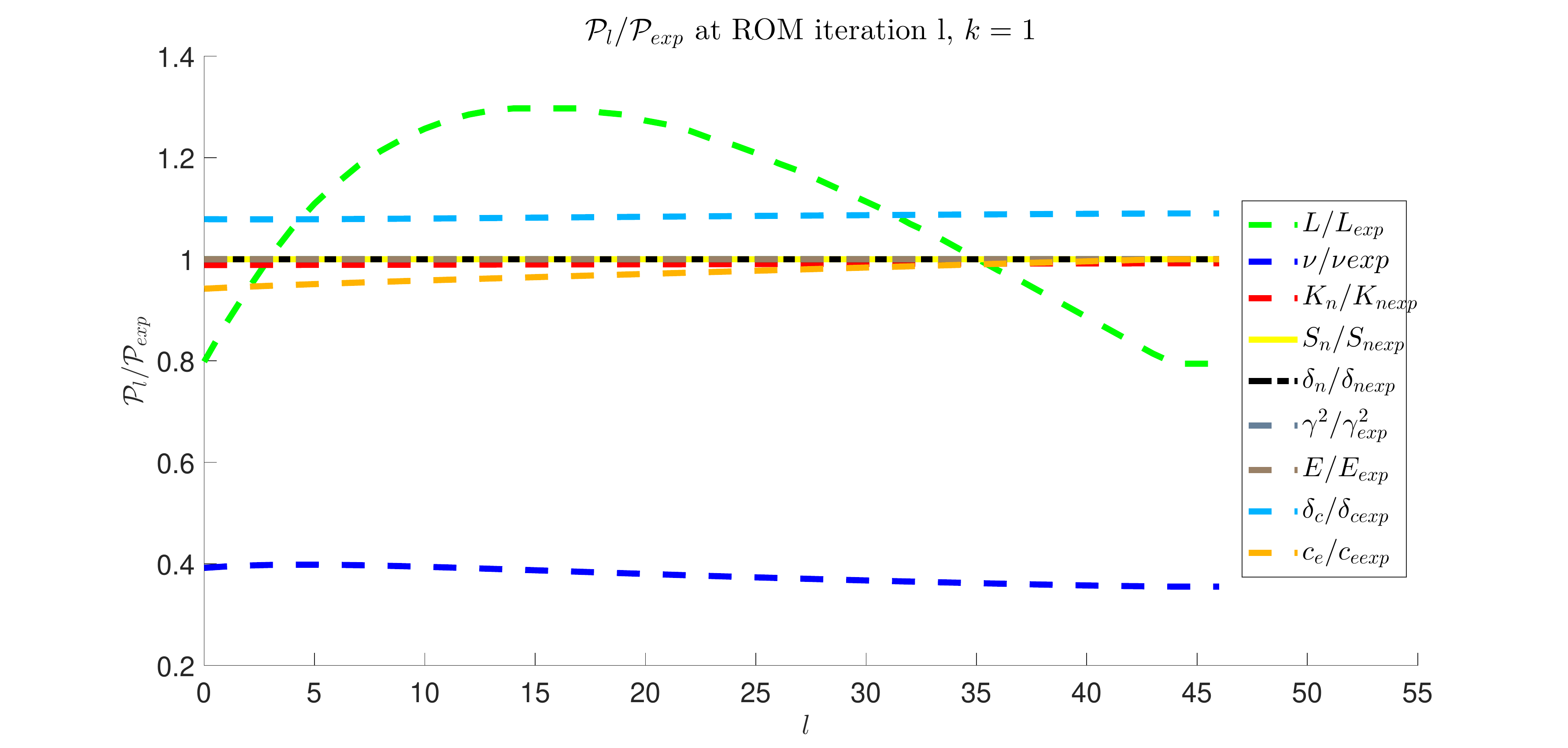}}&\\ \\
\parbox[t][][t]{5.6cm}{\centering $|\mathcal{P}_l(1)-\mathcal{P}_l|$} & & \\  
  \raisebox{-0.5\height}{\includegraphics[scale=0.23]{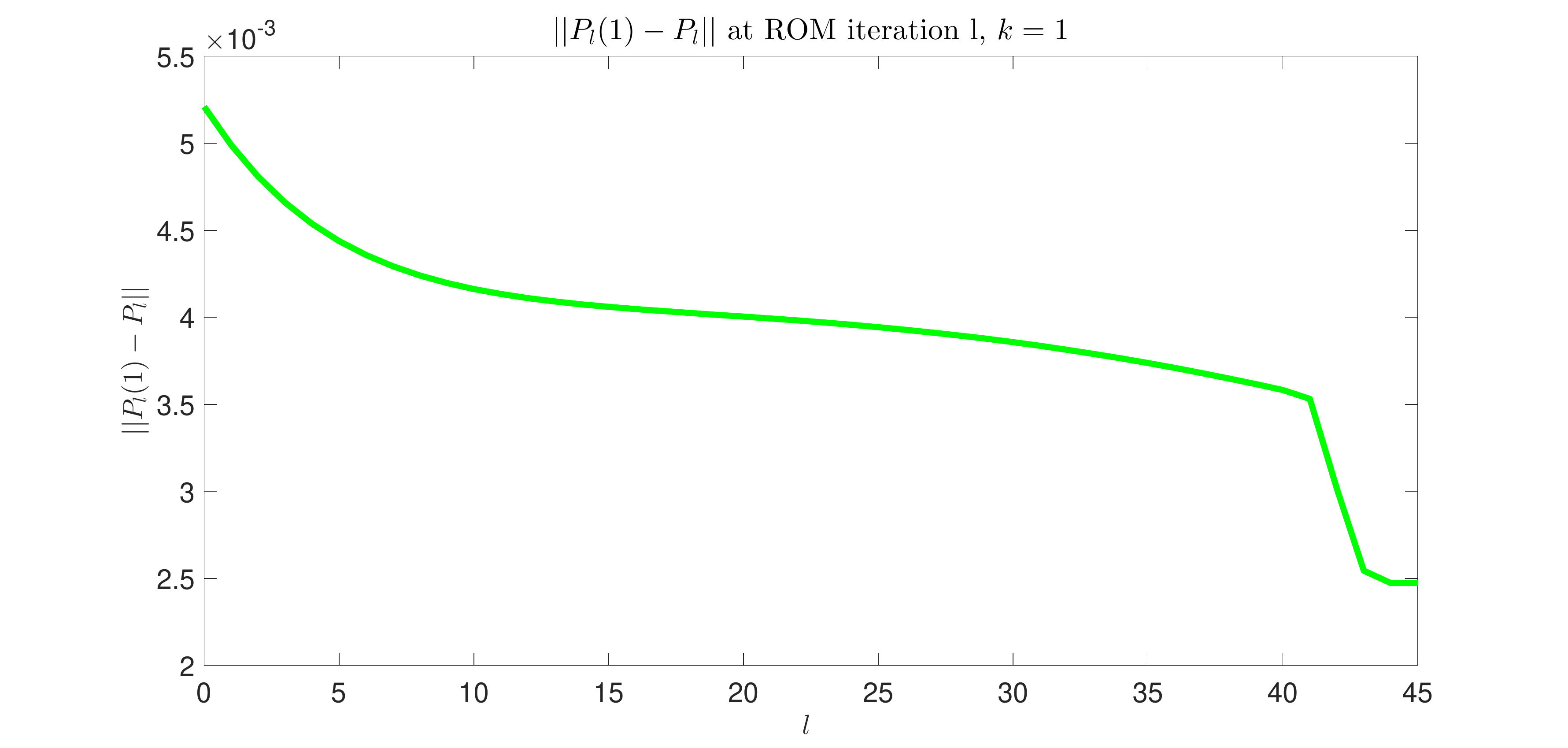}}& &\\ \\  
 \end{tabular}\\
 \end{minipage} }
\end{table}
\newpage
\begin{table}[h!]
\resizebox{0.9\textwidth}{!}{\begin{minipage}{\textwidth}
 \begin{tabular}{c c c}
 \parbox[t][][t]{3.2cm}{\centering\textbf{MRI}} &   \parbox[t][][t]{7.2cm}{\centering\textbf{ROM ($\mathcal{P}_{\bar{l}=45})$}}  &  \parbox[t][][t]{4.2cm}{\centering\textbf{Comparison}}  \\ \\
\end{tabular}
\begin{tabular}{l}
\raisebox{-0.5\height}{\includegraphics[scale=0.8]{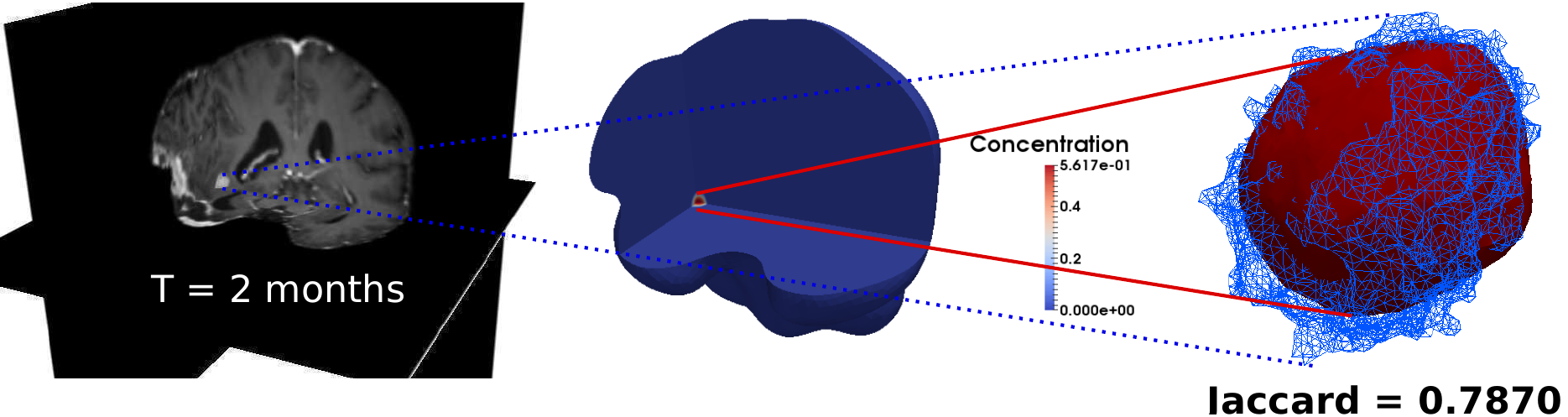}}  \\
  \end{tabular}
   \begin{tabular}{ c c c }
 \hline
 \parbox[t][][t]{5.2cm}{\centering \textbf{Iteration k=2}} & &\\ \hline \\ 
 \parbox[t][][t]{5.6cm}{\centering $\mathbf{J}(\vec{\alpha}_l,\mathcal{P}_l)$} &\parbox[t][][t]{6.6cm}{\centering $|\mathcal{P}_l(1)-\mathcal{P}_l|$} &\\  \
  \raisebox{-0.5\height}{\includegraphics[scale=0.23]{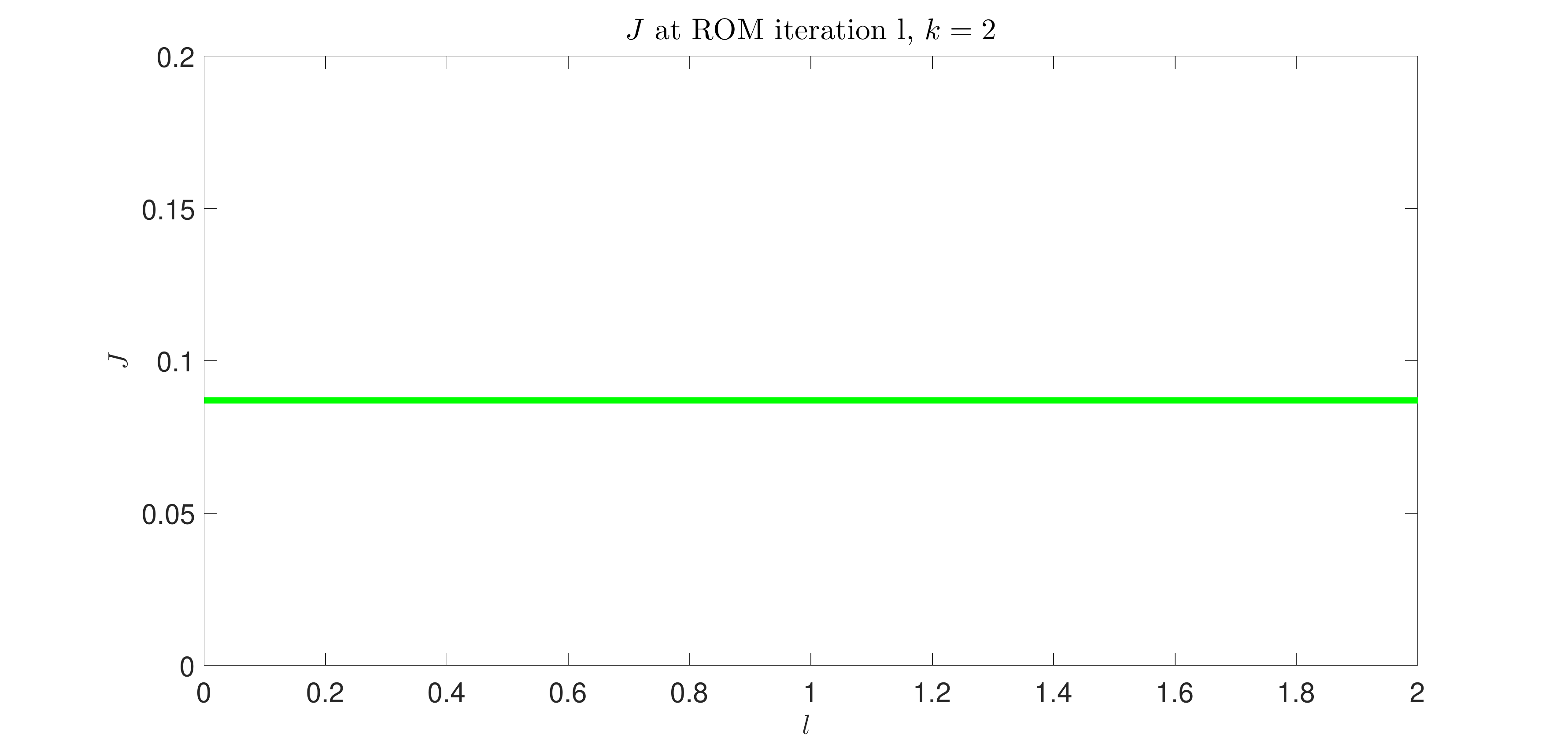}}&\raisebox{-0.5\height}{\includegraphics[scale=0.23]{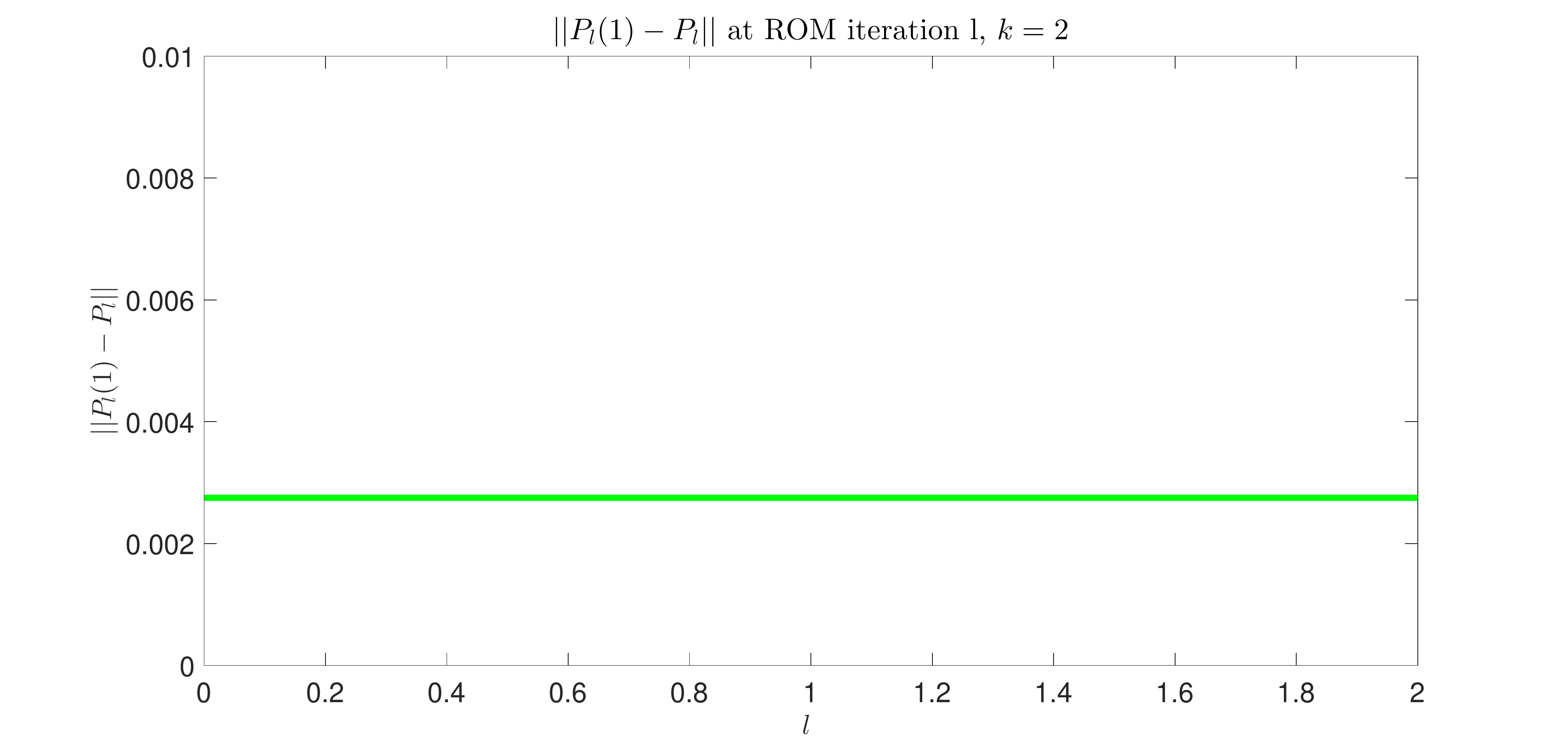}}&\\ \\ 
 \end{tabular}\\
  \end{minipage} }
  \captionof{figure}[]{Values of the functional $J(\vec{\alpha}_l,\mathcal{P}_l)$, of the normalised set of parameters $\mathcal{P}_l/\mathcal{P}_{\text{exp}}$ and of $|\mathcal{P}_l(1)-\mathcal{P}_l|$ for steps $k=0,1,2$ of \textbf{Algorithm 1}, together with a comparison between the isosurfaces $\phi_{\text{data}}(T)=0$ (highlighted in blue color) and $\sum_{i=1}^{N_{\text{POD}}}\alpha_{i\bar{l}}^N\xi_{i}^{\phi}=\phi_e/2$ (highlighted in red color).}
  \label{fig:7}
\end{table}

\subsubsection{Computational cost}
We conclude by reporting in Table \ref{tab:3} the elapsed time (in seconds) for the computation of steps $1-4$ of \textbf{Algorithm 1} for the different values of the iteration step $k$.
\begin{table}[h!]
\caption{Elapsed time (in seconds) for the computation of steps $1-4$, for the different values of the iteration step $k$.}
\label{tab:3}    
\begin{center}
\begin{tabular}{lllll}
\hline
\parbox[t][][t]{2.2cm}{\centering \textbf{Iteration}\\\textbf{k=0}} & \parbox[t][][t]{2.2cm}{\centering \textbf{Step 1}\\$201914.27$} & \parbox[t][][t]{2.2cm}{\centering \textbf{Step 2}\\$12.42$} & \parbox[t][][t]{2.2cm}{\centering \textbf{Step 3}\\$6281.31$} & \parbox[t][][t]{2.2cm}{\centering \textbf{Step 4}\\$157.8289$}\\ \\
\hline
 \parbox[t][][t]{2.2cm}{\centering \textbf{Iteration}\\\textbf{k=1}} & \parbox[t][][t]{2.2cm}{\centering \textbf{Step 1}\\$201425.71$} & \parbox[t][][t]{2.2cm}{\centering \textbf{Step 2}\\$12.89$} & \parbox[t][][t]{2.2cm}{\centering \textbf{Step 3}\\$2600.37$} & \parbox[t][][t]{2.2cm}{\centering \textbf{Step 4}\\$587.07$}\\ \\  
 \hline
 \parbox[t][][t]{2.2cm}{\centering \textbf{Iteration}\\\textbf{k=2}} & \parbox[t][][t]{2.2cm}{\centering \textbf{Step 1}\\$213040.53$} & \parbox[t][][t]{2.2cm}{\centering \textbf{Step 2}\\11.44} & \parbox[t][][t]{2.2cm}{\centering \textbf{Step 3}\\2544.43} & \parbox[t][][t]{2.2cm}{\centering \textbf{Step 4}\\42.53}\\ \\
\hline
\end{tabular}
\end{center}
\end{table}
We observe that the computational time for the projected gradient iterations at the ROM level is $3$ to $4$ order of magnitude smaller than the time needed to solve the FOM problem. We thus can very efficiently minimise the functional $J$ at the ROM level, checking at the FOM level the effective functional decrease and updating the POD basis to span the parameters space. In order to have light ROM systems, we need to have higher order tensors in \textbf{Problem 4} \eqref{eqn:assemble} with low dimension. Indeed, with only $4$ or $5$ basis functions \textbf{Step 3} of Algorithm \eqref{alg:opt} requires a heavy computational cost, which is only $2$ orders of magnitude smaller than the time needed to solve the FOM problem.
\subsubsection{Some remarks on benchmark results}
\begin{figure}[!t]
\centering 
\includegraphics[width=0.9 \linewidth]{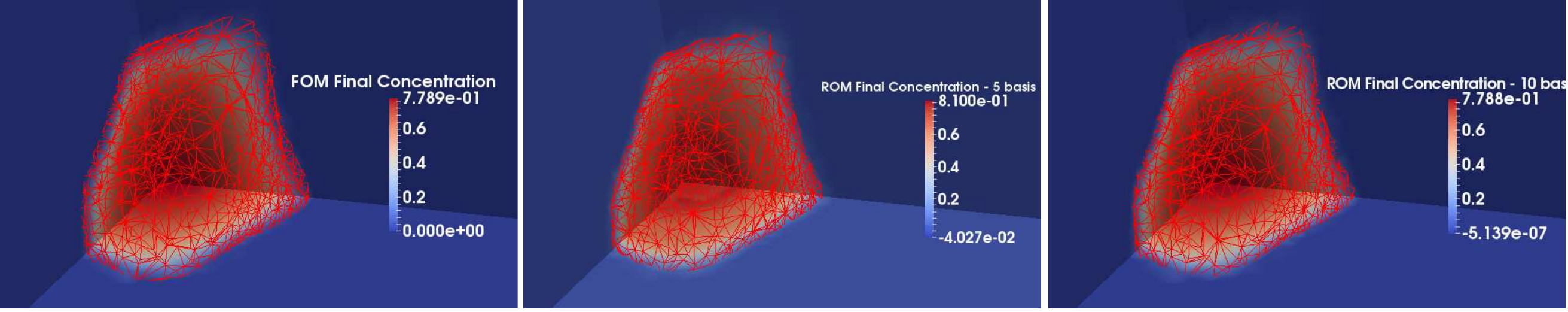}
\caption{Comparison between the final state $\phi_h^N$, solution of the FOM simulation with parameter set $\mathcal{P}_0^{\text{bis}}$ and $N=980$, and the corresponding final state $\sum_{i=1}^{N_{\text{POD}}}\alpha_{i0}^N\xi_{i}^{\phi}$, solution of the ROM system \eqref{eqn:6}. The iso--surfaces $\phi_h^N=\phi_e/2$ and $\sum_{i=1}^{N_{\text{POD}}}\alpha_{i0}^N\xi_{i}^{\phi}=\phi_e/2$ are highlighted in red colors.
}
\label{fig:6bis}
\end{figure}

In this paragraph we report some numerical results to show  how the POD analysis in \textbf{Step 2} of Algorithm \eqref{alg:opt} varies when the tumour concentration in the FOM simulations in \textbf{Step 1} is spreading on a larger region than the one observed in Test Case $1$. This happens if we consider a tumour dynamics over a larger time interval and also if we consider an initial set of parameters which induces a larger and more anisotropic tumour expansion. Moreover, we report numerical results about the performance of \textbf{Step 4} of Algorithm \eqref{alg:opt} when a different initial set of parameter is considered and when a lower threshold of POD significance (namely $99.9\%$) is considered. \\
In Table \ref{tab:2bis} we report the POD analysis of the snapshot matrices obtained from the FOM solution \eqref{eqn:fom} at $k=0$ of \textbf{Algorithm 1} with $N=980$, i.e. when the tumour dynamics span a time interval of $120$ days. 
\begin{table}[h!]
\caption{Values of the cumulated fraction of $\text{tr}\mathbf{F}_1^T\mathbf{F}_1$,$\text{tr}\mathbf{F}_2^T\mathbf{F}_2$,$\text{tr}\mathbf{F}_3^T\mathbf{F}_3$,$\text{tr}(\psi_1'(\mathbf{F}_1))^T\psi_1'(\mathbf{F}_1)$,$\text{tr}(\psi_1''(\mathbf{F}_1))^T\psi_1''(\mathbf{F}_1)$ for the first eigenvalues with the highest magnitude.}
\label{tab:2bis}    
\begin{center}
\begin{tabular}{lllllll}
\hline
\parbox[t][][t]{2.2cm}{\centering \textbf{Iteration}\\\textbf{k=0}\\$N=980$} & Eigenvalue & \% $\text{tr}\mathbf{F}_1^T\mathbf{F}_1$ & \% $\text{tr}\mathbf{F}_2^T\mathbf{F}_2$ & \% $\text{tr}\mathbf{F}_3^T\mathbf{F}_3$ & \% $\text{tr}(\psi_1')^T\psi_1'$ & \% $\text{tr}(\psi_1'')^T\psi_1''$\\
\hline
& First & $90.2061$ & $99.6130$ & $99.8455$ & $99.5709$ & $97.7895$ \\ 
& Second & $98.5523$ & $99.8068$ & $99.9859$ & $99.9554$ & $99.6662$\\
& Third & $99.5893$ & $99.9866$ & $99.9958$ & $99.9836$ & $99.8917$\\ 
& Fourth & $99.8472$ & $99.9967$ & $99.9982$ & $99.9928$ & $99.9446$\\ 
& Fifth & $99.9401$ & $99.9989$ & $99.9992$ & $99.9968$ & $99.9805$\\
& Sixth & $99.9744$ & $99.9995$ & $99.9995$ & $99.9987$ & $99.9941$\\
& Seventh & $99.9884$ & $99.9998$ & $99.9998$ & $99.9995$ & $99.9978$\\
& Eighth & $99.9944$ & $99.9999$ & $99.9999$ & $99.9998$ & $99.9991$\\
 \hline
\end{tabular}
\end{center}
\end{table}\\
We thus have that $N_{\text{POD}}=N_{\phi}^{\text{POD}}=8$ for $k=0$.\\ 
In Table \ref{tab:2tris} we report the same POD analysis for a time span of $120$ days, choosing also an initial set of parameters  
\[
\mathcal{P}_0^{\text{bis}}\equiv \{(1/3205.13),(0.128),(3.00),(9*10^3),(6184),(0.2862),(819.98),(0.215),(0.3364)\};
\]
in order to observe a larger and more anisotropic spread of the initial tumour distribution during the dynamics.
\begin{table}[h!]
\caption{Values of the cumulated fraction of $\text{tr}\mathbf{F}_1^T\mathbf{F}_1$,$\text{tr}\mathbf{F}_2^T\mathbf{F}_2$,$\text{tr}\mathbf{F}_3^T\mathbf{F}_3$,$\text{tr}(\psi_1'(\mathbf{F}_1))^T\psi_1'(\mathbf{F}_1)$,$\text{tr}(\psi_1''(\mathbf{F}_1))^T\psi_1''(\mathbf{F}_1)$ for the first eigenvalues with the highest magnitude.}
\label{tab:2tris}    
\begin{center}
\begin{tabular}{lllllll}
\hline
\parbox[t][][t]{2.2cm}{\centering \textbf{Iteration}\\\textbf{k=0}\\$N=980$\\$\mathcal{P}_0^{\text{bis}}$} & Eigenvalue & \% $\text{tr}\mathbf{F}_1^T\mathbf{F}_1$ & \% $\text{tr}\mathbf{F}_2^T\mathbf{F}_2$ & \% $\text{tr}\mathbf{F}_3^T\mathbf{F}_3$ & \% $\text{tr}(\psi_1')^T\psi_1'$ & \% $\text{tr}(\psi_1'')^T\psi_1''$\\
\hline
& First & $89.5942$ & $96.9529$ & $99.2853$ & $97.0274$ & $93.7024$ \\ 
& Second & $97.6348$ & $97.8125$ & $99.9235$ & $99.7625$ & $99.3436$\\
& Third & $99.3231$ & $99.7177$ & $99.9817$ & $99.9318$ & $99.7444$\\ 
& Fourth & $99.7375$ & $99.9534$ & $99.9927$ & $99.9631$ & $99.8768$\\ 
& Fifth & $99.8762$ & $99.9866$ & $99.9962$ & $99.9850$ & $99.9567$\\
& Sixth & $99.9370$ & $99.9944$ & $99.9978$ & $99.9930$ & $99.9849$\\
& Seventh & $99.9670$ & $99.9975$ & $99.9988$ & $99.9969$ & $99.9945$\\
& Eighth & $99.9821$ & $99.9987$ & $99.9992$ & $99.9986$ & $99.9974$\\
& Nineth & $99.9899$ & $99.9993$ & $99.9996$ & $99.9992$ & $99.9985$\\
& Tenth & $99.9942$ & $99.9996$ & $99.9997$ & $99.9996$ & $99.9992$\\
 \hline
\end{tabular}
\end{center}
\end{table}\\
We thus have that $N_{\text{POD}}=N_{\phi}^{\text{POD}}=10$ for $k=0$.\\ 
In Figure \ref{fig:5bis} we also show the basis elements $\xi_{i}^{\phi}$, corresponding to the highest eigenvalues needed to explain the $99.99\%$ variance of the data, for the case of $N=980$ and initial set $\mathcal{P}_0^{\text{bis}}$, superposed with the initial condition and final distribution of cell concentration (highlighted by a distribution of green and red points respectively).
\newpage
\begin{table}[h!]
\resizebox{0.9\textwidth}{!}{\begin{minipage}{\textwidth}
\begin{tabular}{ c }
 \hline
 \parbox[t][][t]{17.4cm}{\centering \textbf{Iteration k=0}, $N=980$, $\mathcal{P}_0^{\text{bis}}$} \\  \hline \\
\end{tabular}\\
\begin{tabular}{l}
\raisebox{-0.5\height}{\centering \includegraphics[scale=0.9]{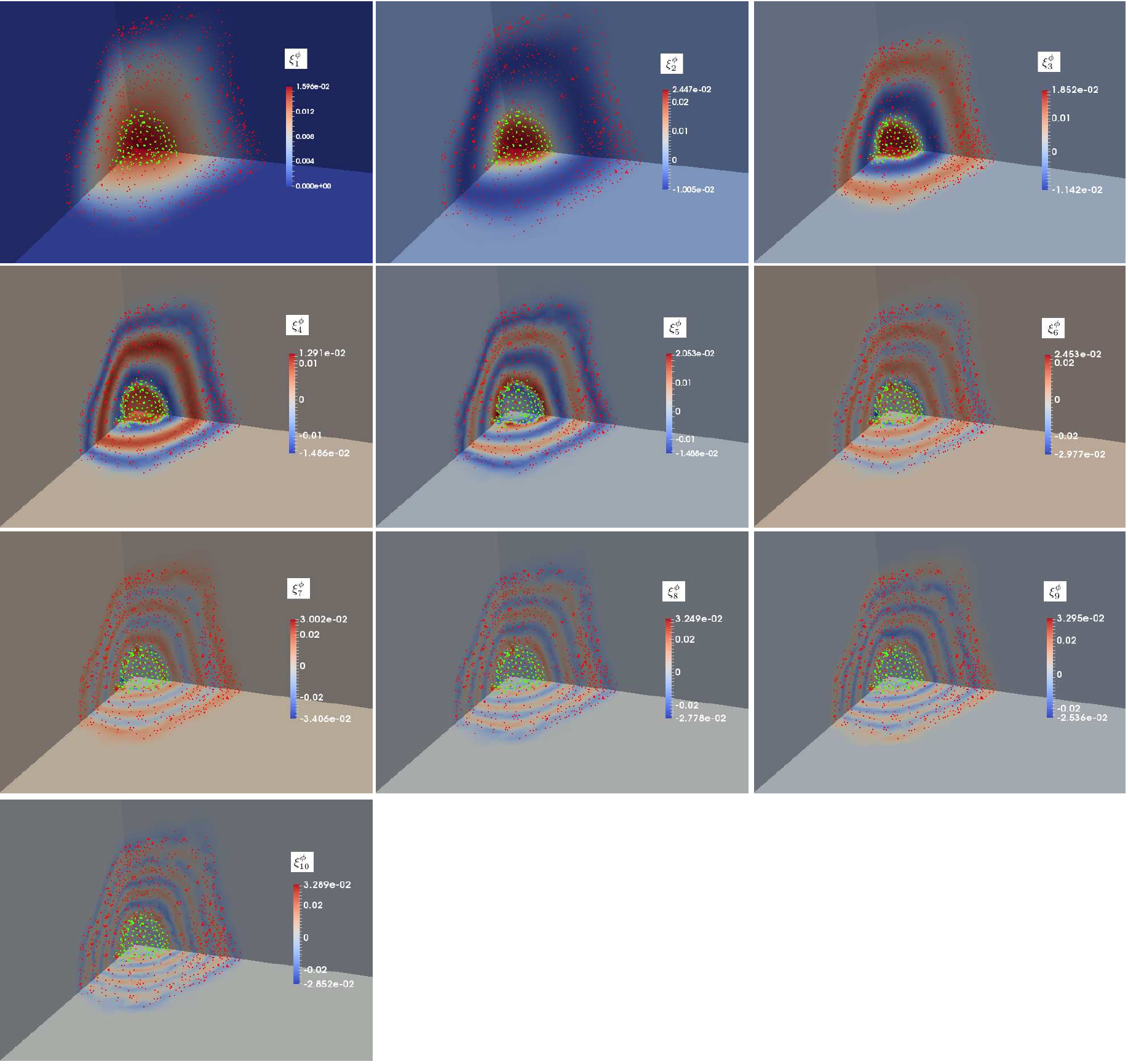}}  \\
  \hline
  \end{tabular}\\
  \end{minipage} }%
  \captionof{figure}[]{Plot of the basis elements $\xi_{i}^{\phi}$ corresponding to the highest eigenvalues needed to explain the $99.99\%$ variance of the data, for $k=0$, $N=980$, $\mathcal{P}_0^{\text{bis}}$.}
  \label{fig:5bis}
  \end{table}%
We can observe that, while $\xi_{1}^{\phi}$ and $\xi_{2}^{\phi}$ are distributed over the final state $\phi_h^N$ and the bulk of the initial condition $\phi_h^0$ respectively, the higher order basis $\xi_{i}^{\phi}$, $i=3,\dots, 10$ are oscillating functions over the set where the tumour is expanding during its temporal evolution, and thus contain the information about the tumour boundary and its expansion. Since the latter set is larger than the case shown in Figure \ref{fig:5}, the number of oscillating functions over this region with a frequency needed to explain the $99.99\%$ of the data is increased. 
\newpage
We show in Figure \ref{fig:6bis} a comparison between the final state $\phi_h^N$ calculated from the FOM simulation through \textbf{Algorithm 2} with parameter set $\mathcal{P}_0^{\text{bis}}$ and the corresponding final state $\sum_{i=1}^{N_{\text{POD}}}\alpha_{i0}^N\xi_{i}^{\phi}$ obtained as a solution of the ROM system \eqref{eqn:6} through \textbf{Algorithm 4}, with $N_{\text{POD}}=5$ and $N_{\text{POD}}=10$. 

We observe that the ROM solution is approximating the FOM solution with a very high fidelity when we consider $10$ POD basis (corresponding to a $99.99\%$ threshold of POD variance), whereas a low fidelity approximation is obtained when considering $5$ basis (corresponding to a $99.87\%$ threshold of POD variance).

In Figure \ref{fig:7bis} we report the values of the functional $J(\vec{\alpha}_l,\mathcal{P}_l)$, of the normalised set of parameters 
\[
\mathcal{P}_l/\mathcal{P}_{\text{exp}}=\{L_l/L_{\text{exp}},\nu_l/\nu_{\text{exp}},k_{nl}/k_{n\text{exp}},S_{nl}/S_{n\text{exp}},\delta_{nl}/\delta_{n\text{exp}},\gamma_l^2/\gamma_{\text{exp}}^2,E_l/E_{\text{exp}},\delta_c/\delta_{c\text{exp}},c_e/c_{e\text{exp}}\},
\]
and of $|\mathcal{P}_l(1)-\mathcal{P}_l|$, computed in \textbf{Steps 3} and \textbf{4}  of \textbf{Algorithm 1}, for $k=0$, obtained when starting from the initial set of parameters $\mathcal{P}_0^{\text{bis}}$ and considering the time span $N=490$. We note that in this case we need $N_{\text{POD}}=6$ to explain the $99.99\%$ of variance of the data. Moreover, in order for the Algorithm \eqref{eqn:pgl} to converge we need to choose $n_w=2$.

\begin{table}[h!]
\resizebox{0.9\textwidth}{!}{\begin{minipage}{\textwidth}
\begin{tabular}{ c c c }
 \hline
 \parbox[t][][t]{5.2cm}{\centering \textbf{Iteration k=0}, $N=490$, $\mathcal{P}_0^{\text{bis}}$} & &\\ \hline \\ 
 \parbox[t][][t]{5.6cm}{\centering $\mathbf{J}(\vec{\alpha}_l,\mathcal{P}_l)$} &\parbox[t][][t]{6.6cm}{\centering $\mathcal{P}_l/\mathcal{P}_{\text{exp}}$}& \\  
  \raisebox{-0.5\height}{\includegraphics[scale=0.23]{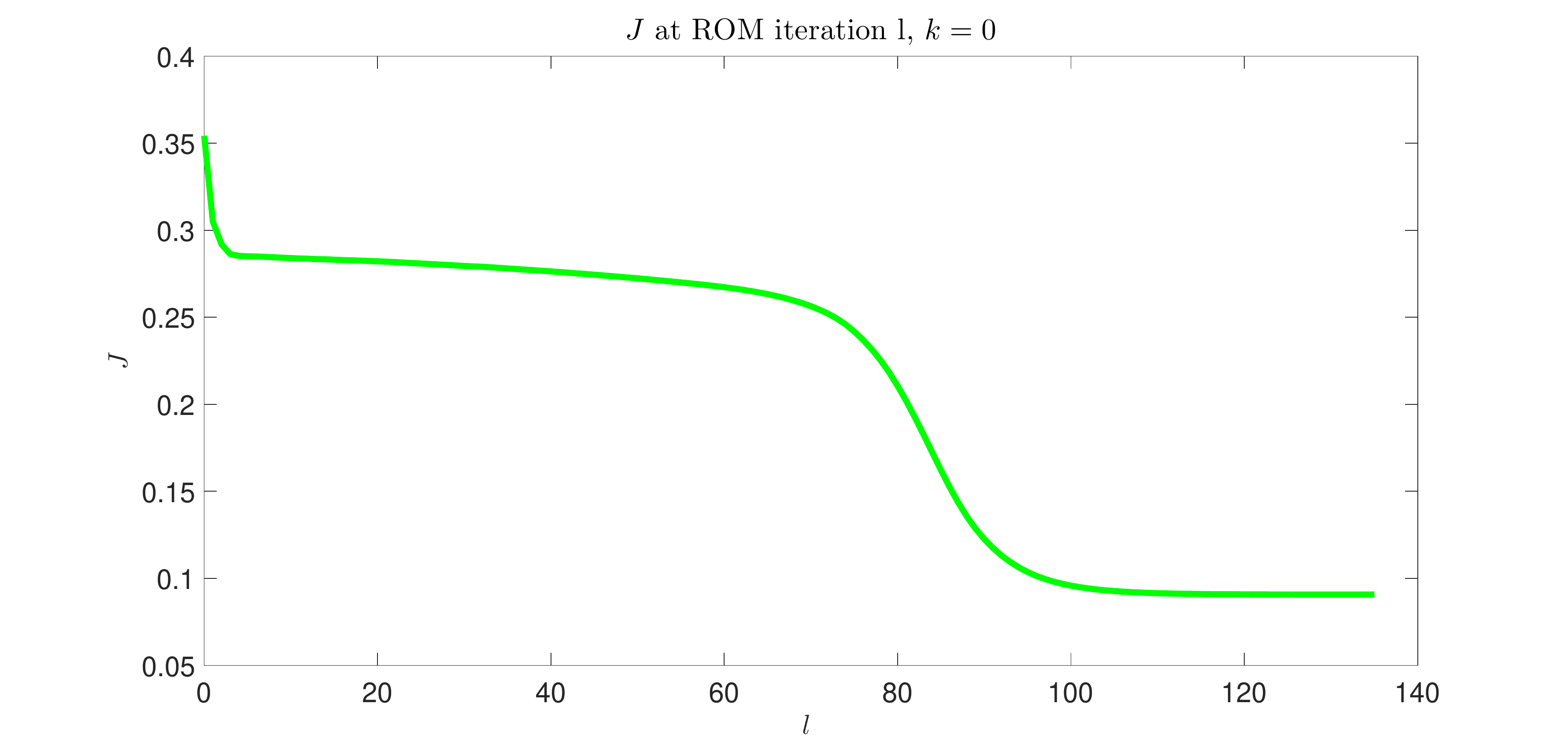}}&\raisebox{-0.5\height}{\includegraphics[scale=0.23]{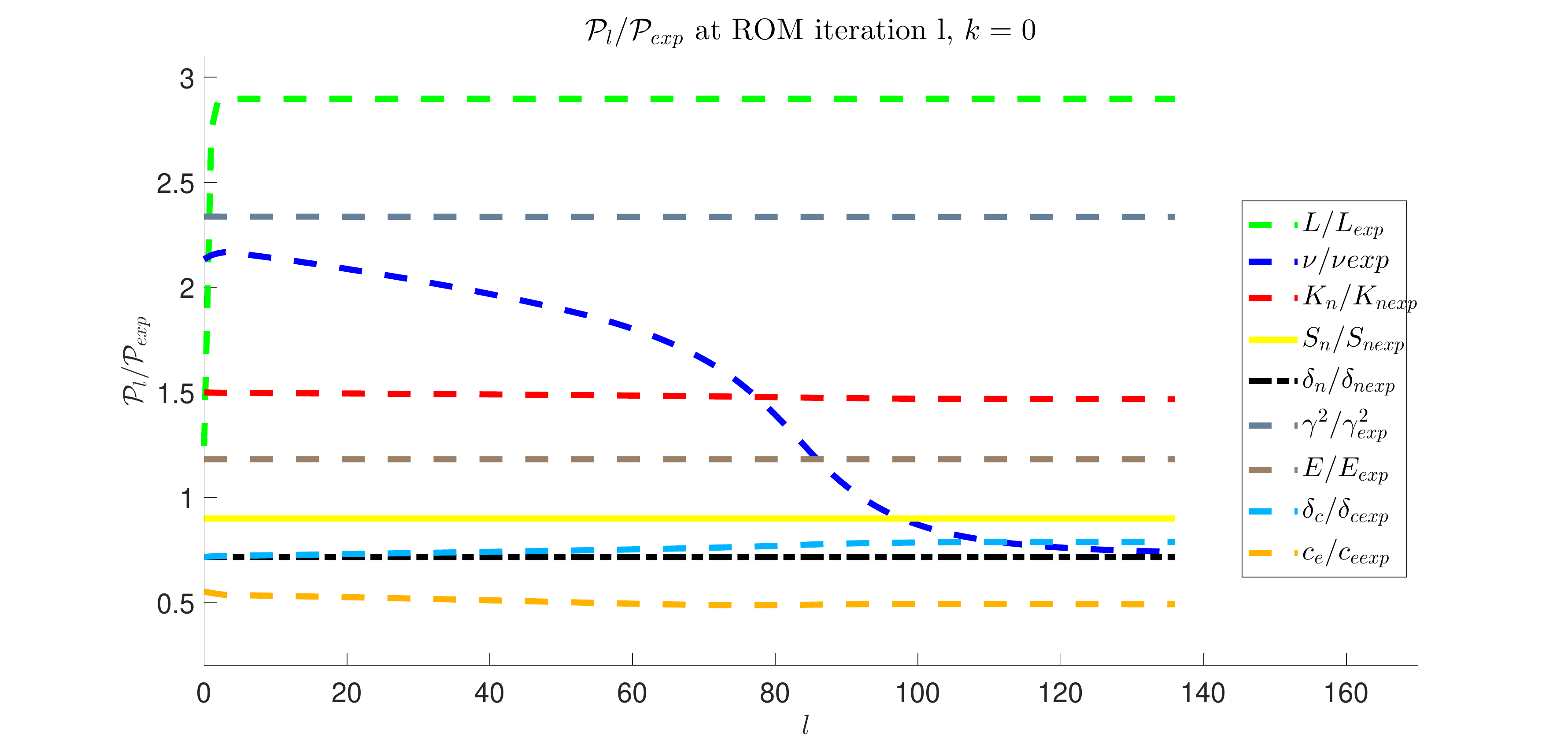}}&\\ \\
\parbox[t][][t]{5.6cm}{\centering $|\mathcal{P}_l(1)-\mathcal{P}_l|$} & & \\  
  \raisebox{-0.5\height}{\includegraphics[scale=0.23]{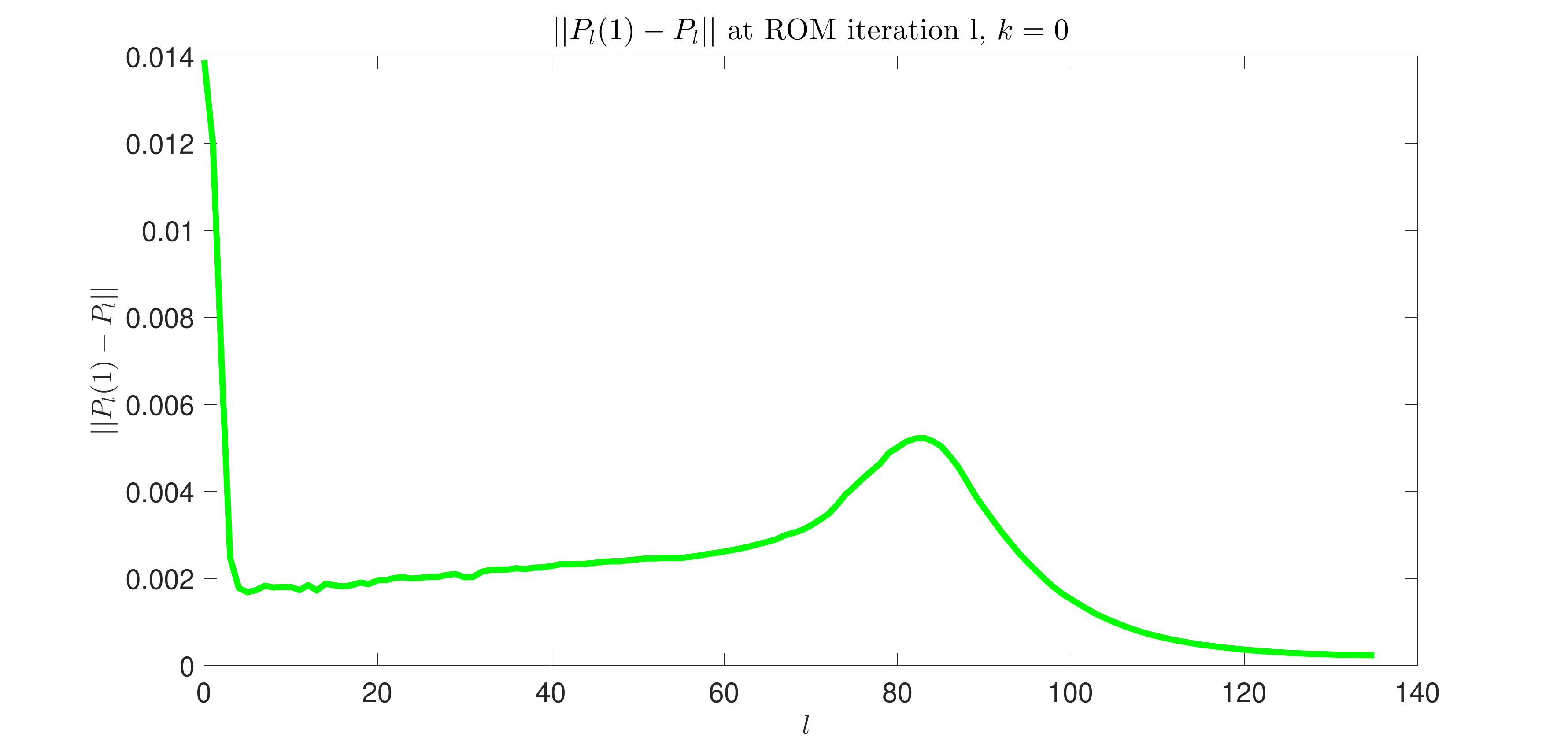}}& &\\ \\  
 \end{tabular}
  \end{minipage} }
  \captionof{figure}[]{Values of the functional $J(\vec{\alpha}_l,\mathcal{P}_l)$, of the normalised set of parameters $\mathcal{P}_l/\mathcal{P}_{\text{exp}}$ and of $|\mathcal{P}_l(1)-\mathcal{P}_l|$ for step $k=0$ of \textbf{Algorithm 1}, with $N=490$ and starting from the set $\mathcal{P}_0^{\text{bis}}$.}
  \label{fig:7bis}
\end{table}
We observe that the parameters $\nu$ and $L$ go through a large excursion. The parameter $\nu$ is decreasing to $1/3$ of its initial value, whereas the parameter $L$ reaches its active value $L_{\text{bio,max}}$ and remains stick to it, differently from the behaviour observed in Figure \ref{fig:7} where, starting from $\mathcal{P}_0$, the parameter $L$ reaches its maximum value and than relaxes to its initial value. The functional $J$, starting from a value $J=0.3543$ higher than $J=0.2847$ reported in Figure \ref{fig:7}, relaxes onto a minimum value $J=0.0906$ which is next to the value $J=0.0870$ reported in Figure \ref{fig:7} in a number of steps $\bar{l}=136$ much higher than $\bar{l}=11$ in Figure \ref{fig:7}.\\
Finally, we consider the results of the \textbf{Steps 2--4} of Algorithm \eqref{alg:opt}, for the first step $k=0$, with starting point $\mathcal{P}_0$, $N=490$ and when a threshold value of $99.9\%$ is considered in the POD analysis. In Figure \ref{fig:7tris} we report the values of the functional $J(\vec{\alpha}_l,\mathcal{P}_l)$, of the normalised set of parameters 
\[
\mathcal{P}_l/\mathcal{P}_{\text{exp}}=\{L_l/L_{\text{exp}},\nu_l/\nu_{\text{exp}},k_{nl}/k_{n\text{exp}},S_{nl}/S_{n\text{exp}},\delta_{nl}/\delta_{n\text{exp}},\gamma_l^2/\gamma_{\text{exp}}^2,E_l/E_{\text{exp}},\delta_c/\delta_{c\text{exp}},c_e/c_{e\text{exp}}\},
\]
and of $|\mathcal{P}_l(1)-\mathcal{P}_l|$, computed in \textbf{Steps 3} and \textbf{4}  of \textbf{Algorithm 1}. We note that in this case we need $N_{\text{POD}}=4$ to explain the $99.9\%$ of variance of the data. 

\begin{table}[h!]
\resizebox{0.9\textwidth}{!}{\begin{minipage}{\textwidth}
\begin{tabular}{ c c c }
 \hline
 \parbox[t][][t]{5.2cm}{\centering \textbf{Iteration k=0}, $99.9\%$ POD threshold, $N=490$, $\mathcal{P}_0$} & &\\ \hline \\ 
 \parbox[t][][t]{5.6cm}{\centering $\mathbf{J}(\vec{\alpha}_l,\mathcal{P}_l)$} &\parbox[t][][t]{6.6cm}{\centering $\mathcal{P}_l/\mathcal{P}_{\text{exp}}$}& \\  
  \raisebox{-0.5\height}{\includegraphics[scale=0.23]{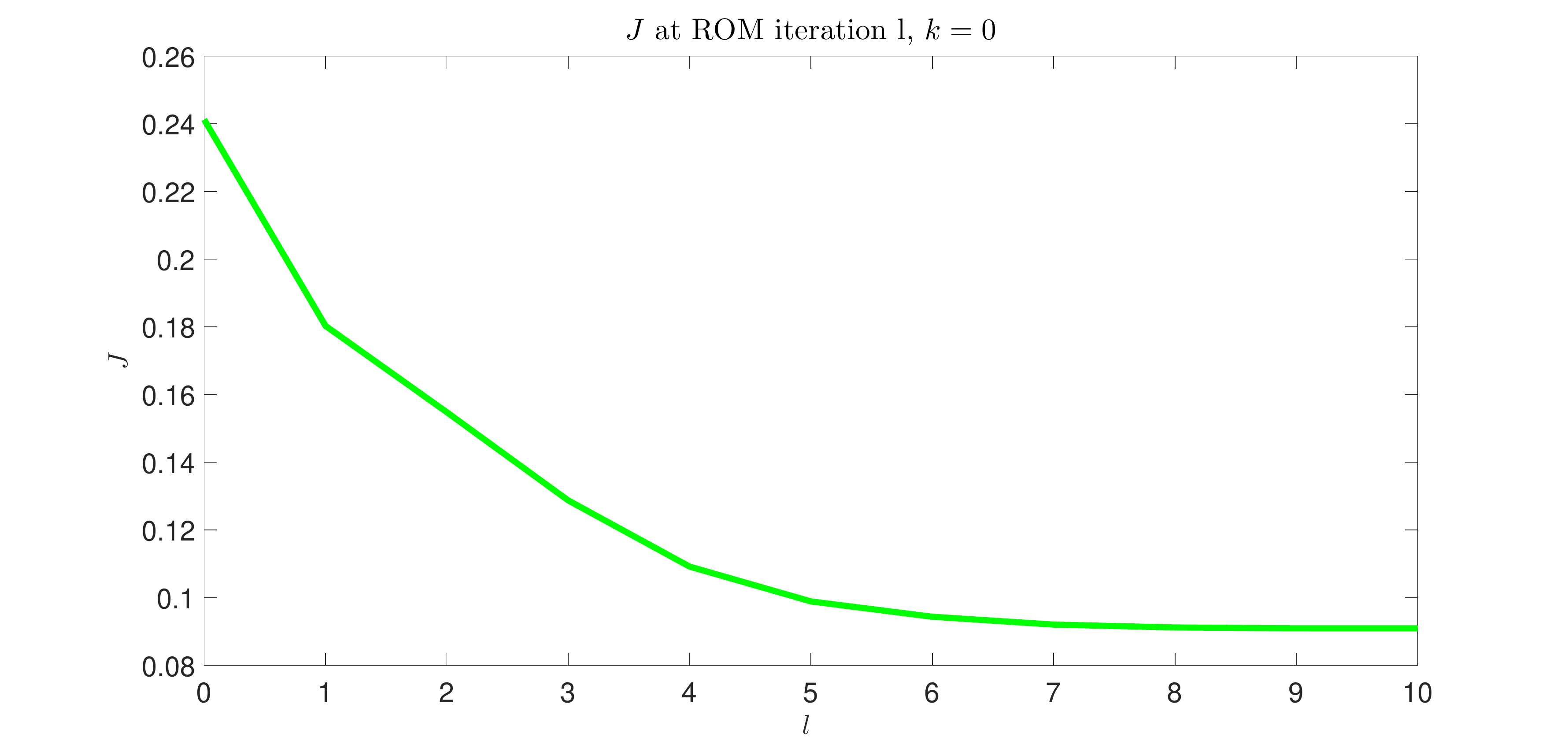}}&\raisebox{-0.5\height}{\includegraphics[scale=0.23]{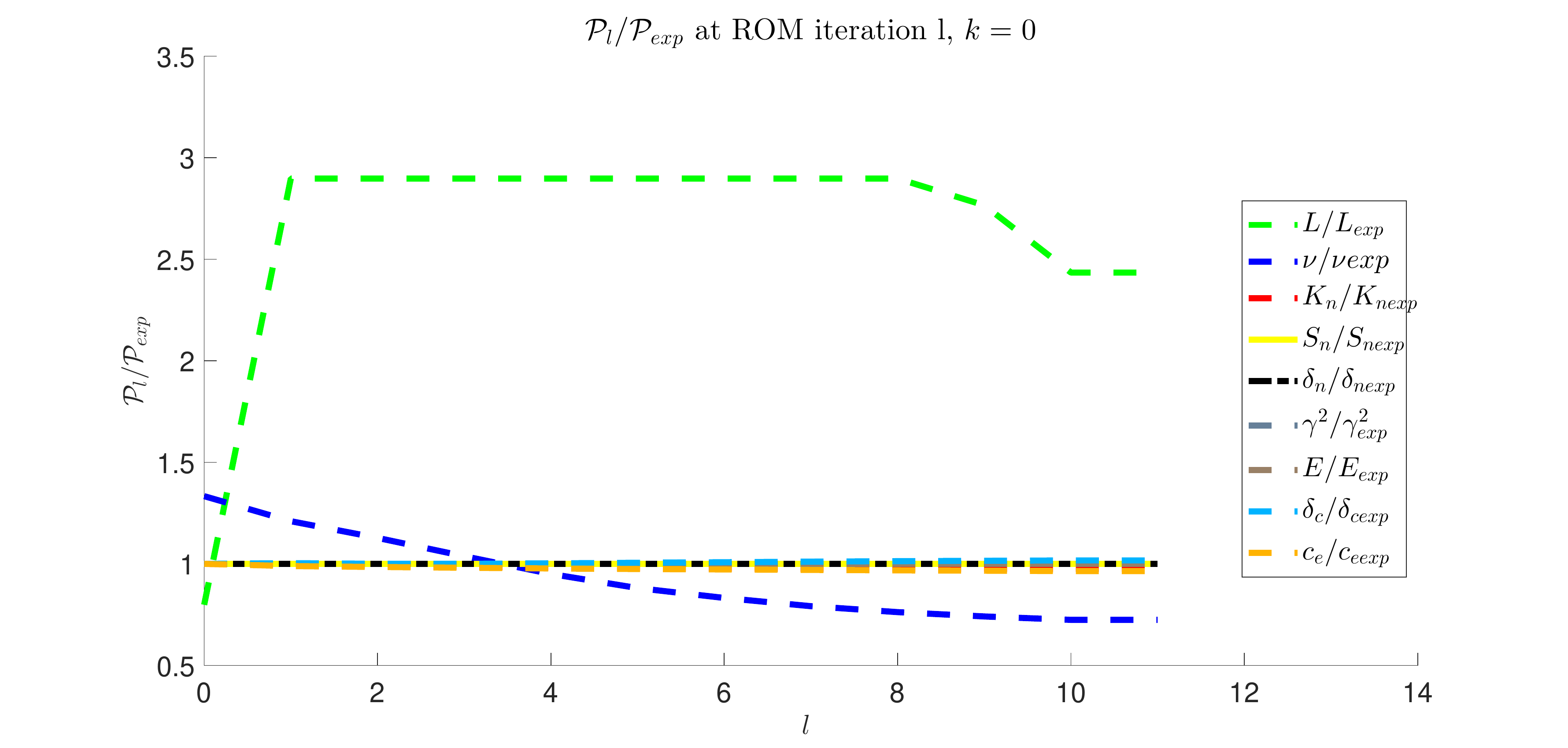}}&\\ \\
\parbox[t][][t]{5.6cm}{\centering $|\mathcal{P}_l(1)-\mathcal{P}_l|$} & & \\  
  \raisebox{-0.5\height}{\includegraphics[scale=0.23]{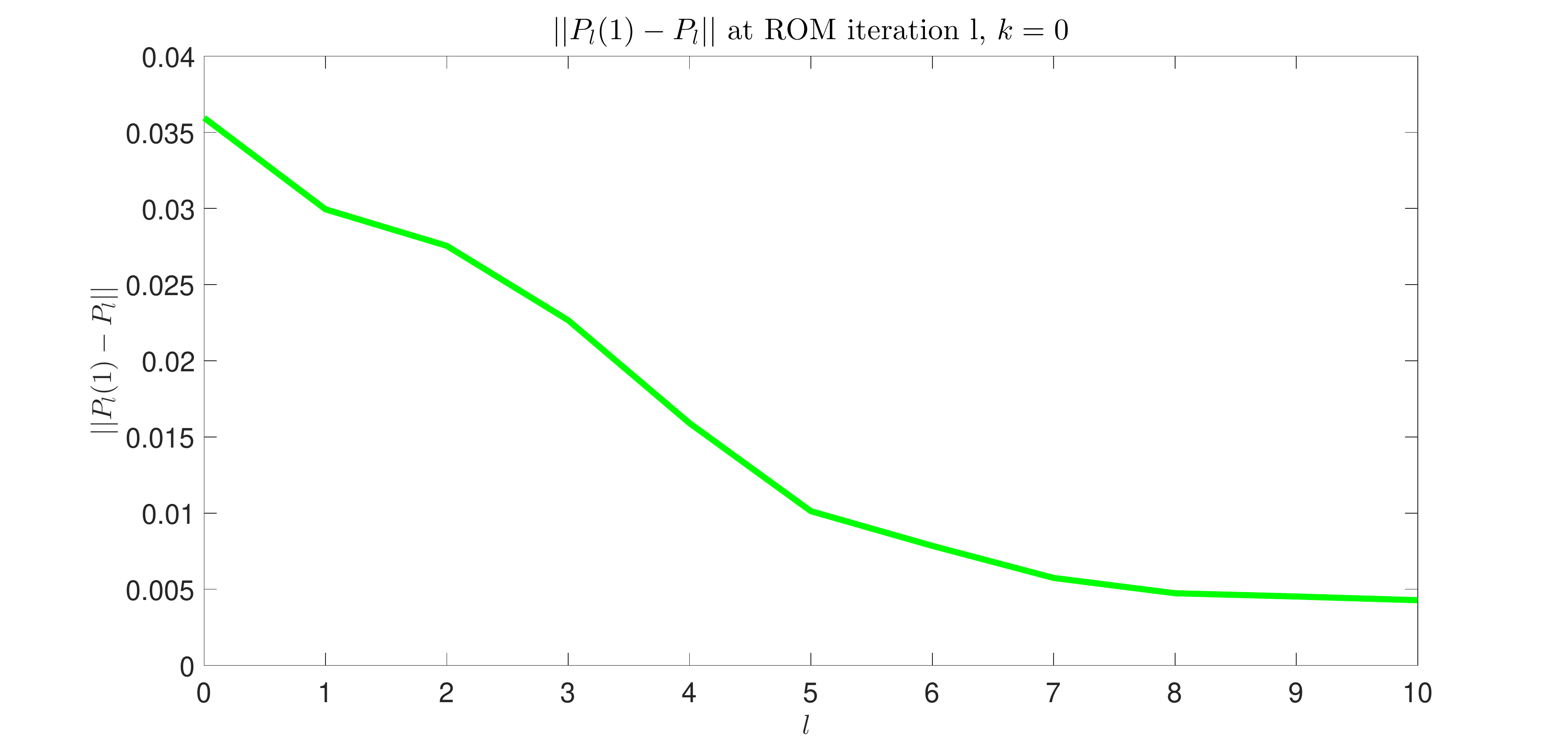}}& &\\ \\  
 \end{tabular}
 \begin{tabular}{c c c}
 \parbox[t][][t]{3.2cm}{\centering\textbf{MRI}} &   \parbox[t][][t]{7.2cm}{\centering\textbf{ROM ($\mathcal{P}_{\bar{l}=10})$}}  &  \parbox[t][][t]{4.2cm}{\centering\textbf{Comparison}}  \\ \\
\end{tabular}\\
\begin{tabular}{l}
\raisebox{-0.5\height}{\includegraphics[scale=0.8]{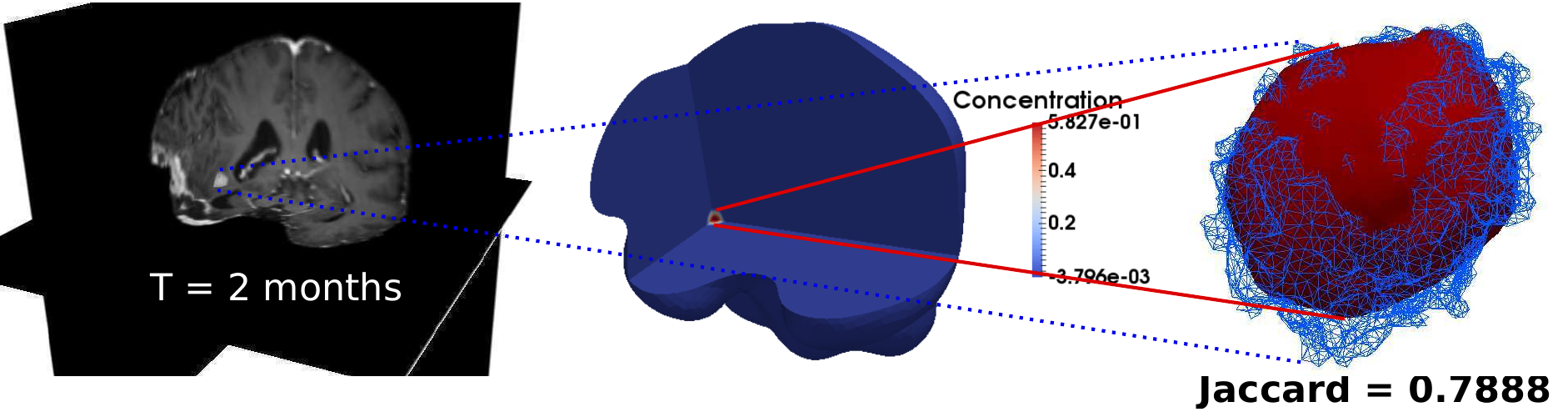}}  \\
  \end{tabular}
  \end{minipage} }
  \captionof{figure}[]{Values of the functional $J(\vec{\alpha}_l,\mathcal{P}_l)$, of the normalised set of parameters $\mathcal{P}_l/\mathcal{P}_{\text{exp}}$ and of $|\mathcal{P}_l(1)-\mathcal{P}_l|$for steps $k=0$ of \textbf{Algorithm 1}, together with a comparison between the isosurfaces $\phi_{\text{data}}(T)=0$ (highlighted in blue color) and $\sum_{i=1}^{N_{\text{POD}}}\alpha_{i\bar{l}}^N\xi_{i}^{\phi}=\phi_e/2$ (highlighted in red color).}
  \label{fig:7tris}
\end{table}
We observe that in this case of a lower POD threshold of $99.9\%$ the value of the functional $J$ reaches a minimum value $0.0909$ which is higher than the minimal value $0.0870$ attained during the optimization algorithm with an higher threshold of $99.99\%$ (see Figure \ref{fig:7}).
Moreover, in the former case the parameter $\nu$ varies by a smaller amount, whereas the parameter $L$, once reached the active value $L_{\text{max}}$, changes of a smaller amount.
\newpage
\subsection{Test case $2$: clinical follow-up after surgical resection and recurrence}
\begin{figure}[!h]
\centering
\includegraphics[width=0.7\linewidth]{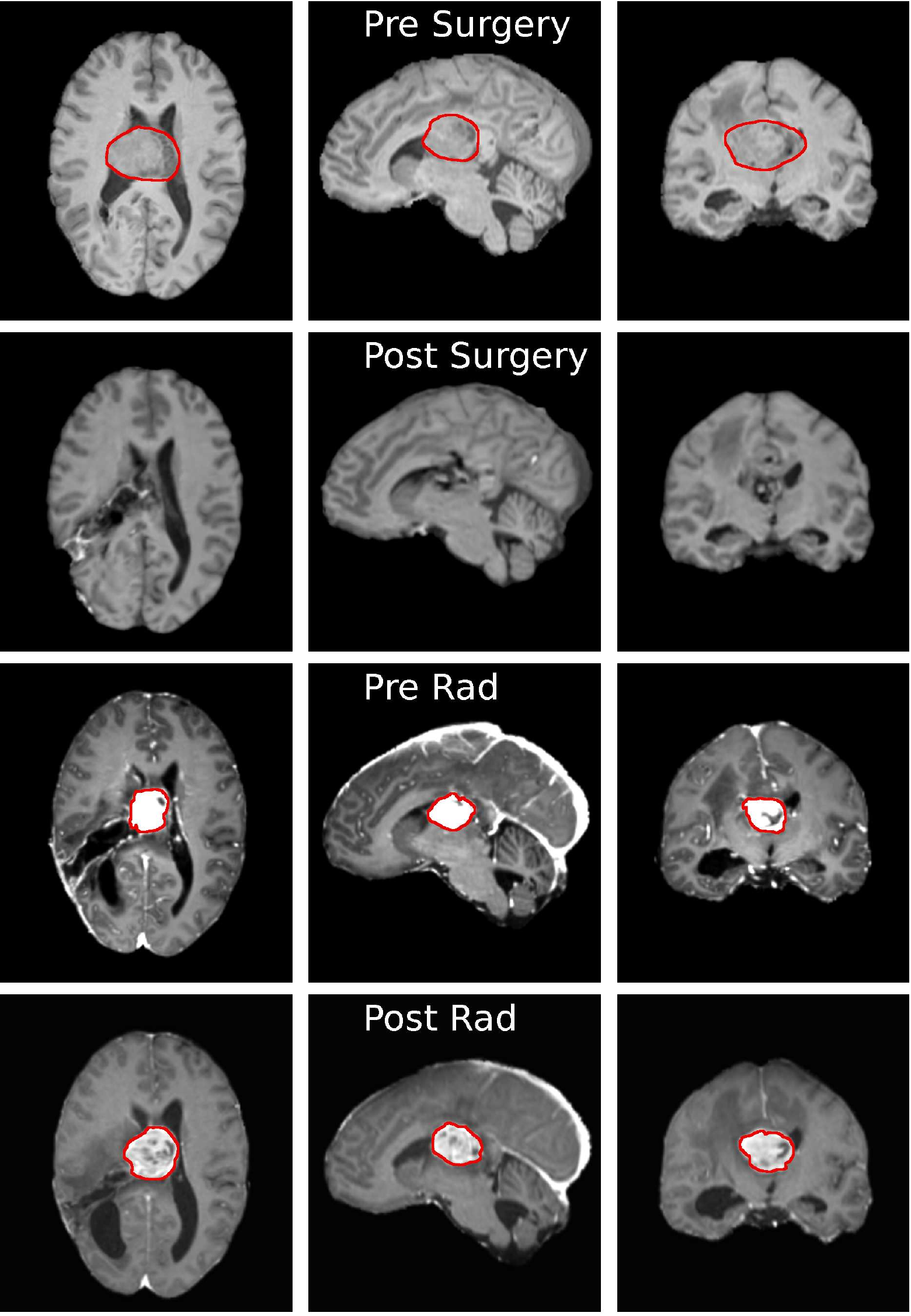}
\caption{ Axial (first column), sagittal (second column) and coronal (third column) slices of the T1-weighted MRI at different temporal stages. First row: before surgery; second raw: after surgery; third row: 34 days after surgery; fourth row: 78 days after surgery. It is possible to appreciate the subtotal resection of the corpus callosum GBM and the early tumour relapse at Pre Rad MRI and at Post Rad MRI. The segmented boundary of the tumour is highlighted in red color.}
\label{fig:1test2}
\end{figure}

Secondly, we apply the Optimization Algorithm so a  clinical test case which followed the surgical resection and the  recurrence pattern of a GBM. \\
A patient diagnosed with giant cell GBM underwent subtotal tumour removal. The patient started radiotherapy with concomitant Temozolomide  42 days after surgery following  to Stupp protocol. The pre-Radiotherapy MRI showed tumour relapse. After 25 doses of RT, the patient had a severe  worsening of the clinical status.\\ 
MRI were taken at the pre-operative, immediate post-operative, pre-radiotherapy ($34$ days after surgery) temporal stages and lastly at $5$ days after the interruption of RT and concomitant CHT due to disease progression. \\
Our numerical simulations  focus on  the period starting  with the surgical removal (initial time $t=0$ of simulations)  up to the first application of radiotherapy (time $t=T=:34$ days).
At $t=T$ we compare data and simulations, searching for the optimal set of parameters $\mathcal{P}_{\text{opt}}$ which locally minimises the functional \eqref{eqn:2}, obtained by solving \textbf{Algorithm 1}, thus estimating the model parameters directing the recurrence growth without any adjuvant therapy.\\
In Figure \ref{fig:1test2} we show the axial, sagittal and coronal slices of the T1-weighted MRI at different temporal stages.

We can observe that, after the application of $25$ fractions of RT, at $t=78$ days after surgery (Post Rad event) the GBM recurrence has grown in volume with respect to the Pre Rad event, inducing a severe
and rapid worsening of the patient's clinical status. The aim of the present work is to study the patient specific optimization of the growth parameters in the temporal range of tumour evolution after surgery and before the application of radiotherapy and chemotherapy.
\newpage
\subsubsection{Initialisation}
\begin{figure}[!ht]
\centering
\includegraphics[width=0.7\linewidth]{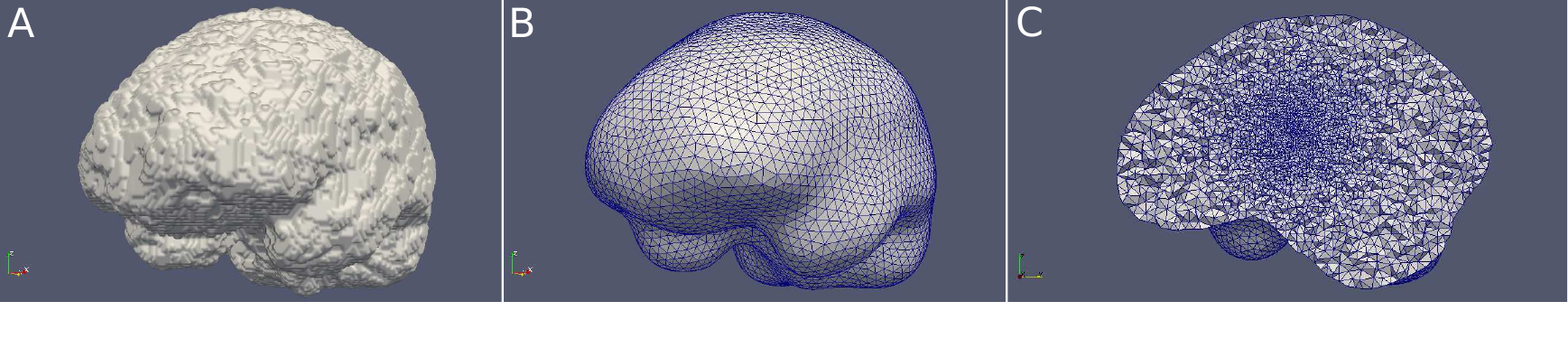}
\caption{ (A) External brain surface extracted from the medical images; (B) Smoothed and re--meshed external surface; (C) Tetrahedral mesh generated within
the external surface, conveniently refined in
the area surrounding the tumour.
}
\label{fig:8}
\end{figure}

\begin{figure}[!ht]
\centering
\includegraphics[width=0.7 \linewidth]{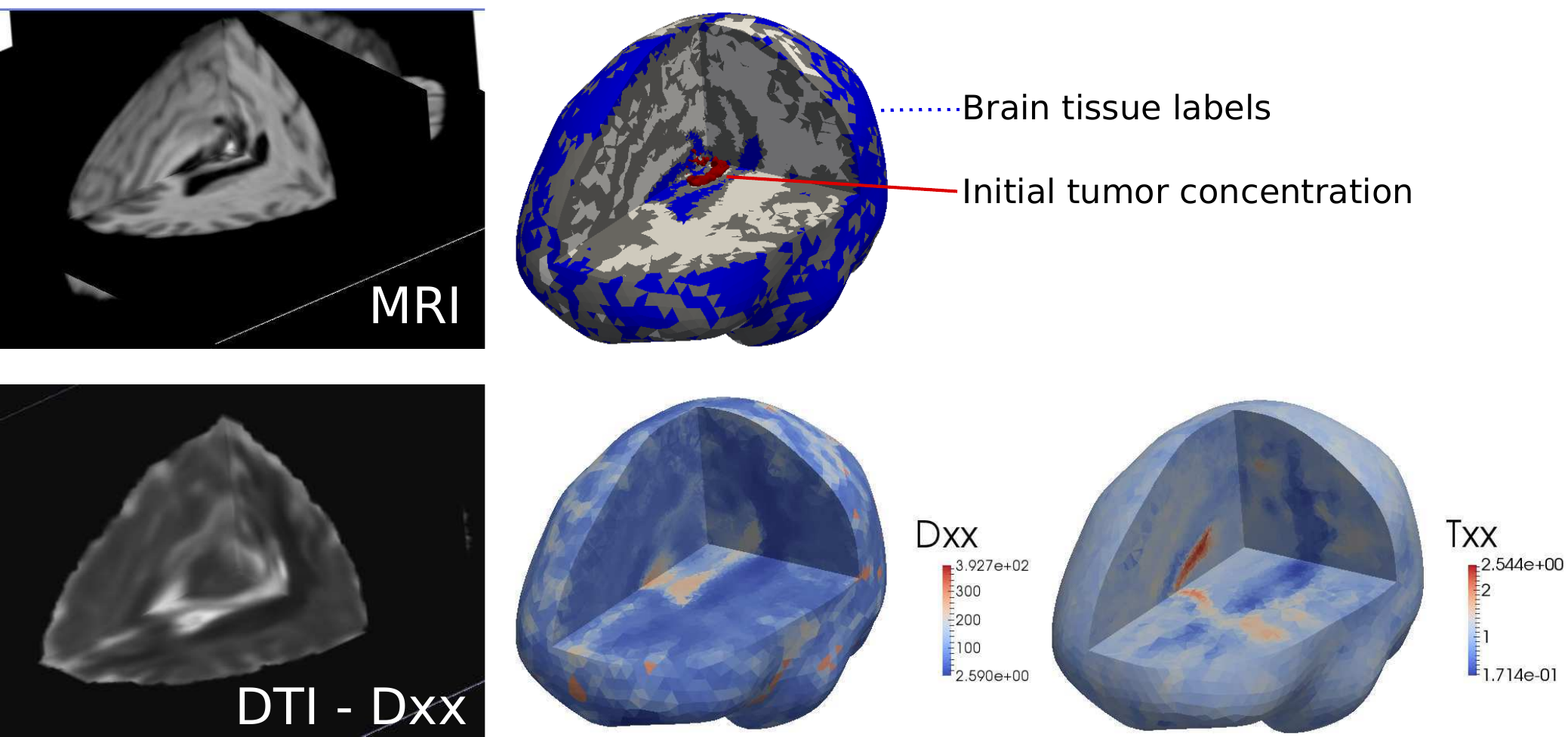}
\caption{ 3D virtual reconstructions of the MRI (top) and DTI (bottom) data, with the corresponding
computational meshes containing the labels of the brain tissues and the xx component of the tensors \textbf{D} and \textbf{T}. White matter, grey matter and CSF are highlighted in white, grey
and blue colors respectively. The initial tumour distribution is also highlighted in brown color, as segmented from the T1 MRI.
}
\label{fig:9}
\end{figure}

In Figures \ref{fig:8} and \ref{fig:9} we represent the results of the \textbf{initialisation} step of \textbf{Algorithm 1}, which defines the domain $\Omega$ (Figure \ref{fig:8}), the map(WM,GM,CSF), the initial condition $\phi_h^0$ and the tensors \textbf{D} and \textbf{T} (Figure \ref{fig:9}). 

\noindent The number of elements and the number of nodes of the mesh $\mathcal{T}_h$ are $305489$ and $51005$ respectively.  Moreover, we choose $\Delta T=0.1225$ (days), than $N=280$. As in Test Case $1$, a good refinement of the mesh in the region of tumour evolution is necessary to obtain ROM systems with sufficiently  low dimensionality to be solved with low computational resources and in highly reduced
computational times.
\\
In a similar way the characteristic function of the tumour extension $\phi_{\text{data}}(T)$ has been obtained from the segmentation of the MR images at $t=34$ days (PreRad event).
\subsubsection{Step 1}
In Figure \ref{fig:10} we report the values of the functional $J(\phi_h^N(\mathcal{P}_k),\mathcal{P}_k)$, calculated in \textbf{step 1} of \textbf{Algorithm 1}, and of the set of parameters $\mathcal{P}_k$, for different values of $k$. We also plot the isosurfaces $\phi_{\text{data}}(T)=0$ and $\phi_h^N(\mathcal{P}_k)=\phi_e/2$ from the FOM simulations, reporting the Jaccard index between the two volumes enclosed by these surfaces.
{\footnotesize
\centering
\begin{table}[p]
\resizebox{0.9\textwidth}{!}{\begin{minipage}{\textwidth}
\begin{tabular}{ c c c c c c c c c c c }
 \hline
 \parbox[t][][t]{2.2cm}{\centering \textbf{Iteration}\\\textbf{k=0}} & \parbox[t][][t]{1.6cm}{\centering $\mathbf{J}(\mathcal{P}_0)$\\0.28323} &\parbox[t][][t]{0.9cm}{\centering $\mathbf{L}_0$\\0.0002} & \parbox[t][][t]{0.9cm}{\centering $\boldsymbol{\nu}_0$\\0.08} &  \parbox[t][][t]{0.9cm}{\centering $\mathbf{k}_{n0}$\\2} & \parbox[t][][t]{0.9cm}{\centering $\mathbf{S}_{n0}$\\10000}  & \parbox[t][][t]{0.9cm}{\centering $\boldsymbol{\delta}_{n0}$\\8640} & \parbox[t][][t]{0.9cm}{\centering $\boldsymbol{\gamma}_0^2$\\0.1225} &  \parbox[t][][t]{0.9cm}{\centering $\mathbf{E}_{0}$\\694} &  \parbox[t][][t]{0.9cm}{\centering $\boldsymbol{\delta}_{0}$\\0.3} &  \parbox[t][][t]{0.9cm}{\centering $\textbf{c}_{e0}$\\0.611}\\  \hline \\
 \textbf{MRI} &  &  & &  \textbf{FOM} &   &  & \textbf{Comparison} &  & & \\ \\
\end{tabular}
\begin{tabular}{l}
\raisebox{-0.5\height}{\includegraphics[scale=0.92]{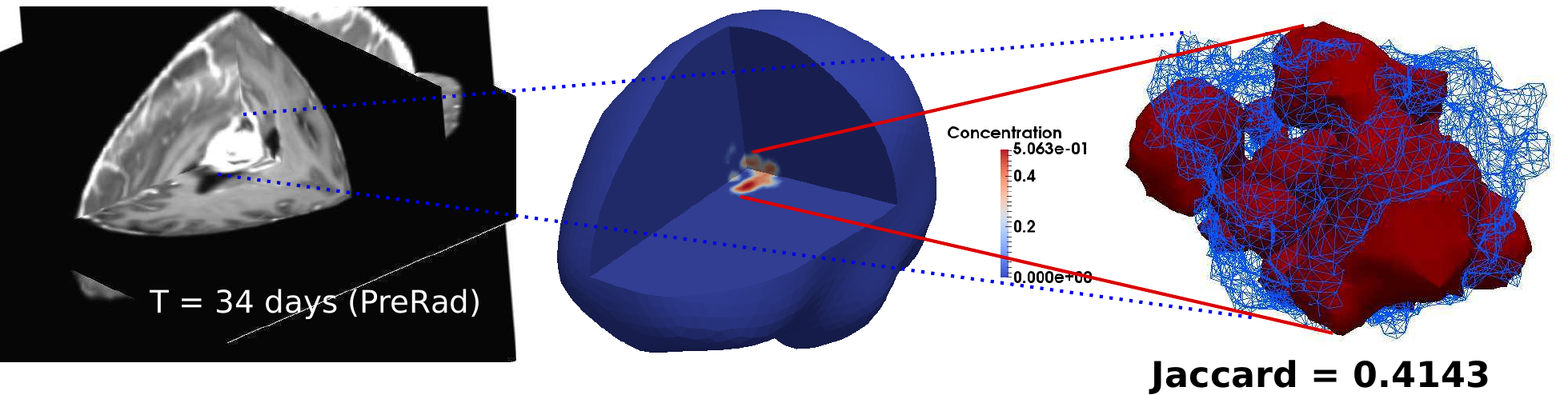}}  \\
  \hline
  \end{tabular}
 \begin{tabular}{ c c c c c c c c c c c }
 \hline
 \parbox[t][][t]{2.2cm}{\centering \textbf{Iteration}\\\textbf{k=1}} & \parbox[t][][t]{1.2cm}{\centering $\mathbf{J}(\mathcal{P}_1)$\\0.2554} &\parbox[t][][t]{1.2cm}{\centering $\mathbf{L}_1$\\0.0002} & \parbox[t][][t]{1.cm}{\centering $\boldsymbol{\nu}_1$\\0.18537} &  \parbox[t][][t]{1.cm}{\centering $\mathbf{k}_{n1}$\\2.1046} & \parbox[t][][t]{1.2cm}{\centering $\mathbf{S}_{n1}$\\10000.06}  & \parbox[t][][t]{1.cm}{\centering $\boldsymbol{\delta}_{n1}$\\8639.94} & \parbox[t][][t]{0.9cm}{\centering $\boldsymbol{\gamma}_1^2$\\0.1225} &  \parbox[t][][t]{0.9cm}{\centering $\mathbf{E}_{1}$\\693.97} &  \parbox[t][][t]{0.9cm}{\centering $\boldsymbol{\delta}_{1}$\\0.2160} &  \parbox[t][][t]{0.9cm}{\centering $\textbf{c}_{e1}$\\0.611}\\  \hline \\
 \textbf{MRI} &  &  & &  \textbf{FOM} &   &  & \textbf{Comparison} &  & & \\ \\
\end{tabular}
\begin{tabular}{l}
\raisebox{-0.5\height}{\includegraphics[scale=0.92]{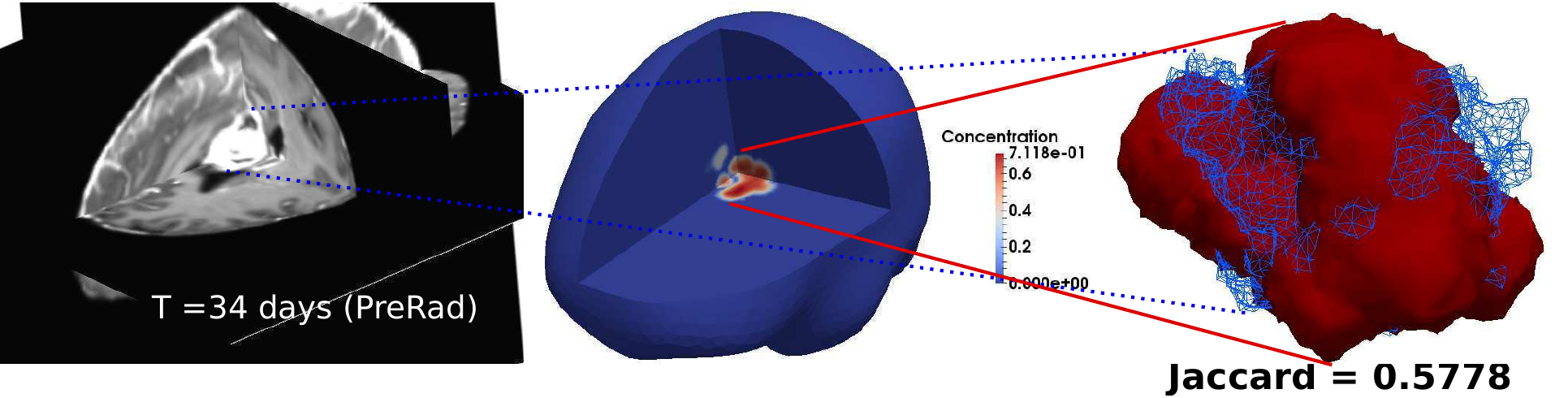}}  \\
  \hline
  \end{tabular}
 \begin{tabular}{ c c c c c c c c c c c }
 \hline
 \parbox[t][][t]{2.2cm}{\centering \textbf{Iteration}\\\textbf{k=2}} & \parbox[t][][t]{1.2cm}{\centering $\mathbf{J}(\mathcal{P}_2)$\\0.2298} &\parbox[t][][t]{1.2cm}{\centering $\mathbf{L}_2$\\0.000532} & \parbox[t][][t]{1.cm}{\centering $\boldsymbol{\nu}_2$\\0.10573} &  \parbox[t][][t]{1.cm}{\centering $\mathbf{k}_{n2}$\\2.0805} & \parbox[t][][t]{1.cm}{\centering $\mathbf{S}_{n2}$\\10000.05}  & \parbox[t][][t]{1.cm}{\centering $\boldsymbol{\delta}_{n2}$\\8639.94} & \parbox[t][][t]{0.9cm}{\centering $\boldsymbol{\gamma}_2^2$\\0.1225} &  \parbox[t][][t]{0.9cm}{\centering $\mathbf{E}_{2}$\\693.97} &  \parbox[t][][t]{0.9cm}{\centering $\boldsymbol{\delta}_{2}$\\0.2376} &  \parbox[t][][t]{0.9cm}{\centering $\textbf{c}_{e2}$\\0.5792}\\  \hline \\
 \textbf{MRI} &  &  & &  \textbf{FOM} &   &  & \textbf{Comparison} &  & & \\ \\
\end{tabular}
\begin{tabular}{l}
\raisebox{-0.5\height}{\includegraphics[scale=0.92]{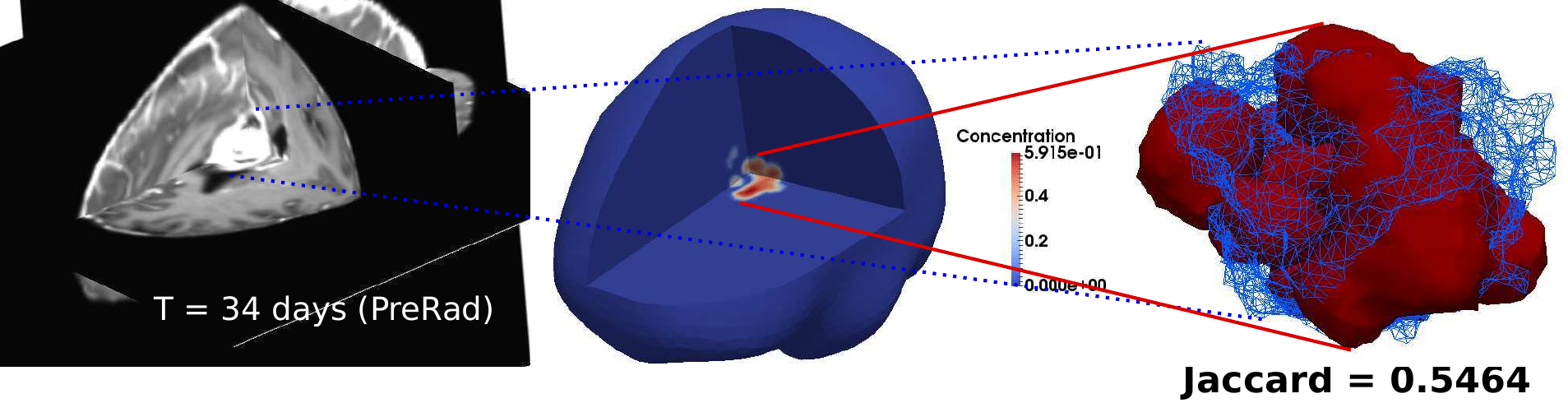}}  \\ \\
  \hline
  \end{tabular}
  \end{minipage} }
  \captionof{figure}[]{Values of $J(\phi_h^N(\mathcal{P}_k),\mathcal{P}_k)$ and of $\mathcal{P}_k$ for different iteration steps $k$ of \textbf{Algorithm 1}, with a comparison between the iso--surfaces $\phi_{\text{data}}(T)=0$ (highlighted in blue color) and $\phi_h^N(\mathcal{P}_k)=\phi_e/2$ (highlighted in red color).}
  \label{fig:10}
\end{table}
}
\newpage
\noindent The \textbf{Algorithm 1} stops since
\[J(\phi_h^N(\mathcal{P}_{k+1}),\mathcal{P}_{k+1})= J(\phi_h^N(\mathcal{P}_{k}),\mathcal{P}_{k})\]
for $k=2$. Indeed, the ROM optimization algorithm for $k=2$ makes no advances in the parameter space (see Figure \ref{fig:13}).
We thus identify 
\begin{align}
\label{eqn:opttest2}
\notag &\mathcal{P}_{\text{opt}}\equiv \mathcal{P}_2=\\
&\notag \{L=0.000532,\nu=0.10573,k_n=2.0805,S_n=10000.05,\delta_n=8639.94,\gamma^2=0.1225,\\
& E=693.97,\delta=0.2376,c_e=0.5792\},
\end{align} 
with the corresponding units. 
We observe that also in this case the overall overlapping between the tumour extensions from FOM simulations and from data is increasing, which is shown by the corresponding increase in the value of the Jaccard index. The Jaccard indexes are anyhow smaller than the values computed for Test Case 1, reported in Figure \ref{fig:4}, due to greater uncertainty in the segmentation of initial left tumour particles after surgery and in the segmentation of the tumour extension before radiotherapy, which consists in peritumoural infiltrations which are not visible in the MRI data.
\subsubsection{Step 2}
\begin{figure}[!h]
\centering
\includegraphics[width=0.8 \linewidth]{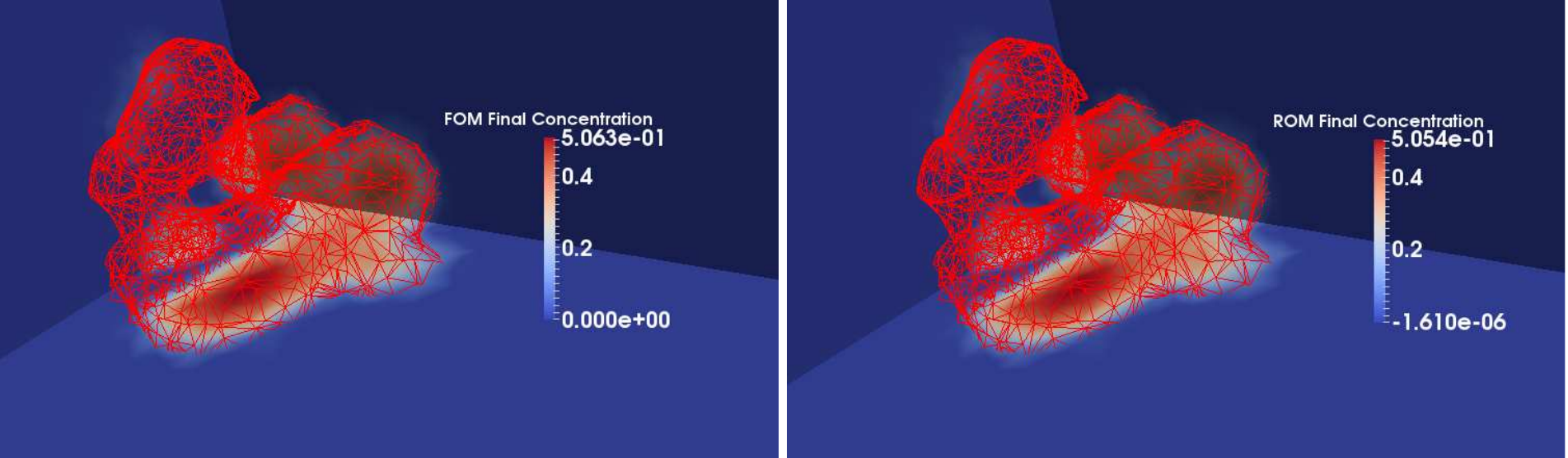}
\caption{Comparison between the final state $\phi_h^N$, solution of the FOM simulation with parameter set $\mathcal{P}_0$ and the corresponding final state $\sum_{i=1}^{N_{\text{POD}}}\alpha_{i0}^N\xi_{i}^{\phi}$, solution of the ROM system \eqref{eqn:6}. The iso--surfaces $\phi_h^N=\phi_e/2$ and $\sum_{i=1}^{N_{\text{POD}}}\alpha_{i0}^N\xi_{i}^{\phi}=\phi_e/2$ are highlighted in red color.
}
\label{fig:12}
\end{figure}

In Table \ref{tab:4} we report, for each step $k$ of \textbf{Algorithm 1}, the values of the cumulated fractions of $\text{tr}\mathbf{F}_1^T\mathbf{F}_1,\text{tr}\mathbf{F}_2^T\mathbf{F}_2,\text{tr}\mathbf{F}_3^T\mathbf{F}_3,\text{tr}(\psi_1'(\mathbf{F}_1))^T\psi_1'(\mathbf{F}_1),\text{tr}(\psi_1''(\mathbf{F}_1))^T\psi_1''(\mathbf{F}_1)$ associated to the eigenvalues of the corresponding matrices, arranging them starting from the eigenvalue with the highest magnitude and following a decreasing order.
\begin{table}[h!]
\caption{Values of the cumulated fraction of $\text{tr}\mathbf{F}_1^T\mathbf{F}_1$,$\text{tr}\mathbf{F}_2^T\mathbf{F}_2$,$\text{tr}\mathbf{F}_3^T\mathbf{F}_3$,$\text{tr}(\psi_1'(\mathbf{F}_1))^T\psi_1'(\mathbf{F}_1)$,$\text{tr}(\psi_1''(\mathbf{F}_1))^T\psi_1''(\mathbf{F}_1)$ for the first eigenvalues with the highest magnitude.}
\label{tab:4}    
\begin{center}
\begin{tabular}{lllllll}
\hline
\parbox[t][][t]{2.2cm}{\centering \textbf{Iteration}\\\textbf{k=0}} & Eigenvalue & \% $\text{tr}\mathbf{F}_1^T\mathbf{F}_1$ & \% $\text{tr}\mathbf{F}_2^T\mathbf{F}_2$ & \% $\text{tr}\mathbf{F}_3^T\mathbf{F}_3$ & \% $\text{tr}(\psi_1')^T\psi_1'$ & \% $\text{tr}(\psi_1'')^T\psi_1''$\\
\hline
& First & $93.5098$ & $99.9731$ & $99.9904$ & $99.9714$ & $99.7753$ \\ 
& Second & $99.5129$ & $99.9973$ & $99.9993$ & $99.9960$ & $99.9577$\\
& Third & $99.9574$ & $99.9998$ & $99.9999$ & $99.9998$ & $99.9984$\\ 
& Fourth & $99.9957$ & $99.9999$ & $99.9999$ & $99.9999$ & $99.9998$\\ 
 \hline
 \parbox[t][][t]{2.2cm}{\centering \textbf{Iteration}\\\textbf{k=1}} & Eigenvalue & \% $\text{tr}\mathbf{F}_1^T\mathbf{F}_1$ & \% $\text{tr}\mathbf{F}_2^T\mathbf{F}_2$ & \% $\text{tr}\mathbf{F}_3^T\mathbf{F}_3$ & \% $\text{tr}(\psi_1')^T\psi_1'$ & \% $\text{tr}(\psi_1'')^T\psi_1''$\\
\hline
& First & $92.6867$ & $99.7216$ & $99.9210$ & $99.5998$ & $99.3884$ \\ 
& Second & $98.8774$ & $99.9860$ & $99.9951$ & $99.9878$ & $99.9182$\\
& Third & $99.8253$ & $99.9980$ & $99.9993$ & $99.9980$ & $99.9844$\\ 
& Fourth & $99.9702$ & $99.9996$ & $99.9998$ & $99.9997$ & $99.9983$\\ 
& Fifth & $99.9945$ & $99.9999$ & $99.9999$ & $99.9999$ & $99.9997$\\
 \hline
  \parbox[t][][t]{2.2cm}{\centering \textbf{Iteration}\\\textbf{k=2}} & Eigenvalue & \% $\text{tr}\mathbf{F}_1^T\mathbf{F}_1$ & \% $\text{tr}\mathbf{F}_2^T\mathbf{F}_2$ & \% $\text{tr}\mathbf{F}_3^T\mathbf{F}_3$ & \% $\text{tr}(\psi_1')^T\psi_1'$ & \% $\text{tr}(\psi_1'')^T\psi_1''$\\
\hline
& First & $93.2361$ & $99.9348$ & $99.9793$ & $99.9325$ & $99.4953$ \\ 
& Second & $99.3359$ & $99.9958$ & $99.9988$ & $99.9946$ & $99.9383$\\
& Third & $99.9291$ & $99.9995$ & $99.9998$ & $99.9996$ & $99.9971$\\ 
& Fourth & $99.9915$ & $99.9999$ & $99.9999$ & $99.9999$ & $99.9996$\\ 
 \hline
\end{tabular}
\end{center}
\end{table}\\
We thus have that
\[
N_{\text{POD}}=N_{\phi}^{\text{POD}}=4, \quad \text{for} \; k=0,2, \; N_{\text{POD}}=N_{\phi}^{\text{POD}}=5, \quad \text{for} \; k=1.
\]
In Figure \ref{fig:11} we show the basis elements $\xi_{i}^{\phi}$, corresponding to the highest eigenvalues needed to explain the $99.99\%$ variance of the data, for $k=0,1$.
\begin{table}[ht!]
\resizebox{0.9\textwidth}{!}{\begin{minipage}{\textwidth}
\begin{tabular}{ c }
 \hline
 \parbox[t][][t]{14.2cm}{\centering \textbf{Iteration k=0}} \\  \hline \\
\end{tabular}\\
\begin{tabular}{l}
\raisebox{-0.5\height}{\includegraphics[scale=0.9]{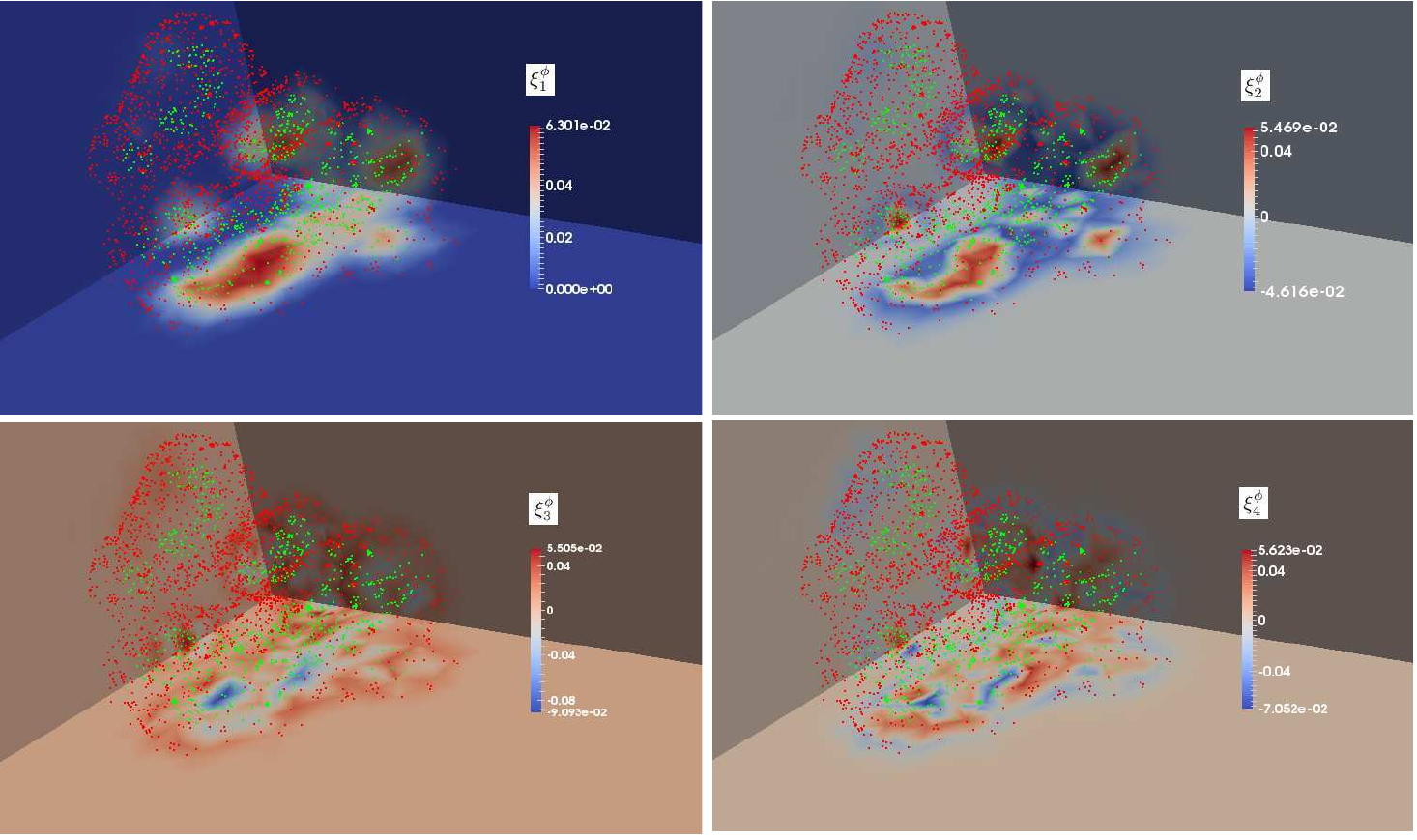}}  \\
  \hline
  \end{tabular}\\
  \begin{tabular}{ c }
 \hline
 \parbox[t][][t]{14.2cm}{\centering \textbf{Iteration k=1}} \\  \hline \\
\end{tabular}\\
\begin{tabular}{l}
\raisebox{-0.5\height}{\includegraphics[scale=0.9]{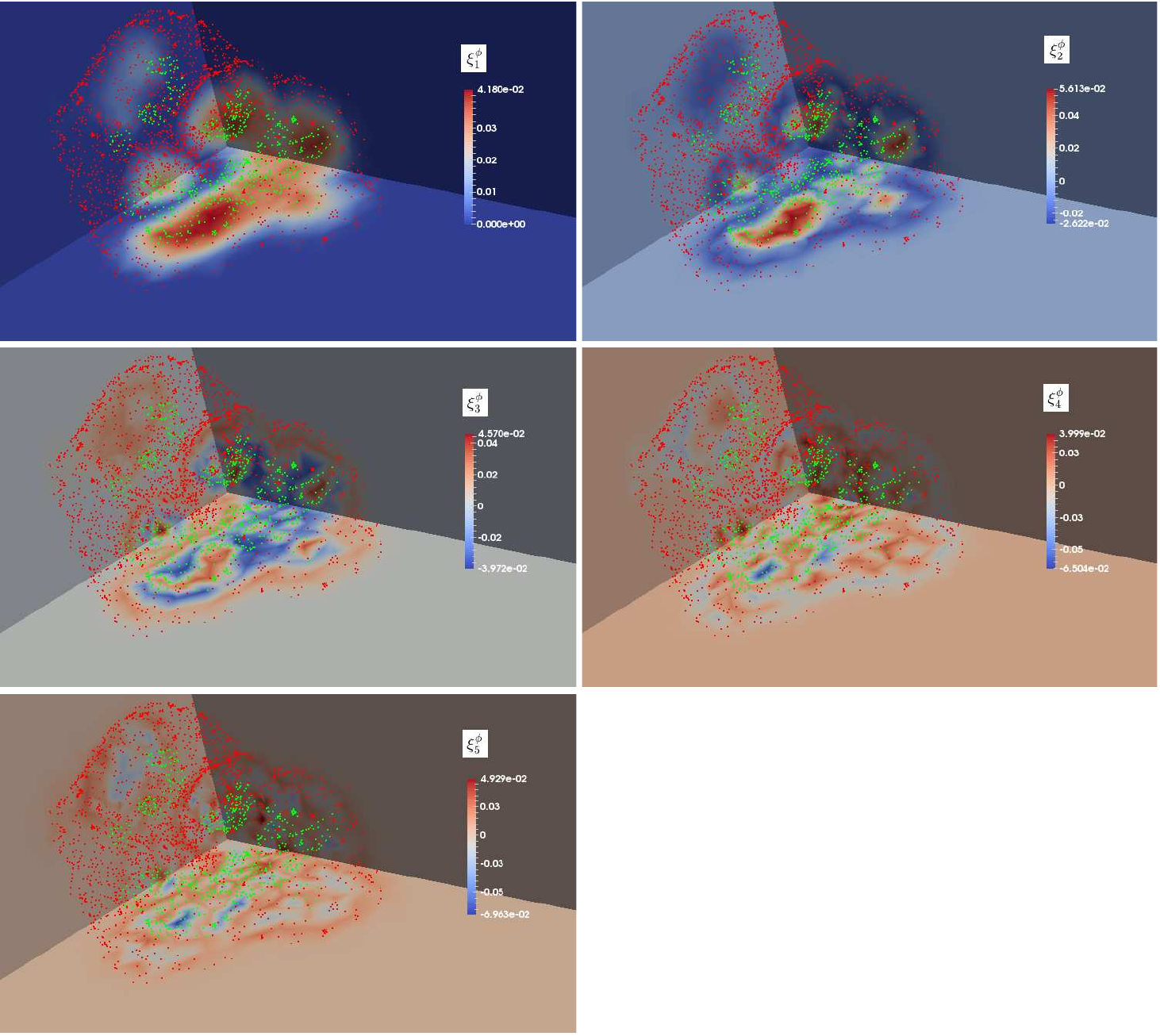}}  \\
  \hline
  \end{tabular}
  \end{minipage} }
  \captionof{figure}[]{Plot of the basis elements $\xi_{i}^{\phi}$, corresponding to the highest eigenvalues needed to explain the $99.99\%$ of the data, for $k=0,1$. Green and red points are
distributed over the initial condition and final distribution of cell concentration, respectively.}
  \label{fig:11}
  \end{table}
  \\ \\ \\ \\ \\ \\ \\ \\ \\ \\ \\ \\ \\
\noindent
We observe, as in the previous Test Case $1$, that $\xi_{1}^{\phi}$ and $\xi_{2}^{\phi}$ are distributed over the core of the final state $\phi_h^N$ and the initial condition $\phi_h^0$ respectively, whereas $\xi_{3}^{\phi}$, $\xi_{4}^{\phi}$ and $\xi_{5}^{\phi}$ are oscillating functions over the set where the tumour is expanding during its temporal evolution, and thus contain the information about the tumour boundary and its expansion. We observe that the low dimensionality of the ROM systems is preserved also in the case of tumour dynamics with sparse tumour particles and infiltrations, if the mesh for the FOM system is sufficiently well refined in the region of the tumour core and infiltrations.\\
\noindent We finally show in Figure \ref{fig:12} a comparison between the final state $\phi_h^N$ calculated from the FOM simulation through \textbf{Algorithm 2} with parameter set $\mathcal{P}_0$ and the corresponding final state $\sum_{i=1}^{N_{\text{POD}}}\alpha_{i0}^N\xi_{i}^{\phi}$ obtained as a solution of the ROM system \eqref{eqn:6} through \textbf{Algorithm 4}. 

\noindent We observe also in this test case that the ROM solution is approximating the FOM solution with a very high fidelity.

\subsubsection{Steps 3 and 4}
In Figure \ref{fig:13} we report the values of the functional $J(\vec{\alpha}_l,\mathcal{P}_l)$, of the normalised set of parameters 
\[
\mathcal{P}_l/\mathcal{P}_{\text{exp}}=\{L_l/L_{\text{exp}},\nu_l/\nu_{\text{exp}},k_{nl}/k_{n\text{exp}},S_{nl}/S_{n\text{exp}},\delta_{nl}/\delta_{n\text{exp}},\gamma_l^2/\gamma_{\text{exp}}^2,E_l/E_{\text{exp}},\delta_c/\delta_{c\text{exp}},c_e/c_{e\text{exp}}\},
\]
and of $|\mathcal{P}_l(1)-\mathcal{P}_l|$, computed in \textbf{Steps 3} and \textbf{4}  of \textbf{Algorithm 1}, for $k=0,1,2$. We also plot the isosurfaces $\phi_{\text{data}}(T)=0$ from the MRI data and $\sum_{i=1}^{N_{\text{POD}}}\alpha_{i\bar{l}}^N\xi_{i}^{\phi}=\phi_e/2$ from the ROM simulations, where $\bar{l}$ is the value of the last iteration of \textbf{Step 4}, reporting the value of the Jaccard index between the two sets enclosed by these surfaces.

\begin{table}[p]
\resizebox{0.9\textwidth}{!}{\begin{minipage}{\textwidth}
\begin{tabular}{ c c c }
 \hline
 \parbox[t][][t]{5.2cm}{\centering \textbf{Iteration k=0}} & &\\ \hline \\ 
 \parbox[t][][t]{5.6cm}{\centering $\mathbf{J}(\vec{\alpha}_l,\mathcal{P}_l)$} &\parbox[t][][t]{6.6cm}{\centering $\mathcal{P}_l/\mathcal{P}_{\text{exp}}$}& \\  
  \raisebox{-0.5\height}{\includegraphics[scale=0.23]{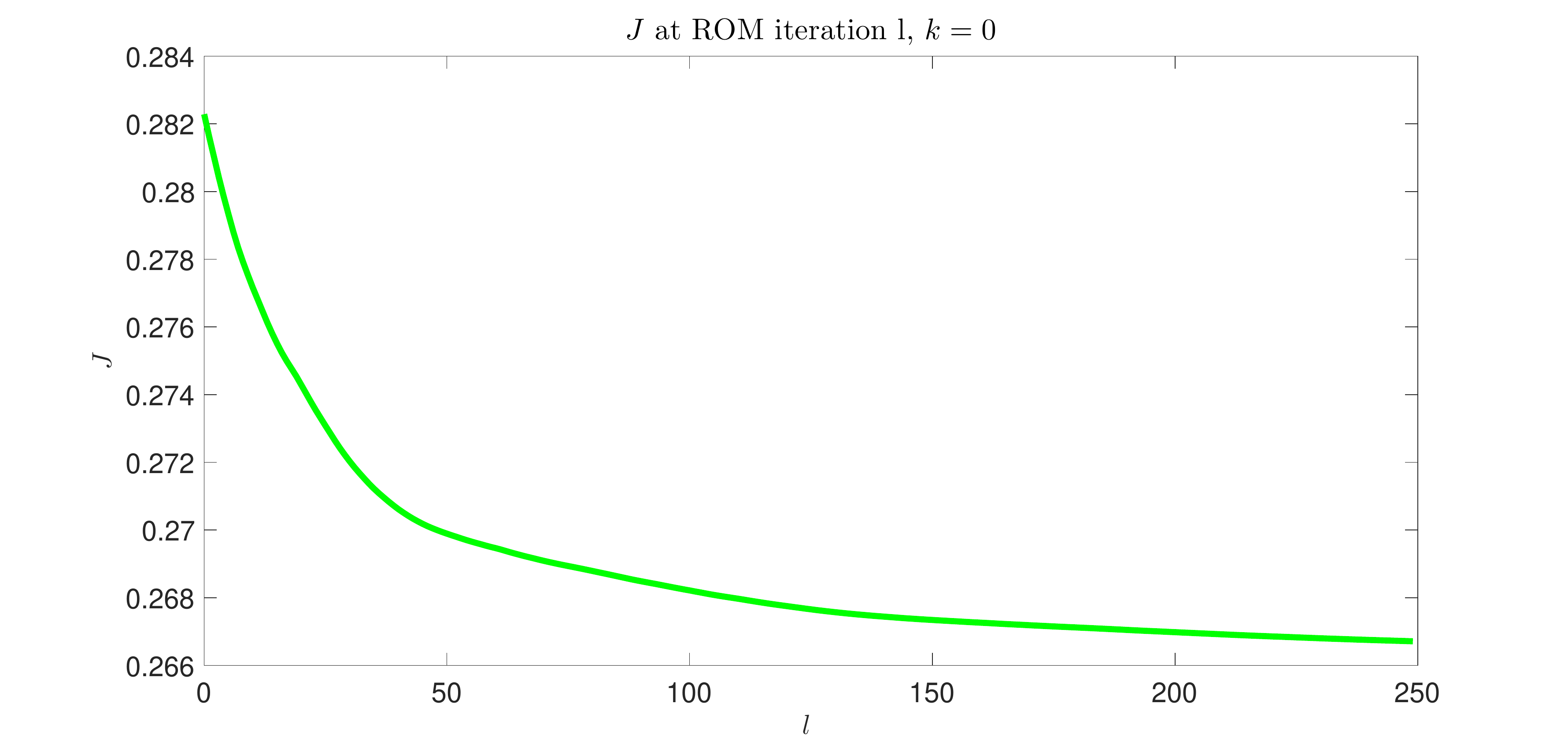}}&\raisebox{-0.5\height}{\includegraphics[scale=0.23]{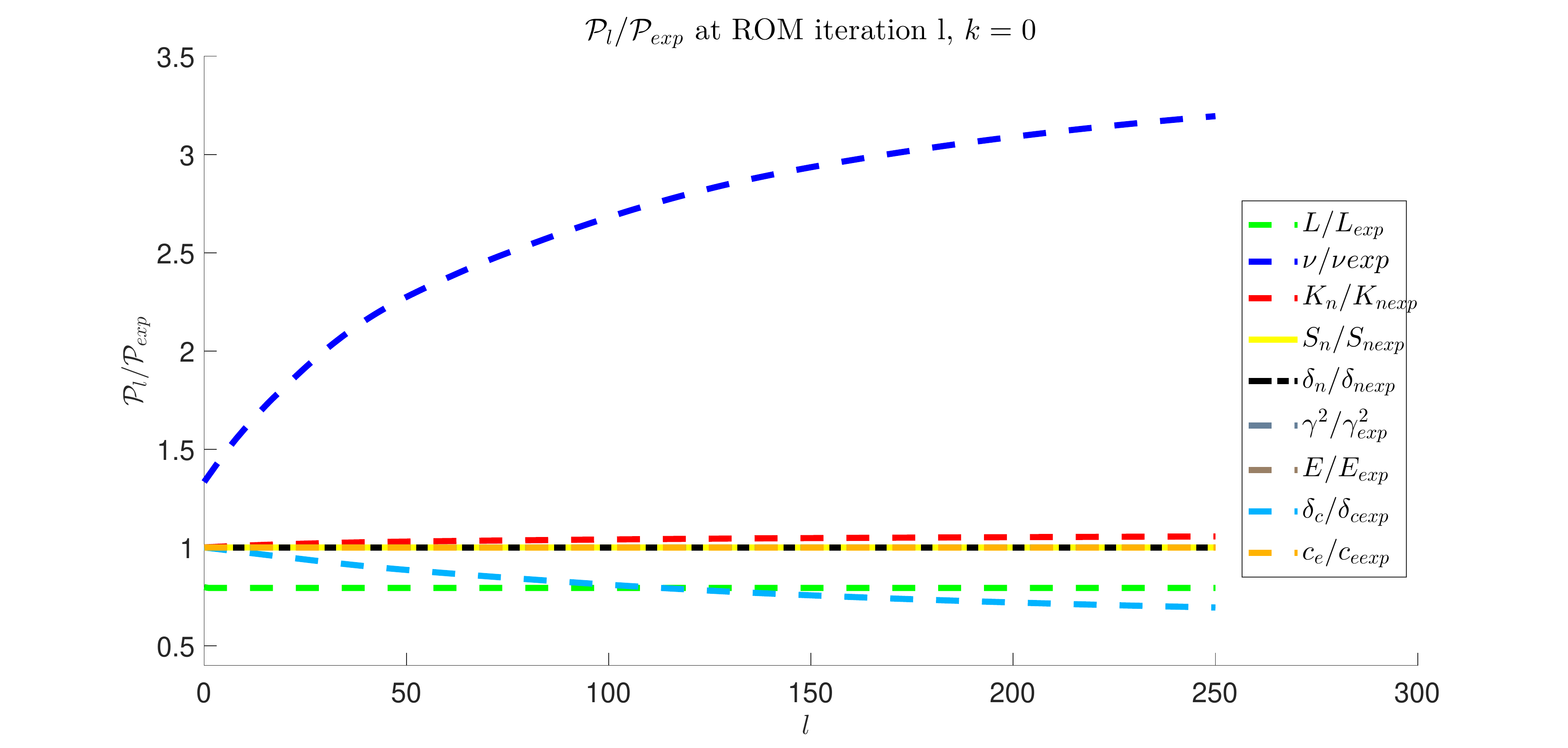}}&\\ \\
\parbox[t][][t]{5.6cm}{\centering $|\mathcal{P}_l(1)-\mathcal{P}_l|$} & & \\  
  \raisebox{-0.5\height}{\includegraphics[scale=0.23]{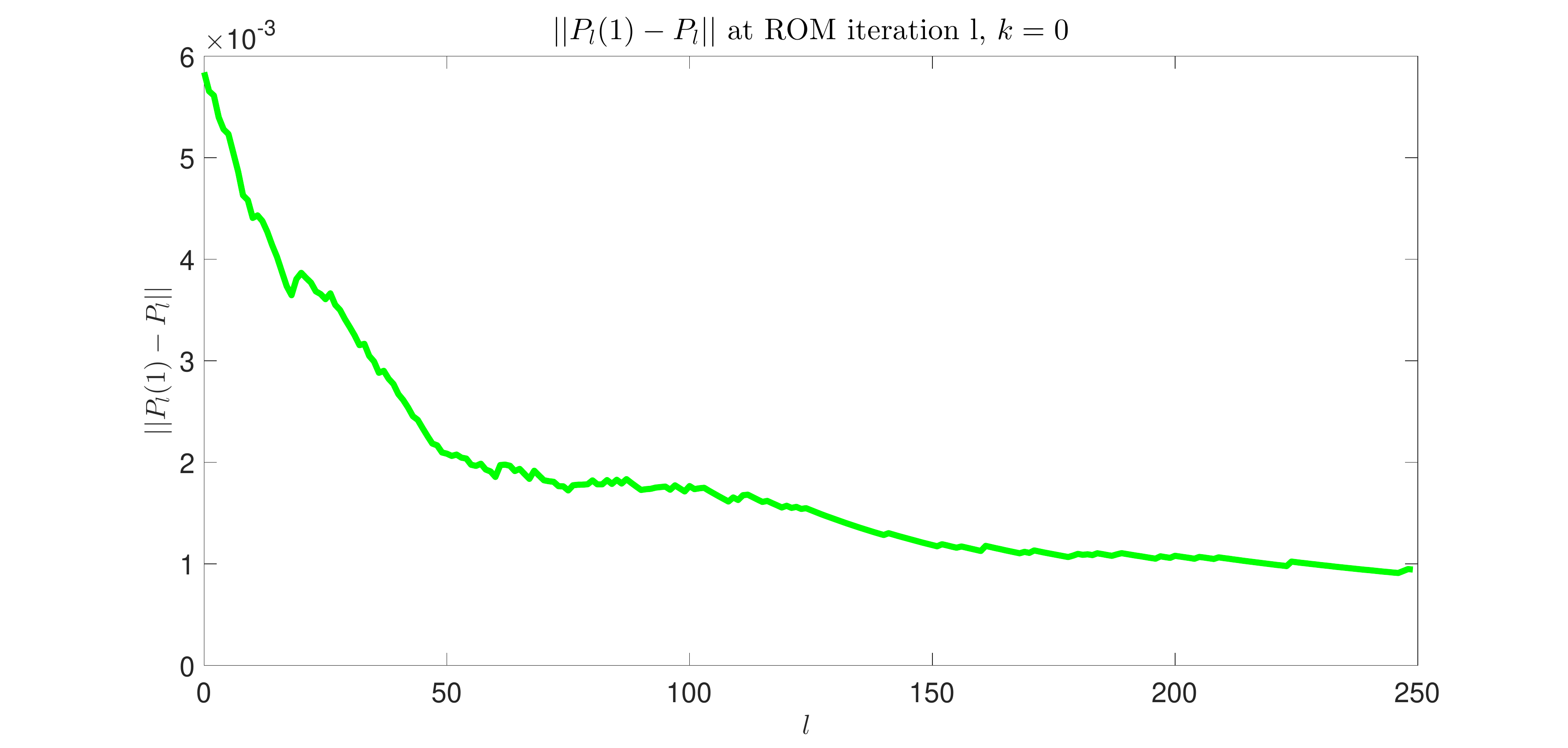}}& &\\ \\  
 \end{tabular}
 \begin{tabular}{c c c}
 \parbox[t][][t]{3.2cm}{\centering\textbf{MRI}} &   \parbox[t][][t]{7.2cm}{\centering\textbf{ROM ($\mathcal{P}_{\bar{l}=249})$}}  &  \parbox[t][][t]{4.2cm}{\centering\textbf{Comparison}}  \\ \\
\end{tabular}\\
\begin{tabular}{l}
\raisebox{-0.5\height}{\includegraphics[scale=0.8]{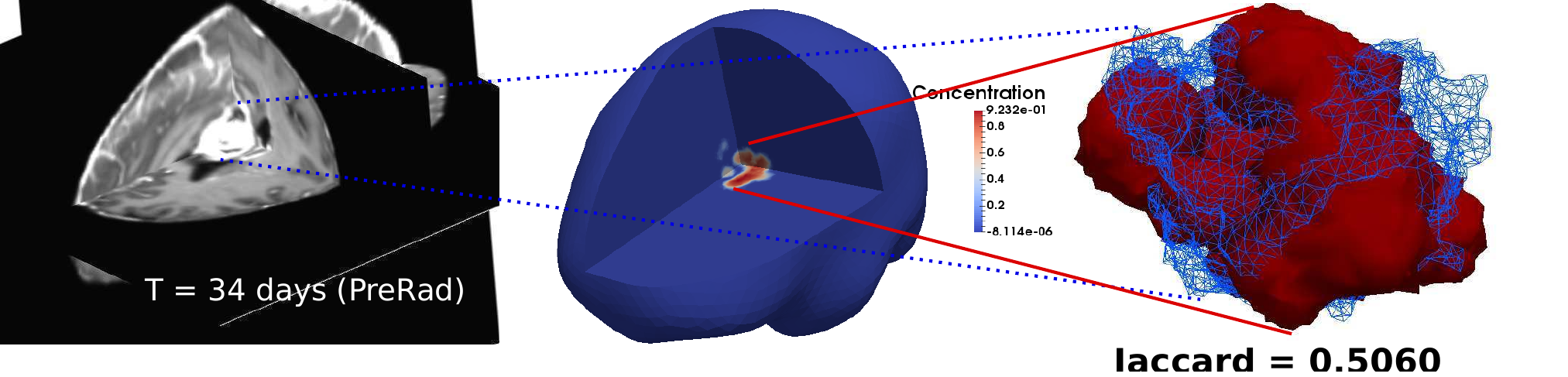}}  \\
  \end{tabular}
  \begin{tabular}{ c c c }
 \hline
 \parbox[t][][t]{5.2cm}{\centering \textbf{Iteration k=1}} & &\\ \hline \\ 
 \parbox[t][][t]{5.6cm}{\centering $\mathbf{J}(\vec{\alpha}_l,\mathcal{P}_l)$} &\parbox[t][][t]{6.6cm}{\centering $\mathcal{P}_l/\mathcal{P}_0$} &\\  \
  \raisebox{-0.5\height}{\includegraphics[scale=0.23]{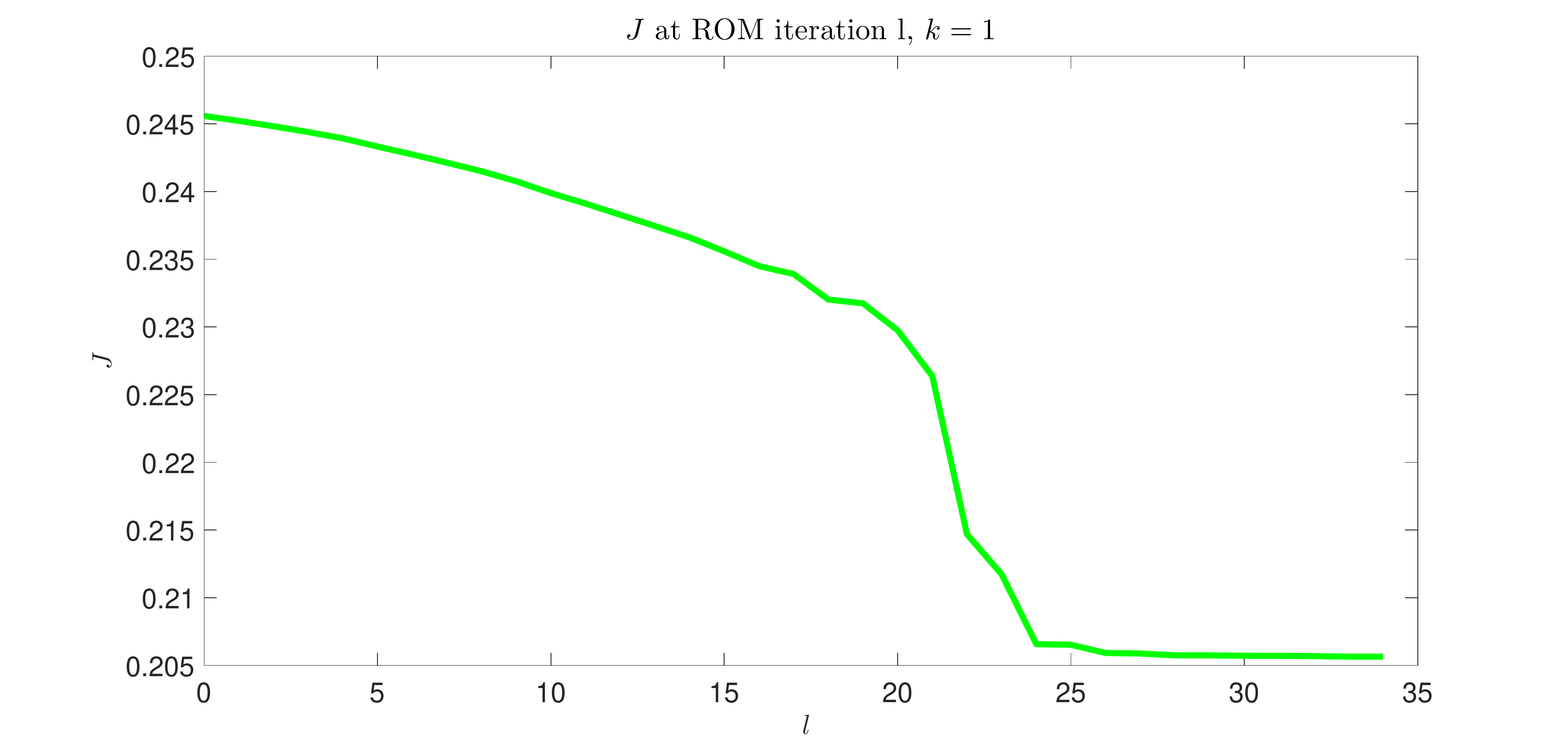}}&\raisebox{-0.5\height}{\includegraphics[scale=0.23]{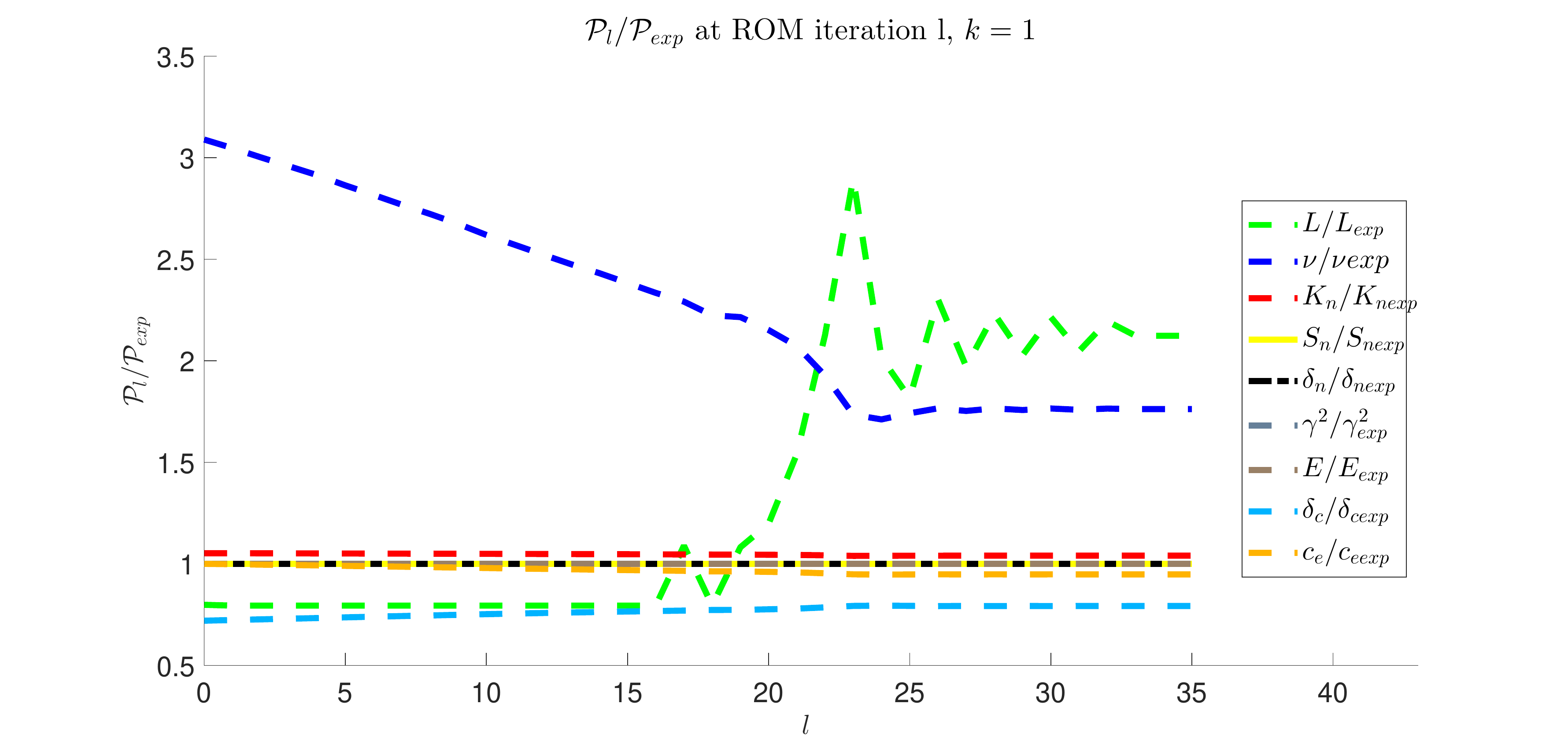}}&\\ \\
\parbox[t][][t]{5.6cm}{\centering $|\mathcal{P}_l(1)-\mathcal{P}_l|$} & & \\  
  \raisebox{-0.5\height}{\includegraphics[scale=0.23]{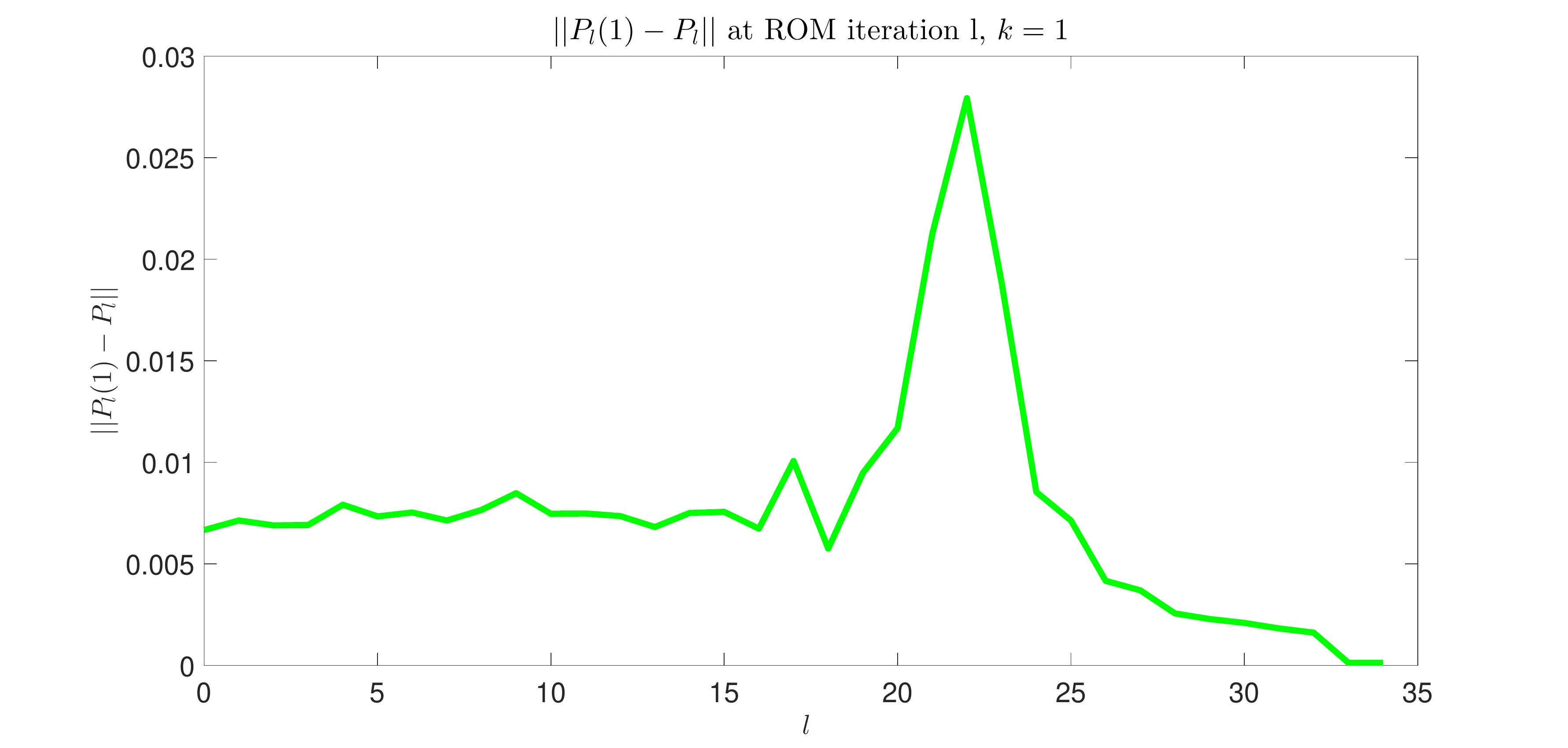}}& &\\ \\  
 \end{tabular}\\
 \end{minipage} }
\end{table}
\newpage
\begin{table}[h!]
\resizebox{0.98\textwidth}{!}{\begin{minipage}{\textwidth}
 \begin{tabular}{c c c}
 \parbox[t][][t]{3.2cm}{\centering\textbf{MRI}} &   \parbox[t][][t]{7.2cm}{\centering\textbf{ROM ($\mathcal{P}_{\bar{l}=34})$}}  &  \parbox[t][][t]{4.2cm}{\centering\textbf{Comparison}}  \\ \\
\end{tabular}
\begin{tabular}{l}
\raisebox{-0.5\height}{\includegraphics[scale=0.8]{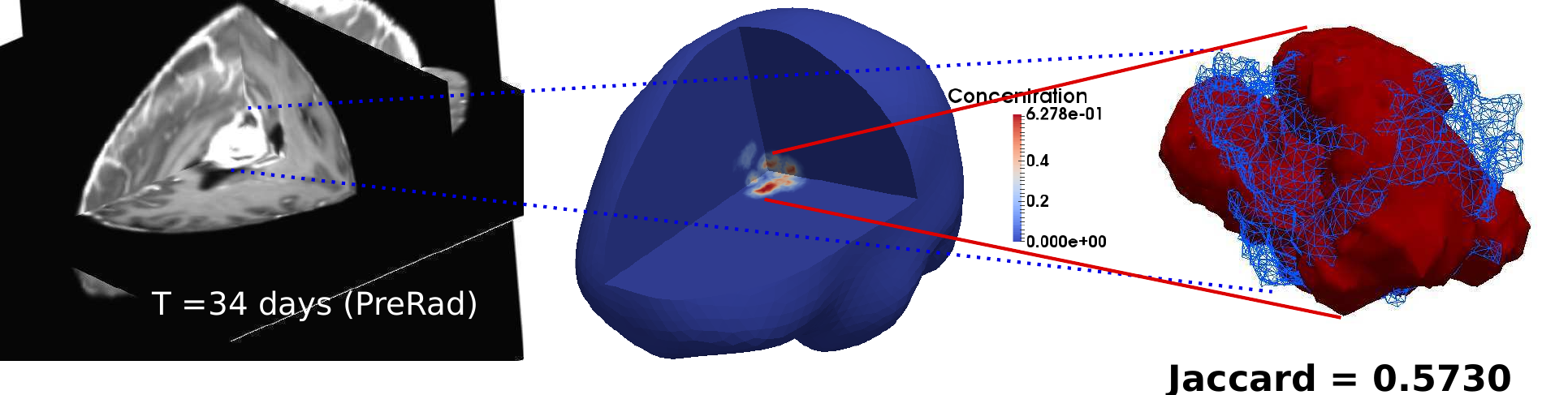}}  \\
  \end{tabular}
   \begin{tabular}{ c c c }
 \hline
 \parbox[t][][t]{5.2cm}{\centering \textbf{Iteration k=2}} & &\\ \hline \\ 
 \parbox[t][][t]{5.6cm}{\centering $\mathbf{J}(\vec{\alpha}_l,\mathcal{P}_l)$} &\parbox[t][][t]{6.6cm}{\centering $|\mathcal{P}_l(1)-\mathcal{P}_l|$} &\\  \
  \raisebox{-0.5\height}{\includegraphics[scale=0.23]{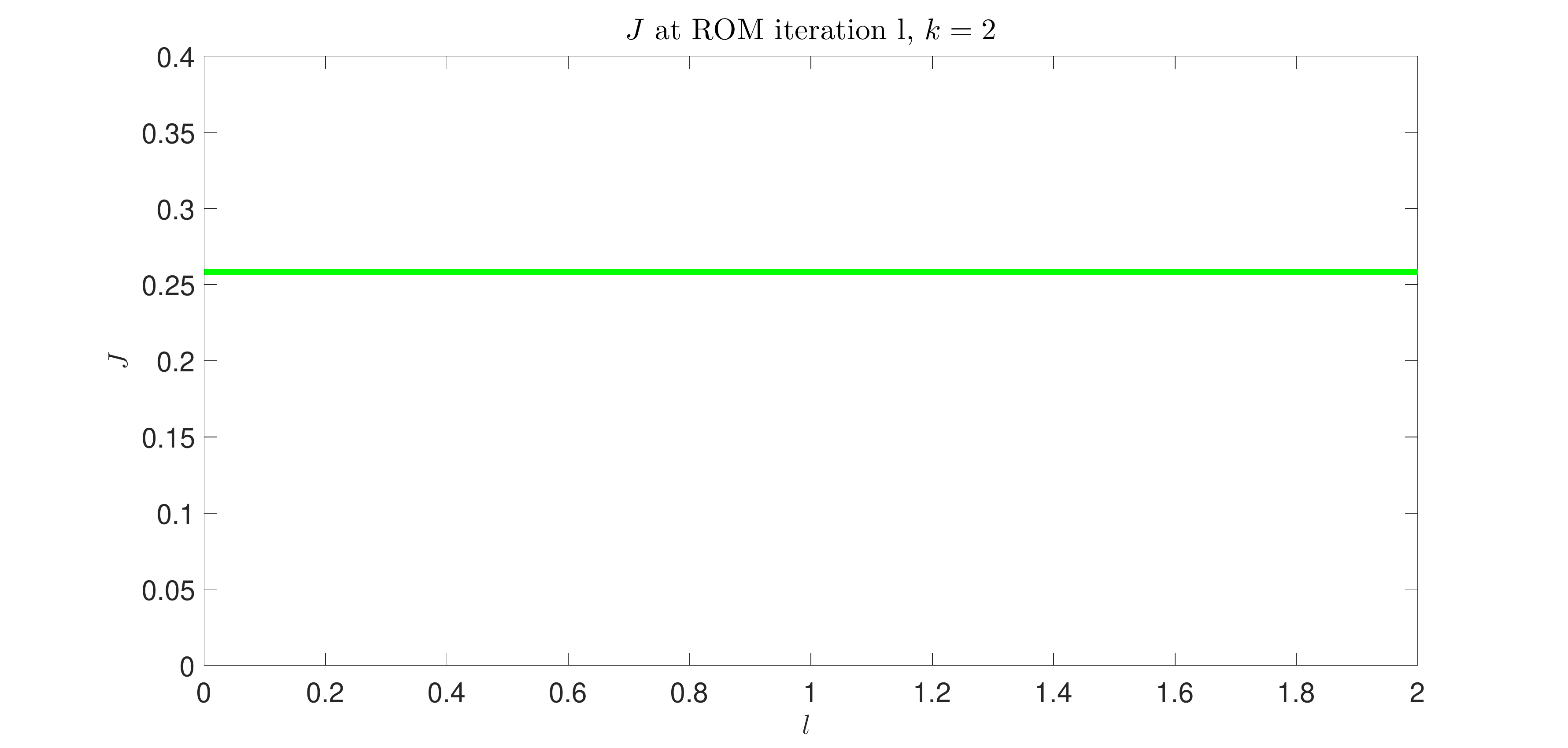}}&\raisebox{-0.5\height}{\includegraphics[scale=0.23]{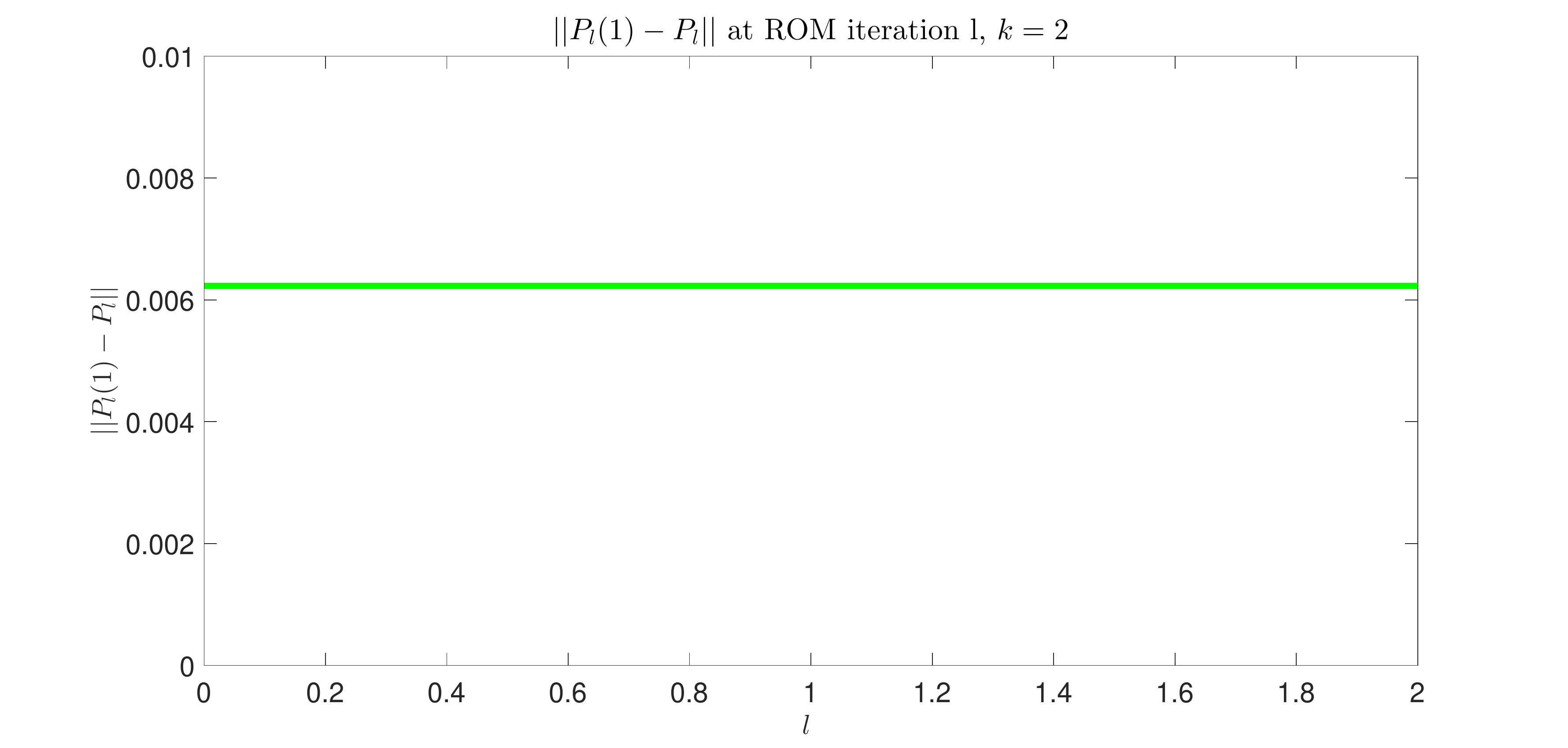}}&\\ \\ 
 \end{tabular}\\
  \end{minipage} }
 \captionof{figure}[]{Values of the functional $J(\vec{\alpha}_l,\mathcal{P}_l)$ and of the normalised set of parameters $\mathcal{P}_l/\mathcal{P}_0$ for steps $k=0,1,2$ of \textbf{Algorithm 1}, together with a comparison between the iso--surfaces $\phi_{\text{data}}(T)=0$(highlighted in blue color) and $\sum_{i=1}^{N_{\text{POD}}}\alpha_{i\bar{l}}^N\xi_{i}^{\phi}=\phi_e/2$ (highlighted in red color).}
  \label{fig:13}
\end{table}
For $k=0$ the ROM optimization process in \textbf{Steps 3} and \textbf{4} goes through $249$ steps before matching the termination conditions. Thus the ROM Optimization Algorithm is much slower in Test Case $2$, corresponding to irregular initial data and target function, then in Test Case $1$, which is characterised by more regular data.
We also observe that, like in Test Case $1$, only the model parameters $L,\nu,\delta,c_e$ change significantly from their starting values during the optimization process, being the system quite insensitive to changes of the remaining parameters $k_n,S_n,\delta_n,\gamma^2,E$. The proliferation rate $\nu$ is the most sensitive parameter for $k=0$.
For $k=1$ the ROM optimization process in \textbf{Steps 3} and \textbf{4} goes through $34$ steps. We observe that for $k=1$ also the parameter $L$ varies by a large excursion along the functional minimisation and relaxes in an oscillating manner onto its local equilibrium value.  Finally, for $k=2$ the ROM optimization process is making no progress.

Comparing the output of the Optimization Algorithm \eqref{alg:opt} with the ones reported in Figure \ref{fig:7} for Test Case $1$, we observe that in Test Case $2$ we need more iterations to converge to an optimal state, that the functional overall decreases of a much smaller amount (due to the uncertainty in identifying the target function) and that the computation of new POD basis associated to different set of parameters introduces a higher degree of variability in the tumour dynamics, causing the minimising sequences of parameters to show oscillations instead of monotone relaxing to the optimal state like in Test Case $1$. 
\subsubsection{Computational cost}
We conclude by reporting in Table \ref{tab:5} the elapsed time (in seconds) for the computation of steps $1-4$ of \textbf{Algorithm 1} for the different values of the iteration step $k$.
\begin{table}[h!]
\caption{Elapsed time (in seconds) for the computation of steps $1-4$, for the different values of the iteration step $k$.}
\label{tab:5}    
\begin{center}
\begin{tabular}{lllll}
\hline
\parbox[t][][t]{2.2cm}{\centering \textbf{Iteration}\\\textbf{k=0}} & \parbox[t][][t]{2.2cm}{\centering \textbf{Step 1}\\$144575.36$} & \parbox[t][][t]{2.2cm}{\centering \textbf{Step 2}\\$8.58$} & \parbox[t][][t]{2.2cm}{\centering \textbf{Step 3}\\$2438.42$} & \parbox[t][][t]{2.2cm}{\centering \textbf{Step 4}\\$3431.52$}\\ \\
\hline
 \parbox[t][][t]{2.2cm}{\centering \textbf{Iteration}\\\textbf{k=1}} & \parbox[t][][t]{2.2cm}{\centering \textbf{Step 1}\\$164645.5$} & \parbox[t][][t]{2.2cm}{\centering \textbf{Step 2}\\$11.75$} & \parbox[t][][t]{2.2cm}{\centering \textbf{Step 3}\\$6453.26$} & \parbox[t][][t]{2.2cm}{\centering \textbf{Step 4}\\$336.09$}\\ \\  
 \hline
 \parbox[t][][t]{2.2cm}{\centering \textbf{Iteration}\\\textbf{k=2}} & \parbox[t][][t]{2.2cm}{\centering \textbf{Step 1}\\$135222.43$} & \parbox[t][][t]{2.2cm}{\centering \textbf{Step 2}\\$7.75$} & \parbox[t][][t]{2.2cm}{\centering \textbf{Step 3}\\$2416.61$} & \parbox[t][][t]{2.2cm}{\centering \textbf{Step 4}\\$45.14$}\\ \\
\hline
\end{tabular}
\end{center}
\end{table}
Comparing Table \ref{tab:5} with Table \ref{tab:3} we observe that in the case of tumour recurrence with sparse particles and infiltrations the computational time for the projected gradient iterations at the ROM level is $2$ to $3$ order of magnitude smaller than the time needed to solve the FOM problem, provided to properly refine the mesh in the FOM simulations. We also note that the FOM computations require a comparable computational time with respect to Test Case $1$. The same is valid for the computational time required to assemble the ROM systems in \textbf{Step 3}.\\
We finally conclude that the computational efficiency of the Optimization Algorithm \eqref{alg:opt} is unaffected by the degree of regularity of the tumour dynamics, at least in the test cases analysed here where no morphological transition happens during the evolution. The degree of convergence of the optimization algorithm and the degree of variability of parameters along the projected gradient directions introduced by exploring the parameter space through different basis functions is instead affected by the tumour dynamics regularity.

\subsubsection{Some remarks on benchmark results}
In this paragraph we report some numerical results to show  how the POD analysis in \textbf{Step 2} of Algorithm \eqref{alg:opt} varies when the tumour concentration in the FOM simulations in \textbf{Step 1} is spreading on a larger region than the one observed in Test Case $1$, considering a tumour dynamics over a larger time interval. Moreover, we report numerical results about the performance of \textbf{Step 4} of Algorithm \eqref{alg:opt} when a lower threshold of POD significance (namely $99.9\%$) is considered. \\
In Table \ref{tab:2bisbis} we report the POD analysis of the snapshot matrices obtained from the FOM solution \eqref{eqn:fom} at $k=0$ of \textbf{Algorithm 1} with $N=980$ and initial set $\mathcal{P}_0$, i.e. when the tumour dynamics span a time interval of $120$ days. 
\begin{table}[h!]
\caption{Values of the cumulated fraction of $\text{tr}\mathbf{F}_1^T\mathbf{F}_1$,$\text{tr}\mathbf{F}_2^T\mathbf{F}_2$,$\text{tr}\mathbf{F}_3^T\mathbf{F}_3$,$\text{tr}(\psi_1'(\mathbf{F}_1))^T\psi_1'(\mathbf{F}_1)$,$\text{tr}(\psi_1''(\mathbf{F}_1))^T\psi_1''(\mathbf{F}_1)$ for the first eigenvalues with the highest magnitude.}
\label{tab:2bisbis}    
\begin{center}
\begin{tabular}{lllllll}
\hline
\parbox[t][][t]{2.2cm}{\centering \textbf{Iteration}\\\textbf{k=0}\\$N=980$, $\mathcal{P}_0$} & Eigenvalue & \% $\text{tr}\mathbf{F}_1^T\mathbf{F}_1$ & \% $\text{tr}\mathbf{F}_2^T\mathbf{F}_2$ & \% $\text{tr}\mathbf{F}_3^T\mathbf{F}_3$ & \% $\text{tr}(\psi_1')^T\psi_1'$ & \% $\text{tr}(\psi_1'')^T\psi_1''$\\
\hline
& First & $89.4596$ & $99.9168$ & $99.9707$ & $99.9142$ & $97.3884$ \\ 
& Second & $98.1194$ & $99.9783$ & $99.9937$ & $99.9737$ & $99.7634$\\
& Third & $99.6551$ & $99.9931$ & $99.9973$ & $99.9925$ & $99.9445$\\ 
& Fourth & $99.9223$ & $99.9991$ & $99.9996$ & $99.9989$ & $99.9914$\\ 
& Fifth & $99.9818$ & $99.9998$ & $99.9999$ & $99.9998$ & $99.9985$\\
& Sixth & $99.9954$ & $99.9999$ & $99.9999$ & $99.9999$ & $99.9995$\\
 \hline
\end{tabular}
\end{center}
\end{table}\\
We thus have that $N_{\text{POD}}=N_{\phi}^{\text{POD}}=6$ for $k=0$. We note that in this case the number of basis functions needed to explain $99.99\%$ of the data on a time window of $120$ days is lower than in Test Case $1$ (see Table \ref{tab:2bis}). This is due to the fact that the tumour expansion through time in Test Case $2$ is limited by the presence of the ventricle' walls, and thus the region where the tumour is spreading from the initial tumour distribution is contained.\\ 
Finally, we consider the results of the \textbf{Steps 2--4} of Algorithm \eqref{alg:opt}, for the first step $k=0$, with $N=490$ and starting set $\mathcal{P}_0$ and when a threshold value of $99.9\%$ is considered in the POD analysis. In Figure \ref{fig:7tristris} we report the values of the functional $J(\vec{\alpha}_l,\mathcal{P}_l)$, of the normalised set of parameters 
\[
\mathcal{P}_l/\mathcal{P}_{\text{exp}}=\{L_l/L_{\text{exp}},\nu_l/\nu_{\text{exp}},k_{nl}/k_{n\text{exp}},S_{nl}/S_{n\text{exp}},\delta_{nl}/\delta_{n\text{exp}},\gamma_l^2/\gamma_{\text{exp}}^2,E_l/E_{\text{exp}},\delta_c/\delta_{c\text{exp}},c_e/c_{e\text{exp}}\},
\]
and of $|\mathcal{P}_l(1)-\mathcal{P}_l|$, computed in \textbf{Steps 3} and \textbf{4}  of \textbf{Algorithm 1}, for $k=0$. We note that in this case we need $N_{\text{POD}}=3$ to explain the $99.9\%$ of variance of the data. 

\begin{table}[h!]
\resizebox{0.9\textwidth}{!}{\begin{minipage}{\textwidth}
\begin{tabular}{ c c c }
 \hline
 \parbox[t][][t]{5.2cm}{\centering \textbf{Iteration k=0}, $99.9\%$ POD threshold, $N=490$, $\mathcal{P}_0$} & &\\ \hline \\ 
 \parbox[t][][t]{5.6cm}{\centering $\mathbf{J}(\vec{\alpha}_l,\mathcal{P}_l)$} &\parbox[t][][t]{6.6cm}{\centering $\mathcal{P}_l/\mathcal{P}_{\text{exp}}$}& \\  
  \raisebox{-0.5\height}{\includegraphics[scale=0.23]{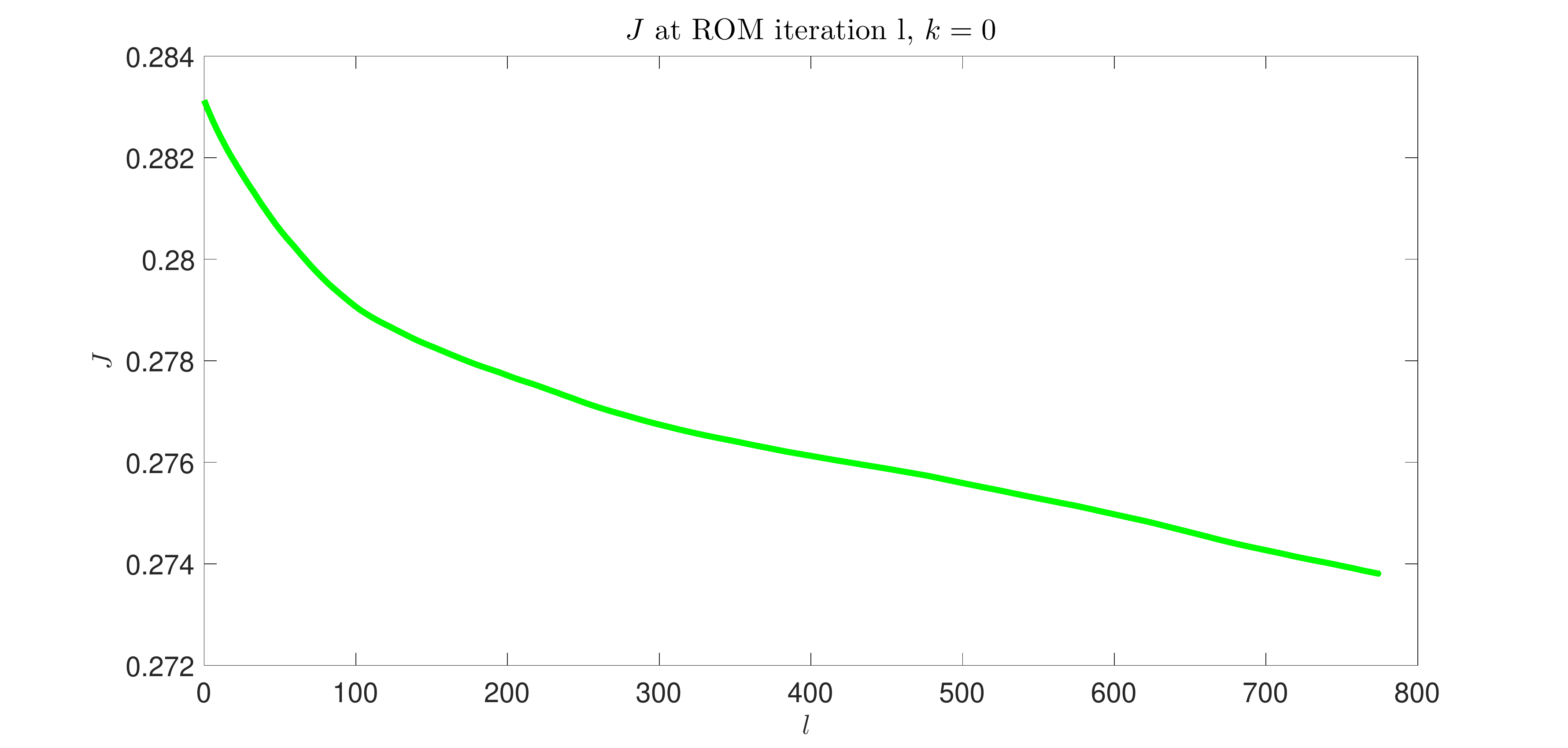}}&\raisebox{-0.5\height}{\includegraphics[scale=0.23]{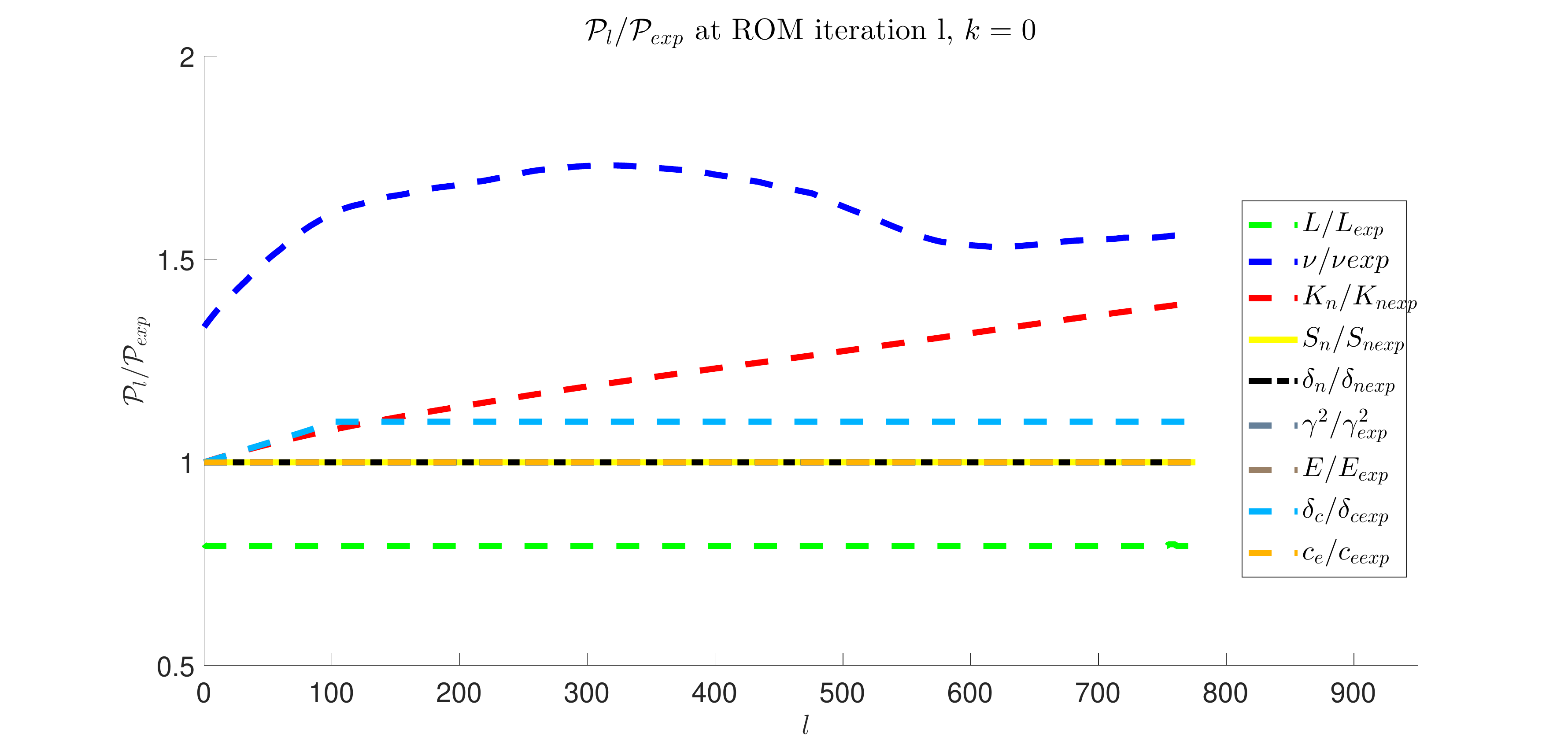}}&\\ \\
\parbox[t][][t]{5.6cm}{\centering $|\mathcal{P}_l(1)-\mathcal{P}_l|$} & & \\  
  \raisebox{-0.5\height}{\includegraphics[scale=0.23]{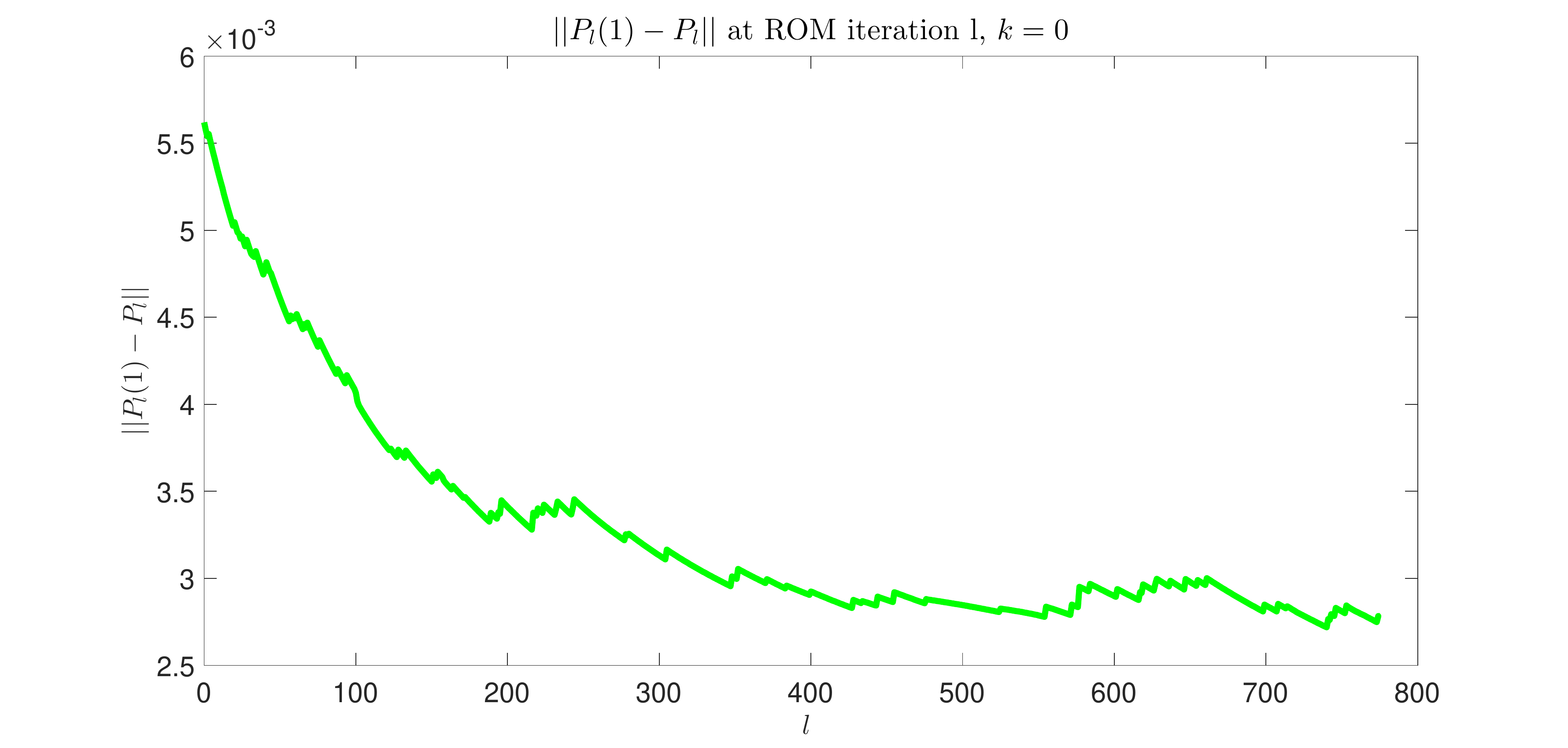}}& &\\ \\  
 \end{tabular}
 \begin{tabular}{c c c}
 \parbox[t][][t]{3.2cm}{\centering\textbf{MRI}} &   \parbox[t][][t]{7.2cm}{\centering\textbf{ROM ($\mathcal{P}_{\bar{l}=776})$}}  &  \parbox[t][][t]{4.2cm}{\centering\textbf{Comparison}}  \\ \\
\end{tabular}\\
\begin{tabular}{l}
\raisebox{-0.5\height}{\includegraphics[scale=0.8]{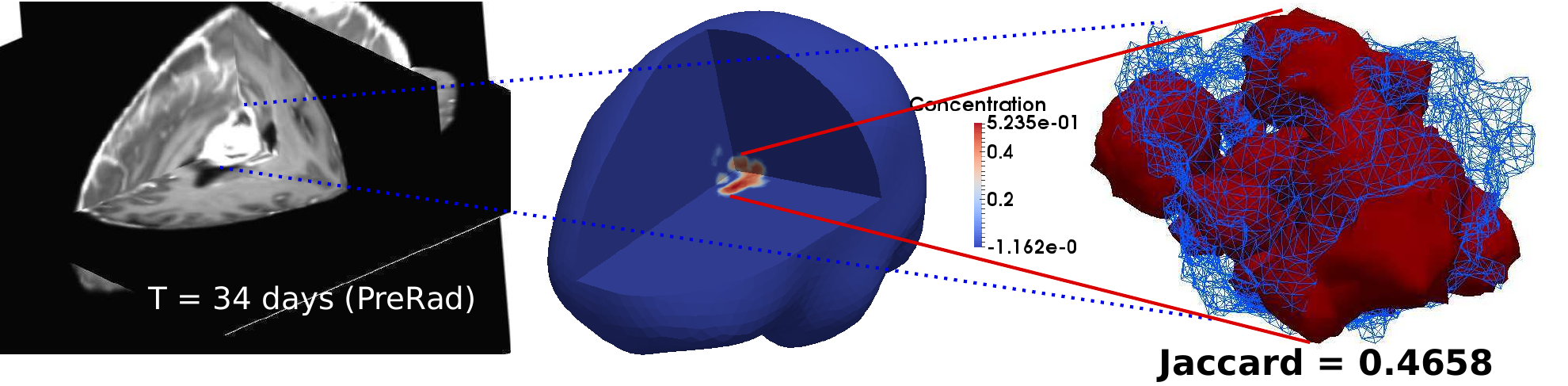}}  \\
  \end{tabular}
  \end{minipage} }
  \captionof{figure}[]{Values of the functional $J(\vec{\alpha}_l,\mathcal{P}_l)$, of the normalised set of parameters $\mathcal{P}_l/\mathcal{P}_{\text{exp}}$ and of $|\mathcal{P}_l(1)-\mathcal{P}_l|$ for steps $k=0$ of \textbf{Algorithm 1}, with $N=490$, $\mathcal{P}_0$ and $99.9\%$ POD threshold, together with a comparison between the isosurfaces $\phi_{\text{data}}(T)=0$ (highlighted in blue color) and $\sum_{i=1}^{N_{\text{POD}}}\alpha_{i\bar{l}}^N\xi_{i}^{\phi}=\phi_e/2$ (highlighted in red color).}
  \label{fig:7tristris}
\end{table}
We observe that in this case of a lower POD threshold of $99.9\%$ the value of the functional $J$ reaches a minimum value $0.2738$ which is higher than the minimal value $0.2667$ attained during the optimization algorithm with an higher threshold of $99.99\%$ (see Figure \ref{fig:7}). The local minimum is attained in a large number of steps $\bar{l}=776$.
Moreover, in the former case the parameter $\nu$ varies by a much smaller amount, the parameters $L$ and $\delta$ remain unchanged and the parameter $k_n$ varies through the functional minimisation. We also see more noise in the plot of $|\mathcal{P}_l(1)-\mathcal{P}_l|$ with respect to the case with $99.99\%$ threshold.

\section{Conclusions}
This work proposed a novel Optimization Algorithm \eqref{alg:opt} for the parameter estimation of a nonlinear  diffuse interface model  of GBM evolution from  neuroimaging data. The algorithm solves a constrained optimization problem in the form of a MPEC defined in \textbf{Problem 8} \eqref{eqn:optfom}, where the objective functional is defined in \eqref{eqn:2fom} by a proper distance between the tumour domain set in numerical simulations and the corresponding clinical data at a key time of clinical interest. Due to the high computational cost of solving the FOM (defined in \textbf{Problem 2} \eqref{eqn:3}), the algorithm iteratively computes the optimization functional at the FOM level and minimises it at the ROM level through sensitivity analysis. \\
First, an effective MOR is defined through \textbf{Problems 3, 4, 5} (see
\eqref{eqn:pod}, \eqref{eqn:assemble} and \eqref{eqn:romnewton}), by
projecting the equations onto the POD basis obtained from the time
snapshots of the FOM solutions and approximating the nonlinear terms
through DEIM interpolation. The nonlinearity of the ROM system is
solved by a Newton algorithm,  interpolating the
nonlinear terms associated to the singular potential and its first derivative on the same interpolation nodes to enforce there the separation property of the FOM solution. Moreover, the degenerate mobility and chemotactic terms are transformed as higher order tensors acting on the reduced order space,  taking into account for both  the degeneracy of the corresponding terms at the FOM level and  the  heterogeneity and anisotropy of the tumour invasion pattern.\\
A sensitivity analysis is performed at the ROM level to minimise the optimization functional by solving linearised systems defined in \textbf{Problem 6} \eqref{eqn:romlinearized} and updating the parameters along weighted gradient directions solving \textbf{Problem 9} \eqref{eqn:pgl}.\\
The algorithm ensures that the ROM solutions do not violate the physical constraints satisfied by the FOM solutions while dynamically spanning the space of parameters throughout the iterative steps. In particular, new POD basis functions associated to the updated set of parameters are calculated when the ROM minimisation problem has converged. This is an alternative way to consider parameters variability in the MOR of evolution equations with respect to the local reduced basis method used e.g. in \cite{manzoni}.

Second,  we applied the proposed  algorithm using the neuroimaging data provided by two clinical test cases: the growth of a primary GBM and a recurrent GBM  after surgical resection.\\
In both cases we observed the convergence of the algorithm to an optimal state, represented by the optimal sets of parameters \eqref{eqn:opttest1} and \eqref{eqn:opttest2}. These two sets display different optimal values, since the GBM spreading dynamics during primary and recurrent growth are controlled by intrinsically different biological processes. Moreover, the difference in the two set of parameters reflects also the higher aggressiveness of giant GBM cells in Test Case $2$, which is reflected in higher values of motility and proliferation.\\
In both cases we observed that the ROM solution approximates the FOM solution with a very high fidelity, and that the accuracy in reproducing the tumour domain from neuroimaging data increases at each step during the optimization process. \\
The number of basis functions needed to explain the $99.99\%$ variance of the data and the kind of morphological informations contained in them are the same in both cases, probably due to the  fact that there are no significant topological changes in the evolution dynamics, that occurs smoothly   thanks to chosen mesh refinement in the tumour evolution region. This turned out to be an important feature in order to deal with low dimensional higher--order tensors in \textbf{Problem 4} \eqref{eqn:assemble}:  ROM systems become indeed solvable with low computational resources and in highly reduced computational times.\\
We also observe that the computational time for the projected gradient iterations at the ROM level (\textbf{Step 4} of Algorithm \eqref{alg:opt}) is $2$ to $4$ orders of magnitude smaller than the time needed to solve the FOM problem, while the time for the assembly of the ROM systems (\textbf{Step 3}) is $2$ orders of magnitude smaller. This makes the whole algorithm very efficient in reducing the computational cost of the optimization process, both in the cases of the regular growth of a primary tumour and in the case of tumour recurrence with sparse particles and
infiltrations. The overall time of a sensitivity analysis performed at the FOM level is finally reduced by $2$ order of magnitude.\\
Finally we note that the Jaccard indexes calculated at the optimal
states for Test Case $2$ are smaller than the values computed for Test
Case $1$, and that  we need more iterations to converge to an optimal
state in the former case, with a smaller overall decrease of the
objective functional and a higher degree of oscillations in the
minimising sequences of parameters. This is due  partly to the sparse
distribution of tumour cells in the initial segmentation  after surgery,
and partly to the peritumoral infiltrations which may not be detectable by MRI data at the final time.\\

In conclusion, we proposed an optimization algorithm that allows a
robust parameter estimation of a nonlinear diffuse interface model
from neuroimaging data. The optimization is performed at low
computational cost using an automated procedure that allows to solve the
difficulties related the presence of a Cahn--Hilliard-type equation with single-well potential, non-conserved order parameter and degenerate mobility.  These features open the path to the creation of an automated computational platform that may be integrated in clinical practice to run simulations from neuroimaging data, thus to assist medical doctors in evaluating patient-specific therapeutic options.  Future developments will concern the  patient-specific therapy optimization in a given  temporal range and the assessment of uncertainty quantification of the underlying model. 

\section*{Acknowledgements}
This work was partly supported by MIUR, PRIN 2017 Research Project "Mathematics of active materials" and by the AIRC grant MFAG 17412 awarder to PC. The authors are grateful to Dr. F. Acerbi and Dr. A. Bizzi (Istituto Neurologico Besta, Milan) for providing the neuroimaging data and their guidance on key aspects of the clinical practice.
\newpage
\section{Appendix}
We report here the DEIM algorithm for the nonlinear term $\psi_1'(\phi)$ used to calculate the associated POD basis and row selection operator (see \cite{deim} for details):
\begin{algorithm}
\caption{DEIM Algorithm}\label{alg:deim}
\begin{algorithmic}
\Require Assemble the matrix $\Psi_{1,ij}:=\psi_1'^T(\phi_h^i)\psi_1'(\phi_h^j)$;\\
\textbf{Step A} Compute the POD basis $\xi_{1}^{\psi_1'}, \dots, \xi_{N_{\psi'}^{\text{POD}}}^{\psi_1'}$ for $\Psi_1$; \\
\textbf{Step B} $l\leftarrow \text{arg max}_{j=1,\dots,N_h}|\xi_{1}^{\psi_1'}(\mathbf{x}_j)|$; \\
\textbf{Step C} $U_{\psi_1'}=[\xi_{1}^{\psi_1'}]$, $i=l$, $P_1=[\vec{l}]$;
\For{$j=2, \dots, N_{\psi'}^{\text{POD}}$}  \\
$u \leftarrow \xi_{j}^{\psi_1'}$;\\
Solve $U_{\psi_1'i}c=u_i$; $r=u-U_{\psi_1'}c$;\\
$l\leftarrow \text{arg max}_{j=1,\dots,N_h}|r(x_j)|$;
$U_{\psi_1'}\leftarrow [U_{\psi_1'},u]$; $P_1=[P_1,\vec{l}]$,
\EndFor
\end{algorithmic}
\end{algorithm}
\\
where $\vec{l}$ is the finite element vector with value $1$ on the node $l$ and zero otherwise.
In the same way, we obtain $U_{\psi_1''}=(\xi_{1}^{\psi_1''}, \dots, \xi_{N_{\psi''}^{\text{POD}}}^{\psi_1''})$ and $P_2$ for $\psi_1''(\phi)$.\\
We also report here the Newton method $\mathbf{RN}_k$ \eqref{eqn:romnewton} used to solve \eqref{eqn:6}. Let us define the second order tensors
\[
B_{ml}(\vec{\alpha}_{k}^{n-1}):=\alpha_{ik}^{n-1}\alpha_{jk}^{n-1}\alpha_{sk}^{n-1} V_{2,ijsml} -2\alpha_{ik}^{n-1}\alpha_{jk}^{n-1}V_{3,ijml}+\alpha_{ik}^{n-1}V_{4,iml},
\]
\[
K_{ml}(\vec{\alpha}_{k}^{n-1}):=\alpha_{ik}^{n-1}\alpha_{jk}^{n-1}\alpha_{sk}^{n-1} V_{8,ijsml} -2\alpha_{ik}^{n-1}\alpha_{jk}^{n-1}V_{9,ijml}+\alpha_{ik}^{n-1}V_{10,iml},
\]
\[
V_{5,ml}(\vec{\alpha}_{k}^{n-1}):=\alpha_{ik}^{n-1} V_{5,iml} \quad V_{6,ml}(\vec{\alpha}_{k}^{n-1}):= \alpha_{ik}^{n-1} \alpha_{jk}^{n-1} V_{6,ijml},
\]
\[
V_{7,ml}(\vec{\alpha}_{k}^{n-1}):= \alpha_{ik}^{n-1} V_{7,iml} \quad U_{3,ml}(\vec{\alpha}_{k}^{n-1}):=\alpha_{ik}^{n-1}U_{3,iml} \quad W_{3,ml}(\vec{\alpha}_{k}^{n-1}):=\alpha_{ik}^{n-1}W_{3,iml}.
\]
We use the following algorithm to solve the  {\bf ROM Newton } problem.
\begin{algorithm}[!h]
\caption{ROM Newton Algorithm}\label{alg:rom}
\begin{algorithmic}
\State \textbf{Initialisation}\[\vec{\eta}_k^n=\biggl[\frac{1}{\Delta t}W_1+W_2+(\delta_{nk}-S_{nk})W_3(\vec{\alpha}_{k}^{n-1})+S_{nk}W_1\biggr]^{-1}\biggl(\frac{1}{\Delta t}W_1\vec{\eta}_{k}^{n-1}+S_{nk}W_4-S_{nk}W_5\vec{\alpha}_{k}^{n-1}\biggr),\]
\[\vec{\alpha}_{k}^{n,0}=\vec{\alpha}_{k}^{n-1},\vec{\beta}_{k}^{n,0}=\vec{\beta}_{k}^{n-1};\]
\[\text{error}=1, \; p=0.\]
\While{$\text{error}>10^{-3}$ and $p<1000$}
    \State \textbf{Assemble second order tensors:}
    \[U_{22,ml}(\vec{\alpha}_{k}^{n,p}):=(P_{2,ij}^TU_{\psi'',js})^{-1}\psi_1''(P_{2,sn}^T\Phi_{nh}\alpha_{hk}^{n,p})U_{22,iml},
    \]
    \[
    C(\vec{\alpha}_{k}^{n,p}):=\gamma_k^2U_1^{-1}U_6+E_kc_{ek}U_1^{-1}U_{22}(\vec{\alpha}_{k}^{n,p}),
    \]
    \begin{align*}
    D(\vec{\alpha}_{k}^{n-1},\vec{\alpha}_{k}^{n,p},\vec{\beta}_{k}^{n,p}):=&-U_1\vec{\beta}_{k}^{n,p}+\gamma_k^2U_6\vec{\alpha}_{k}^{n,p}+E_kc_{ek}U_{21}(P_{2}^TU_{\psi'})^{-1}\psi_1'(P_{2}^T\Phi\alpha_{k}^{n,p})-\\ & E_kU_3(\vec{\alpha}_{k}^{n-1})\vec{\alpha}_{k}^{n-1}-E_kc_{ek}U_4\vec{\alpha}_{k}^{n-1}-E_kc_{ek}U_5,
    \end{align*}
    \State \textbf{Solve for first order variations:}
 \begin{align*}
 d\vec{\alpha}_k=&\biggl[\frac{V_1}{\Delta t}+L_kB(\vec{\alpha}_{k}^{n-1})C(\vec{\alpha}_{k}^{n,p})\biggr]^{-1}\biggl(-L_kB(\vec{\alpha}_{k}^{n-1})U_1^{-1}D(\vec{\alpha}_{k}^{n-1},\vec{\alpha}_{k}^{n,p},\vec{\beta}_{k}^{n,p})-\frac{V_1}{\Delta t}\vec{\alpha}_{k}^{n,p}-\\
 &L_kB(\vec{\alpha}_{k}^{n-1})\vec{\beta}_{k}^{n,p}+k_{nk}K(\vec{\alpha}_{k}^{n-1})\vec{\eta}_k^n+\nu_k(V_5(\vec{\alpha}_{k}^{n-1})-V_6(\vec{\alpha}_{k}^{n-1}))\vec{\eta}_{k}^{n}+\\
 &\biggl(\nu_k\delta_k (V_7(\vec{\alpha}_{k}^{n-1})- V_1)+\frac{V_1}{\Delta t}-k_TV_1\biggr)\vec{\alpha}_{k}^{n-1}\biggr),
 \end{align*}
 \begin{equation*}
 d\vec{\beta}_k=C(\vec{\alpha}_{k}^{n,p})d\vec{\alpha}_k+U_1^{-1}D(\vec{\alpha}_{k}^{n-1},\vec{\alpha}_{k}^{n,p},\vec{\beta}_{k}^{n,p}),
 \end{equation*}
 \State \textbf{Update the Newton iterations:}
 \State \[\vec{\alpha}_{k}^{n,p+1}=\vec{\alpha}_{k}^{n,p}+d\vec{\alpha}_k, \quad \vec{\beta}_{k}^{n,p+1}=\vec{\beta}_{k}^{n,p}+d\vec{\beta}_k\quad \]
 \State \[\text{error}=\sqrt{d\vec{\alpha}_k^Td\vec{\alpha}_k+d\vec{\beta}_k^Td\vec{\beta}_k},\] 
\State \[p=p+1.\]    
\EndWhile
\State \textbf{Update the time step:}
\[
\vec{\alpha}_k^n=\vec{\alpha}_{k}^{n,p+1}, \quad \vec{\beta}_k^n=\vec{\beta}_{k}^{n,p+1}.
\]
\end{algorithmic}
\end{algorithm}\\
\newpage
We finally define the ROM linearised systems $\mathbf{RL}_k$ \eqref{eqn:romlinearized}, obtained by varying $\mathcal{P}_k=\mathcal{P}_k+\delta \mathcal{P}$ in \eqref{eqn:2}. Let us define the second order tensors
\[
B'_{ml,\mathcal{P}_{kh}}\biggl(\vec{\alpha}_{k}^{n-1},\frac{\partial \vec{\alpha}_{k}^{n-1}}{\partial \mathcal{P}_{kh}}\biggr):=3\frac{\partial \alpha_{ik}^{n-1}}{\partial \mathcal{P}_{kh}}\alpha_{jk}^{n-1}\alpha_{sk}^{n-1} V_{2,ijsml} -4\frac{\partial \alpha_{ik}^{n-1}}{\partial \mathcal{P}_{kh}}\alpha_{jk}^{n-1}V_{3,ijml}+\frac{\partial \alpha_{ik}^{n-1}}{\partial \mathcal{P}_{kh}}V_{4,iml},
\]
\[
K'_{ml,\mathcal{P}_{kh}}\biggl(\vec{\alpha}_{k}^{n-1},\frac{\partial \vec{\alpha}_{k}^{n-1}}{\partial \mathcal{P}_{kh}}\biggr):=3\frac{\partial \alpha_{ik}^{n-1}}{\partial \mathcal{P}_{kh}}\alpha_{jk}^{n-1}\alpha_{sk}^{n-1} V_{8,ijsml} -4\frac{\partial \alpha_{ik}^{n-1}}{\partial \mathcal{P}_{kh}}\alpha_{jk}^{n-1}V_{9,ijml}+\frac{\partial \alpha_{ik}^{n-1}}{\partial \mathcal{P}_{kh}}V_{10,iml},
\]
\[
V'_{5,ml}\biggl(\frac{\partial \vec{\alpha}_{k}^{n-1}}{\partial \mathcal{P}_{kh}}\biggr):=\frac{\partial \alpha_{ik}^{n-1}}{\partial \mathcal{P}_{kh}} V_{5,iml}, \quad V'_{6,ml}\biggl(\vec{\alpha}_{k}^{n-1},\frac{\partial \vec{\alpha}_{k}^{n-1}}{\partial \mathcal{P}_{kh}}\biggr):=2 \frac{\partial \alpha_{ik}^{n-1}}{\partial \mathcal{P}_{kh}} \alpha_{jk}^{n-1} V_{6,ijml},
\]
\[
V'_{7,ml}\biggl(\frac{\partial \vec{\alpha}_{k}^{n-1}}{\partial \mathcal{P}_{kh}}\biggr):= \frac{\partial \alpha_{ik}^{n-1}}{\partial \mathcal{P}_{kh}} V_{7,iml}, \quad U'_{3,ml}\biggl(\frac{\partial \vec{\alpha}_{k}^{n-1}}{\partial \mathcal{P}_{kh}}\biggr):=\frac{\partial \alpha_{ik}^{n-1}}{\partial \mathcal{P}_{kh}}U_{3,iml}, \quad W'_{3,ml}\biggl(\frac{\partial \vec{\alpha}_{k}^{n-1}}{\partial \mathcal{P}_{kh}}\biggr):=\frac{\partial \alpha_{ik}^{n-1}}{\partial \mathcal{P}_{kh}}W_{3,iml}.
\]
where $i,j,k,s,l,m=1,\dots,N_{\text{POD}}$ and $\mathcal{P}_{kh}$ is the $h-$th element of the set $\mathcal{P}_{k}$, with $h=1,\dots,|\mathcal{P}_k|$. Then we obtain the following \textbf{ROM} linearised systems, for $n=0,\dots,N,$\newline \\
\textbf{Linearised system for} $\mathbf{L}_k\to \mathbf{L}_k+\delta \mathbf{L}$:
\begin{equation}
\label{eqn:7}
\begin{cases}
\displaystyle \frac{V_{1}}{\Delta t}\frac{\partial \vec{\alpha}_{k}^{n}}{\partial L_k}&=-L_kB(\vec{\alpha}_{k}^{n-1})\frac{\partial \vec{\beta}_{k}^{n}}{\partial L_k}-L_kB'_{L_k}\biggl(\vec{\alpha}_{k}^{n-1},\frac{\partial \vec{\alpha}_{k}^{n-1}}{\partial L_k}\biggr)\vec{\beta}_{k}^{n}-B(\vec{\alpha}_{k}^{n-1})\vec{\beta}_{k}^{n}+k_{nk}K(\vec{\alpha}_{k}^{n-1})\frac{\partial \vec{\eta}_{k}^{n}}{\partial L_k}+\\
& \displaystyle k_{nk}K'\biggl(\vec{\alpha}_{k}^{n-1},\frac{\partial \vec{\alpha}_{k}^{n-1}}{\partial L_k}\biggr)\vec{\eta}_{k}^{n}+\nu_k \biggl(V'_{5}\biggl(\frac{\partial \vec{\alpha}_{k}^{n-1}}{\partial L_k}\biggr)-V'_{6}\biggl(\vec{\alpha}_{k}^{n-1},\frac{\partial \vec{\alpha}_{k}^{n-1}}{\partial L_k}\biggr)\biggr)\vec{\eta}_{k}^{n} +\\
\displaystyle & \nu_k \bigl(V_{5}(\vec{\alpha}_{k}^{n-1})-V_{6}(\vec{\alpha}_{k}^{n-1})\bigr)\frac{\partial \vec{\eta}_{k}^{n}}{\partial L_k} +\biggl(\nu_k\delta_k (V_7(\vec{\alpha}_{k}^{n-1})-V_1)+\biggl(\frac{1}{\Delta t} - K_T\biggr)V_{1}\biggr)\frac{\partial \vec{\alpha}_{k}^{n-1}}{\partial L_k}+\\
\displaystyle &\nu_k \delta_k V'_{7}\biggl(\frac{\partial \vec{\alpha}_{k}^{n-1}}{\partial L_k}\biggr)\vec{\alpha}_{k}^{n-1},\\ \\
\displaystyle U_{1}\frac{\partial \vec{\beta}_{k}^{n}}{\partial L_k}&=\gamma_k^2U_{6}\frac{ \partial \vec{\alpha}_{k}^{n}}{\partial L_k}+E_kc_{ek}U_{22}(\vec{\alpha}_{k}^{n})\frac{ \partial \vec{\alpha}_{k}^{n}}{\partial L_k}-E_kU'_{3}\biggl(\frac{ \partial \vec{\alpha}_{k}^{n-1}}{\partial L_k}\biggr)\vec{\alpha}_{k}^{n-1}-E_kU_3(\vec{\alpha}_{k}^{n-1})\frac{ \partial \vec{\alpha}_{k}^{n-1}}{\partial L_k}-\\
\displaystyle &E_kc_{ek}U_{4}\frac{ \partial \vec{\alpha}_{k}^{n-1}}{\partial L_k},\\ \\
\displaystyle \frac{W_{1}}{\Delta t}\frac{\partial \vec{\eta}_{k}^{n}}{\partial L_k}&=-W_{2}\frac{\partial \vec{\eta}_{k}^{n}}{\partial L_k}+(S_{nk}-\delta_{nk})W'_3\biggl(\frac{\partial \vec{\alpha}_{k}^{n-1}}{\partial L_k}\biggr)\vec{\eta}_{k}^{n}+\biggl((S_{nk}-\delta_{nk}) W_{3}(\vec{\alpha}_{k}^{n-1})-S_{nk}W_{1}\biggr)\frac{\partial \vec{\eta}_{k}^{n}}{\partial L_k}\\
\displaystyle &-S_{nk}W_5\frac{\partial \vec{\alpha}_{k}^{n-1}}{\partial L_k}+\frac{W_{1}}{\Delta t}\frac{\partial \vec{\eta}_{k}^{n-1}}{\partial L_k}.
\end{cases}
\end{equation}\\
\textbf{Linearised system for} $\boldsymbol{\nu}_k\to \boldsymbol{\nu}_k+\delta \boldsymbol{\nu}$:
\begin{equation}
\label{eqn:8}
\begin{cases}
\displaystyle \frac{V_{1}}{\Delta t}\frac{\partial \vec{\alpha}_{k}^{n}}{\partial \nu_k}&=-L_kB(\vec{\alpha}_{k}^{n-1})\frac{\partial \vec{\beta}_{k}^{n}}{\partial \nu_k}-L_kB'_{L_k}\biggl(\vec{\alpha}_{k}^{n-1},\frac{\partial \vec{\alpha}_{k}^{n-1}}{\partial \nu_k}\biggr)\vec{\beta}_{k}^{n}+k_{nk}K(\vec{\alpha}_{k}^{n-1})\frac{\partial \vec{\eta}_{k}^{n}}{\partial \nu_k}+\\
& \displaystyle k_{nk}K'\biggl(\vec{\alpha}_{k}^{n-1},\frac{\partial \vec{\alpha}_{k}^{n-1}}{\partial \nu_k}\biggr)\vec{\eta}_{k}^{n}+\nu_k \biggl(V'_{5}\biggl(\frac{\partial \vec{\alpha}_{k}^{n-1}}{\partial \nu_k}\biggr)-V'_{6}\biggl(\vec{\alpha}_{k}^{n-1},\frac{\partial \vec{\alpha}_{k}^{n-1}}{\partial \nu_k}\biggr)\biggr)\vec{\eta}_{k}^{n} +\\
\displaystyle & \nu_k \bigl(V_{5}(\vec{\alpha}_{k}^{n-1})-V_{6}(\vec{\alpha}_{k}^{n-1})\bigr)\frac{\partial \vec{\eta}_{k}^{n}}{\partial \nu_k} +\biggl(\nu_k\delta_k (V_7(\vec{\alpha}_{k}^{n-1})-V_1)+\biggl(\frac{1}{\Delta t} - K_T\biggr)V_{1}\biggr)\frac{\partial \vec{\alpha}_{k}^{n-1}}{\partial \nu_k}+\\
\displaystyle &\nu_k \delta_k V'_{7}\biggl(\frac{\partial \vec{\alpha}_{k}^{n-1}}{\partial \nu_k}\biggr)\vec{\alpha}_{k}^{n-1}+\bigl(V_{5}(\vec{\alpha}_{k}^{n-1})-V_{6}(\vec{\alpha}_{k}^{n-1})\bigr)\vec{\eta}_{k}^{n}+\delta_k (V_7(\vec{\alpha}_{k}^{n-1})-V_1)\vec{\alpha}_{k}^{n-1},\\ \\
\displaystyle U_{1}\frac{\partial \vec{\beta}_{k}^{n}}{\partial \nu_k}&=\gamma_k^2U_{6}\frac{ \partial \vec{\alpha}_{k}^{n}}{\partial \nu_k}+E_kc_{ek}U_{22}(\vec{\alpha}_{k}^{n})\frac{ \partial \vec{\alpha}_{k}^{n}}{\partial \nu_k}-E_kU'_{3}\biggl(\frac{ \partial \vec{\alpha}_{k}^{n-1}}{\partial \nu_k}\biggr)\vec{\alpha}_{k}^{n-1}-E_kU_3(\vec{\alpha}_{k}^{n-1})\frac{ \partial \vec{\alpha}_{k}^{n-1}}{\partial \nu_k}-\\
\displaystyle &E_kc_{ek}U_{4}\frac{ \partial \vec{\alpha}_{k}^{n-1}}{\partial \nu_k},\\ \\
\displaystyle \frac{W_{1}}{\Delta t}\frac{\partial \vec{\eta}_{k}^{n}}{\partial \nu_k}&=-W_{2}\frac{\partial \vec{\eta}_{k}^{n}}{\partial \nu_k}+(S_{nk}-\delta_{nk})W'_3\biggl(\frac{\partial \vec{\alpha}_{k}^{n-1}}{\partial \nu_k}\biggr)\vec{\eta}_{k}^{n}+\biggl((S_{nk}-\delta_{nk}) W_{3}(\vec{\alpha}_{k}^{n-1})-S_{nk}W_{1}\biggr)\frac{\partial \vec{\eta}_{k}^{n}}{\partial \nu_k}\\
\displaystyle &-S_{nk}W_5\frac{\partial \vec{\alpha}_{k}^{n-1}}{\partial \nu_k}+\frac{W_{1}}{\Delta t}\frac{\partial \vec{\eta}_{k}^{n-1}}{\partial \nu_k}.
\end{cases}
\end{equation}\\
\textbf{Linearised system for} $\mathbf{k}_{nk}\to \mathbf{k}_{nk}+\delta \mathbf{k}_{n}$:
\begin{equation}
\label{eqn:9}
\begin{cases}
\displaystyle \frac{V_{1}}{\Delta t}\frac{\partial \vec{\alpha}_{k}^{n}}{\partial k_{nk}}&=-L_kB(\vec{\alpha}_{k}^{n-1})\frac{\partial \vec{\beta}_{k}^{n}}{\partial k_{nk}}-L_kB'_{L_k}\biggl(\vec{\alpha}_{k}^{n-1},\frac{\partial \vec{\alpha}_{k}^{n-1}}{\partial k_{nk}}\biggr)\vec{\beta}_{k}^{n}+k_{nk}K(\vec{\alpha}_{k}^{n-1})\frac{\partial \vec{\eta}_{k}^{n}}{\partial k_{nk}}+\\
& \displaystyle k_{nk}K'\biggl(\vec{\alpha}_{k}^{n-1},\frac{\partial \vec{\alpha}_{k}^{n-1}}{\partial k_{nk}}\biggr)\vec{\eta}_{k}^{n}+K(\vec{\alpha}_{k}^{n-1})\vec{\eta}_{k}^{n}+\nu_k \biggl(V'_{5}\biggl(\frac{\partial \vec{\alpha}_{k}^{n-1}}{\partial k_{nk}}\biggr)-V'_{6}\biggl(\vec{\alpha}_{k}^{n-1},\frac{\partial \vec{\alpha}_{k}^{n-1}}{\partial k_{nk}}\biggr)\biggr)\vec{\eta}_{k}^{n} +\\
\displaystyle & \nu_k \bigl(V_{5}(\vec{\alpha}_{k}^{n-1})-V_{6}(\vec{\alpha}_{k}^{n-1})\bigr)\frac{\partial \vec{\eta}_{k}^{n}}{\partial k_{nk}} +\biggl(\nu_k\delta_k (V_7(\vec{\alpha}_{k}^{n-1})-V_1)+\biggl(\frac{1}{\Delta t} - K_T\biggr)V_{1}\biggr)\frac{\partial \vec{\alpha}_{k}^{n-1}}{\partial k_{nk}}+\\
\displaystyle &\nu_k \delta_k V'_{7}\biggl(\frac{\partial \vec{\alpha}_{k}^{n-1}}{\partial k_{nk}}\biggr)\vec{\alpha}_{k}^{n-1},\\ \\
\displaystyle U_{1}\frac{\partial \vec{\beta}_{k}^{n}}{\partial k_{nk}}&=\gamma_k^2U_{6}\frac{ \partial \vec{\alpha}_{k}^{n}}{\partial k_{nk}}+E_kc_{ek}U_{22}(\vec{\alpha}_{k}^{n})\frac{ \partial \vec{\alpha}_{k}^{n}}{\partial k_{nk}}-E_kU'_{3}\biggl(\frac{ \partial \vec{\alpha}_{k}^{n-1}}{\partial k_{nk}}\biggr)\vec{\alpha}_{k}^{n-1}-E_kU_3(\vec{\alpha}_{k}^{n-1})\frac{ \partial \vec{\alpha}_{k}^{n-1}}{\partial k_{nk}}-\\
\displaystyle &E_kc_{ek}U_{4}\frac{ \partial \vec{\alpha}_{k}^{n-1}}{\partial k_{nk}},\\ \\
\displaystyle \frac{W_{1}}{\Delta t}\frac{\partial \vec{\eta}_{k}^{n}}{\partial k_{nk}}&=-W_{2}\frac{\partial \vec{\eta}_{k}^{n}}{\partial k_{nk}}+(S_{nk}-\delta_{nk})W'_3\biggl(\frac{\partial \vec{\alpha}_{k}^{n-1}}{\partial k_{nk}}\biggr)\vec{\eta}_{k}^{n}+\biggl((S_{nk}-\delta_{nk}) W_{3}(\vec{\alpha}_{k}^{n-1})-S_{nk}W_{1}\biggr)\frac{\partial \vec{\eta}_{k}^{n}}{\partial k_{nk}}\\
\displaystyle &-S_{nk}W_5\frac{\partial \vec{\alpha}_{k}^{n-1}}{\partial k_{nk}}+\frac{W_{1}}{\Delta t}\frac{\partial \vec{\eta}_{k}^{n-1}}{\partial k_{nk}}.
\end{cases}
\end{equation}
\\
\textbf{Linearised system for} $\mathbf{S}_{nk}\to \mathbf{S}_{nk}+\delta \mathbf{S}_{n}$:
\begin{equation}
\label{eqn:10}
\begin{cases}
\displaystyle \frac{V_{1}}{\Delta t}\frac{\partial \vec{\alpha}_{k}^{n}}{\partial S_{nk}}&=-L_kB(\vec{\alpha}_{k}^{n-1})\frac{\partial \vec{\beta}_{k}^{n}}{\partial S_{nk}}-L_kB'_{L_k}\biggl(\vec{\alpha}_{k}^{n-1},\frac{\partial \vec{\alpha}_{k}^{n-1}}{\partial S_{nk}}\biggr)\vec{\beta}_{k}^{n}+k_{nk}K(\vec{\alpha}_{k}^{n-1})\frac{\partial \vec{\eta}_{k}^{n}}{\partial S_{nk}}+\\
& \displaystyle k_{nk}K'\biggl(\vec{\alpha}_{k}^{n-1},\frac{\partial \vec{\alpha}_{k}^{n-1}}{\partial S_{nk}}\biggr)\vec{\eta}_{k}^{n}+\nu_k \biggl(V'_{5}\biggl(\frac{\partial \vec{\alpha}_{k}^{n-1}}{\partial S_{nk}}\biggr)-V'_{6}\biggl(\vec{\alpha}_{k}^{n-1},\frac{\partial \vec{\alpha}_{k}^{n-1}}{\partial S_{nk}}\biggr)\biggr)\vec{\eta}_{k}^{n} +\\
\displaystyle & \nu_k \bigl(V_{5}(\vec{\alpha}_{k}^{n-1})-V_{6}(\vec{\alpha}_{k}^{n-1})\bigr)\frac{\partial \vec{\eta}_{k}^{n}}{\partial S_{nk}} +\biggl(\nu_k\delta_k (V_7(\vec{\alpha}_{k}^{n-1})-V_1)+\biggl(\frac{1}{\Delta t} - K_T\biggr)V_{1}\biggr)\frac{\partial \vec{\alpha}_{k}^{n-1}}{\partial S_{nk}}+\\
\displaystyle &\nu_k \delta_k V'_{7}\biggl(\frac{\partial \vec{\alpha}_{k}^{n-1}}{\partial S_{nk}}\biggr)\vec{\alpha}_{k}^{n-1},\\ \\
\displaystyle U_{1}\frac{\partial \vec{\beta}_{k}^{n}}{\partial S_{nk}}&=\gamma_k^2U_{6}\frac{ \partial \vec{\alpha}_{k}^{n}}{\partial S_{nk}}+E_kc_{ek}U_{22}(\vec{\alpha}_{k}^{n})\frac{ \partial \vec{\alpha}_{k}^{n}}{\partial S_{nk}}-E_kU'_{3}\biggl(\frac{ \partial \vec{\alpha}_{k}^{n-1}}{\partial S_{nk}}\biggr)\vec{\alpha}_{k}^{n-1}-E_kU_3(\vec{\alpha}_{k}^{n-1})\frac{ \partial \vec{\alpha}_{k}^{n-1}}{\partial S_{nk}}-\\
\displaystyle &E_kc_{ek}U_{4}\frac{ \partial \vec{\alpha}_{k}^{n-1}}{\partial S_{nk}},\\ \\
\displaystyle \frac{W_{1}}{\Delta t}\frac{\partial \vec{\eta}_{k}^{n}}{\partial S_{nk}}&=-W_{2}\frac{\partial \vec{\eta}_{k}^{n}}{\partial S_{nk}}+(S_{nk}-\delta_{nk})W'_3\biggl(\frac{\partial \vec{\alpha}_{k}^{n-1}}{\partial S_{nk}}\biggr)\vec{\eta}_{k}^{n}+\biggl((S_{nk}-\delta_{nk}) W_{3}(\vec{\alpha}_{k}^{n-1})-S_{nk}W_{1}\biggr)\frac{\partial \vec{\eta}_{k}^{n}}{\partial S_{nk}}\\
\displaystyle &-S_{nk}W_5\frac{\partial \vec{\alpha}_{k}^{n-1}}{\partial S_{nk}}+\frac{W_{1}}{\Delta t}\frac{\partial \vec{\eta}_{k}^{n-1}}{\partial S_{nk}}+W_4-W_5\vec{\alpha}_{k}^{n-1}+(W_{3}(\vec{\alpha}_{k}^{n-1})-W_1)\vec{\eta}_{k}^{n}.
\end{cases}
\end{equation}
\\
\textbf{Linearised system for} $\boldsymbol{\delta}_{nk}\to \boldsymbol{\delta}_{nk}+\delta \boldsymbol{\delta}_{n}$:
\begin{equation}
\label{eqn:11}
\begin{cases}
\displaystyle \frac{V_{1}}{\Delta t}\frac{\partial \vec{\alpha}_{k}^{n}}{\partial \delta_{nk}}&=-L_kB(\vec{\alpha}_{k}^{n-1})\frac{\partial \vec{\beta}_{k}^{n}}{\partial \delta_{nk}}-L_kB'_{L_k}\biggl(\vec{\alpha}_{k}^{n-1},\frac{\partial \vec{\alpha}_{k}^{n-1}}{\partial \delta_{nk}}\biggr)\vec{\beta}_{k}^{n}+k_{nk}K(\vec{\alpha}_{k}^{n-1})\frac{\partial \vec{\eta}_{k}^{n}}{\partial \delta_{nk}}+\\
& \displaystyle k_{nk}K'\biggl(\vec{\alpha}_{k}^{n-1},\frac{\partial \vec{\alpha}_{k}^{n-1}}{\partial \delta_{nk}}\biggr)\vec{\eta}_{k}^{n}+\nu_k \biggl(V'_{5}\biggl(\frac{\partial \vec{\alpha}_{k}^{n-1}}{\partial \delta_{nk}}\biggr)-V'_{6}\biggl(\vec{\alpha}_{k}^{n-1},\frac{\partial \vec{\alpha}_{k}^{n-1}}{\partial \delta_{nk}}\biggr)\biggr)\vec{\eta}_{k}^{n} +\\
\displaystyle & \nu_k \bigl(V_{5}(\vec{\alpha}_{k}^{n-1})-V_{6}(\vec{\alpha}_{k}^{n-1})\bigr)\frac{\partial \vec{\eta}_{k}^{n}}{\partial \delta_{nk}} +\biggl(\nu_k\delta_k (V_7(\vec{\alpha}_{k}^{n-1})-V_1)+\biggl(\frac{1}{\Delta t} - K_T\biggr)V_{1}\biggr)\frac{\partial \vec{\alpha}_{k}^{n-1}}{\partial \delta_{nk}}+\\
\displaystyle &\nu_k \delta_k V'_{7}\biggl(\frac{\partial \vec{\alpha}_{k}^{n-1}}{\partial \delta_{nk}}\biggr)\vec{\alpha}_{k}^{n-1},\\ \\
\displaystyle U_{1}\frac{\partial \vec{\beta}_{k}^{n}}{\partial \delta_{nk}}&=\gamma_k^2U_{6}\frac{ \partial \vec{\alpha}_{k}^{n}}{\partial \delta_{nk}}+E_kc_{ek}U_{22}(\vec{\alpha}_{k}^{n})\frac{ \partial \vec{\alpha}_{k}^{n}}{\partial \delta_{nk}}-E_kU'_{3}\biggl(\frac{ \partial \vec{\alpha}_{k}^{n-1}}{\partial \delta_{nk}}\biggr)\vec{\alpha}_{k}^{n-1}-E_kU_3(\vec{\alpha}_{k}^{n-1})\frac{ \partial \vec{\alpha}_{k}^{n-1}}{\partial \delta_{nk}}-\\
\displaystyle &E_kc_{ek}U_{4}\frac{ \partial \vec{\alpha}_{k}^{n-1}}{\partial \delta_{nk}},\\ \\
\displaystyle \frac{W_{1}}{\Delta t}\frac{\partial \vec{\eta}_{k}^{n}}{\partial \delta_{nk}}&=-W_{2}\frac{\partial \vec{\eta}_{k}^{n}}{\partial \delta_{nk}}+(S_{nk}-\delta_{nk})W'_3\biggl(\frac{\partial \vec{\alpha}_{k}^{n-1}}{\partial \delta_{nk}}\biggr)\vec{\eta}_{k}^{n}+\biggl((S_{nk}-\delta_{nk}) W_{3}(\vec{\alpha}_{k}^{n-1})-S_{nk}W_{1}\biggr)\frac{\partial \vec{\eta}_{k}^{n}}{\partial \delta_{nk}}\\
\displaystyle &-S_{nk}W_5\frac{\partial \vec{\alpha}_{k}^{n-1}}{\partial \delta_{nk}}+\frac{W_{1}}{\Delta t}\frac{\partial \vec{\eta}_{k}^{n-1}}{\partial \delta_{nk}}-W_{3}(\vec{\alpha}_{k}^{n-1})\vec{\eta}_{k}^{n}.
\end{cases}
\end{equation}
\\
\textbf{Linearised system for} $\boldsymbol{\gamma}_{k}\to \boldsymbol{\gamma}_{k}+\delta \boldsymbol{\gamma}$:
\begin{equation}
\label{eqn:12}
\begin{cases}
\displaystyle \frac{V_{1}}{\Delta t}\frac{\partial \vec{\alpha}_{k}^{n}}{\partial \gamma_{k}^2}&=-L_kB(\vec{\alpha}_{k}^{n-1})\frac{\partial \vec{\beta}_{k}^{n}}{\partial \gamma_{k}^2}-L_kB'_{L_k}\biggl(\vec{\alpha}_{k}^{n-1},\frac{\partial \vec{\alpha}_{k}^{n-1}}{\partial \gamma_{k}^2}\biggr)\vec{\beta}_{k}^{n}+k_{nk}K(\vec{\alpha}_{k}^{n-1})\frac{\partial \vec{\eta}_{k}^{n}}{\partial \gamma_{k}^2}+\\
& \displaystyle k_{nk}K'\biggl(\vec{\alpha}_{k}^{n-1},\frac{\partial \vec{\alpha}_{k}^{n-1}}{\partial \gamma_{k}^2}\biggr)\vec{\eta}_{k}^{n}+\nu_k \biggl(V'_{5}\biggl(\frac{\partial \vec{\alpha}_{k}^{n-1}}{\partial \gamma_{k}^2}\biggr)-V'_{6}\biggl(\vec{\alpha}_{k}^{n-1},\frac{\partial \vec{\alpha}_{k}^{n-1}}{\partial \gamma_{k}^2}\biggr)\biggr)\vec{\eta}_{k}^{n} +\\
\displaystyle & \nu_k \bigl(V_{5}(\vec{\alpha}_{k}^{n-1})-V_{6}(\vec{\alpha}_{k}^{n-1})\bigr)\frac{\partial \vec{\eta}_{k}^{n}}{\partial \gamma_{k}^2} +\biggl(\nu_k\delta_k (V_7(\vec{\alpha}_{k}^{n-1})-V_1)+\biggl(\frac{1}{\Delta t} - K_T\biggr)V_{1}\biggr)\frac{\partial \vec{\alpha}_{k}^{n-1}}{\partial \gamma_{k}^2}+\\
\displaystyle &\nu_k \delta_k V'_{7}\biggl(\frac{\partial \vec{\alpha}_{k}^{n-1}}{\partial \gamma_{k}^2}\biggr)\vec{\alpha}_{k}^{n-1},\\ \\
\displaystyle U_{1}\frac{\partial \vec{\beta}_{k}^{n}}{\partial \gamma_{k}^2}&=\gamma_k^2U_{6}\frac{ \partial \vec{\alpha}_{k}^{n}}{\partial \gamma_{k}^2}+E_kc_{ek}U_{22}(\vec{\alpha}_{k}^{n})\frac{ \partial \vec{\alpha}_{k}^{n}}{\partial \gamma_{k}^2}-E_kU'_{3}\biggl(\frac{ \partial \vec{\alpha}_{k}^{n-1}}{\partial \gamma_{k}^2}\biggr)\vec{\alpha}_{k}^{n-1}-E_kU_3(\vec{\alpha}_{k}^{n-1})\frac{ \partial \vec{\alpha}_{k}^{n-1}}{\partial \gamma_{k}^2}-\\
\displaystyle &E_kc_{ek}U_{4}\frac{ \partial \vec{\alpha}_{k}^{n-1}}{\partial \gamma_{k}^2}+U_6\vec{\alpha}_{k}^{n},\\ \\
\displaystyle \frac{W_{1}}{\Delta t}\frac{\partial \vec{\eta}_{k}^{n}}{\partial \gamma_{k}^2}&=-W_{2}\frac{\partial \vec{\eta}_{k}^{n}}{\partial \gamma_{k}^2}+(S_{nk}-\delta_{nk})W'_3\biggl(\frac{\partial \vec{\alpha}_{k}^{n-1}}{\partial \gamma_{k}^2}\biggr)\vec{\eta}_{k}^{n}+\biggl((S_{nk}-\delta_{nk}) W_{3}(\vec{\alpha}_{k}^{n-1})-S_{nk}W_{1}\biggr)\frac{\partial \vec{\eta}_{k}^{n}}{\partial \gamma_{k}^2}\\
\displaystyle &-S_{nk}W_5\frac{\partial \vec{\alpha}_{k}^{n-1}}{\partial \gamma_{k}^2}+\frac{W_{1}}{\Delta t}\frac{\partial \vec{\eta}_{k}^{n-1}}{\partial \gamma_{k}^2}.
\end{cases}
\end{equation}
\\
\textbf{Linearised system for} $\mathbf{E}_{k}\to \mathbf{E}_{k}+\delta \mathbf{E}$:
\begin{equation}
\label{eqn:13}
\begin{cases}
\displaystyle \frac{V_{1}}{\Delta t}\frac{\partial \vec{\alpha}_{k}^{n}}{\partial E_k}&=-L_kB(\vec{\alpha}_{k}^{n-1})\frac{\partial \vec{\beta}_{k}^{n}}{\partial E_k}-L_kB'_{L_k}\biggl(\vec{\alpha}_{k}^{n-1},\frac{\partial \vec{\alpha}_{k}^{n-1}}{\partial E_k}\biggr)\vec{\beta}_{k}^{n}+k_{nk}K(\vec{\alpha}_{k}^{n-1})\frac{\partial \vec{\eta}_{k}^{n}}{\partial E_k}+\\
& \displaystyle k_{nk}K'\biggl(\vec{\alpha}_{k}^{n-1},\frac{\partial \vec{\alpha}_{k}^{n-1}}{\partial E_k}\biggr)\vec{\eta}_{k}^{n}+\nu_k \biggl(V'_{5}\biggl(\frac{\partial \vec{\alpha}_{k}^{n-1}}{\partial E_k}\biggr)-V'_{6}\biggl(\vec{\alpha}_{k}^{n-1},\frac{\partial \vec{\alpha}_{k}^{n-1}}{\partial E_k}\biggr)\biggr)\vec{\eta}_{k}^{n} +\\
\displaystyle & \nu_k \bigl(V_{5}(\vec{\alpha}_{k}^{n-1})-V_{6}(\vec{\alpha}_{k}^{n-1})\bigr)\frac{\partial \vec{\eta}_{k}^{n}}{\partial E_k} +\biggl(\nu_k\delta_k (V_7(\vec{\alpha}_{k}^{n-1})-V_1)+\biggl(\frac{1}{\Delta t} - K_T\biggr)V_{1}\biggr)\frac{\partial \vec{\alpha}_{k}^{n-1}}{\partial E_k}+\\
\displaystyle &\nu_k \delta_k V'_{7}\biggl(\frac{\partial \vec{\alpha}_{k}^{n-1}}{\partial E_k}\biggr)\vec{\alpha}_{k}^{n-1},\\ \\
\displaystyle U_{1}\frac{\partial \vec{\beta}_{k}^{n}}{\partial E_k}&=\gamma_k^2U_{6}\frac{ \partial \vec{\alpha}_{k}^{n}}{\partial E_k}+E_kc_{ek}U_{22}(\vec{\alpha}_{k}^{n})\frac{ \partial \vec{\alpha}_{k}^{n}}{\partial E_k}-E_kU'_{3}\biggl(\frac{ \partial \vec{\alpha}_{k}^{n-1}}{\partial E_k}\biggr)\vec{\alpha}_{k}^{n-1}-E_kU_3(\vec{\alpha}_{k}^{n-1})\frac{ \partial \vec{\alpha}_{k}^{n-1}}{\partial E_k}-\\
\displaystyle &E_kc_{ek}U_{4}\frac{ \partial \vec{\alpha}_{k}^{n-1}}{\partial E_k}+c_{ek}U_{21}(P_{2}^TU_{\psi'})^{-1}\psi_1'(P_{2}^T\Phi\alpha_{k}^{n})-(U_3(\vec{\alpha}_{k}^{n-1})+c_{ek}U_{4})\vec{\alpha}_{k}^{n-1}-\\
\displaystyle & c_{ek}U_{5},\\ \\
\displaystyle \frac{W_{1}}{\Delta t}\frac{\partial \vec{\eta}_{k}^{n}}{\partial E_k}&=-W_{2}\frac{\partial \vec{\eta}_{k}^{n}}{\partial E_k}+(S_{nk}-\delta_{nk})W'_3\biggl(\frac{\partial \vec{\alpha}_{k}^{n-1}}{\partial E_k}\biggr)\vec{\eta}_{k}^{n}+\biggl((S_{nk}-\delta_{nk}) W_{3}(\vec{\alpha}_{k}^{n-1})-S_{nk}W_{1}\biggr)\frac{\partial \vec{\eta}_{k}^{n}}{\partial E_k}\\
\displaystyle &-S_{nk}W_5\frac{\partial \vec{\alpha}_{k}^{n-1}}{\partial E_k}+\frac{W_{1}}{\Delta t}\frac{\partial \vec{\eta}_{k}^{n-1}}{\partial E_k}.
\end{cases}
\end{equation}
\\
\textbf{Linearised system for} $\boldsymbol{\delta}_{k}\to \boldsymbol{\delta}_{k}+\delta \boldsymbol{\delta}$:
\begin{equation}
\label{eqn:14}
\begin{cases}
\displaystyle \frac{V_{1}}{\Delta t}\frac{\partial \vec{\alpha}_{k}^{n}}{\partial \delta_k}&=-L_kB(\vec{\alpha}_{k}^{n-1})\frac{\partial \vec{\beta}_{k}^{n}}{\partial \delta_k}-L_kB'_{L_k}\biggl(\vec{\alpha}_{k}^{n-1},\frac{\partial \vec{\alpha}_{k}^{n-1}}{\partial \delta_k}\biggr)\vec{\beta}_{k}^{n}+k_{nk}K(\vec{\alpha}_{k}^{n-1})\frac{\partial \vec{\eta}_{k}^{n}}{\partial \delta_k}+\\
& \displaystyle k_{nk}K'\biggl(\vec{\alpha}_{k}^{n-1},\frac{\partial \vec{\alpha}_{k}^{n-1}}{\partial \delta_k}\biggr)\vec{\eta}_{k}^{n}+\nu_k \biggl(V'_{5}\biggl(\frac{\partial \vec{\alpha}_{k}^{n-1}}{\partial \delta_k}\biggr)-V'_{6}\biggl(\vec{\alpha}_{k}^{n-1},\frac{\partial \vec{\alpha}_{k}^{n-1}}{\partial \delta_k}\biggr)\biggr)\vec{\eta}_{k}^{n} +\\
\displaystyle & \nu_k \bigl(V_{5}(\vec{\alpha}_{k}^{n-1})-V_{6}(\vec{\alpha}_{k}^{n-1})\bigr)\frac{\partial \vec{\eta}_{k}^{n}}{\partial \delta_k} +\biggl(\nu_k\delta_k (V_7(\vec{\alpha}_{k}^{n-1})-V_1)+\biggl(\frac{1}{\Delta t} - K_T\biggr)V_{1}\biggr)\frac{\partial \vec{\alpha}_{k}^{n-1}}{\partial \delta_k}+\\
\displaystyle &\nu_k \delta_k V'_{7}\biggl(\frac{\partial \vec{\alpha}_{k}^{n-1}}{\partial \delta_k}\biggr)\vec{\alpha}_{k}^{n-1}+\nu_k (V_7(\vec{\alpha}_{k}^{n-1})-V_1)\vec{\alpha}_{k}^{n-1},\\ \\
\displaystyle U_{1}\frac{\partial \vec{\beta}_{k}^{n}}{\partial \delta_k}&=\gamma_k^2U_{6}\frac{ \partial \vec{\alpha}_{k}^{n}}{\partial \delta_k}+E_kc_{ek}U_{22}(\vec{\alpha}_{k}^{n})\frac{ \partial \vec{\alpha}_{k}^{n}}{\partial \delta_k}-E_kU'_{3}\biggl(\frac{ \partial \vec{\alpha}_{k}^{n-1}}{\partial \delta_k}\biggr)\vec{\alpha}_{k}^{n-1}-E_kU_3(\vec{\alpha}_{k}^{n-1})\frac{ \partial \vec{\alpha}_{k}^{n-1}}{\partial \delta_k}-\\
\displaystyle &E_kc_{ek}U_{4}\frac{ \partial \vec{\alpha}_{k}^{n-1}}{\partial \delta_k},\\ \\
\displaystyle \frac{W_{1}}{\Delta t}\frac{\partial \vec{\eta}_{k}^{n}}{\partial \delta_k}&=-W_{2}\frac{\partial \vec{\eta}_{k}^{n}}{\partial \delta_k}+(S_{nk}-\delta_{nk})W'_3\biggl(\frac{\partial \vec{\alpha}_{k}^{n-1}}{\partial \delta_k}\biggr)\vec{\eta}_{k}^{n}+\biggl((S_{nk}-\delta_{nk}) W_{3}(\vec{\alpha}_{k}^{n-1})-S_{nk}W_{1}\biggr)\frac{\partial \vec{\eta}_{k}^{n}}{\partial \delta_k}\\
\displaystyle &-S_{nk}W_5\frac{\partial \vec{\alpha}_{k}^{n-1}}{\partial \delta_k}+\frac{W_{1}}{\Delta t}\frac{\partial \vec{\eta}_{k}^{n-1}}{\partial \delta_k}.
\end{cases}
\end{equation}
\\
\textbf{Linearised system for} $\mathbf{c}_{ek}\to \mathbf{c}_{ek}+\delta \mathbf{c}_{e}$:
\begin{equation}
\label{eqn:15}
\begin{cases}
\displaystyle \frac{V_{1}}{\Delta t}\frac{\partial \vec{\alpha}_{k}^{n}}{\partial c_{ek}}&=-L_kB(\vec{\alpha}_{k}^{n-1})\frac{\partial \vec{\beta}_{k}^{n}}{\partial c_{ek}}-L_kB'_{L_k}\biggl(\vec{\alpha}_{k}^{n-1},\frac{\partial \vec{\alpha}_{k}^{n-1}}{\partial c_{ek}}\biggr)\vec{\beta}_{k}^{n}+k_{nk}K(\vec{\alpha}_{k}^{n-1})\frac{\partial \vec{\eta}_{k}^{n}}{\partial c_{ek}}+\\
& \displaystyle k_{nk}K'\biggl(\vec{\alpha}_{k}^{n-1},\frac{\partial \vec{\alpha}_{k}^{n-1}}{\partial c_{ek}}\biggr)\vec{\eta}_{k}^{n}+\nu_k \biggl(V'_{5}\biggl(\frac{\partial \vec{\alpha}_{k}^{n-1}}{\partial c_{ek}}\biggr)-V'_{6}\biggl(\vec{\alpha}_{k}^{n-1},\frac{\partial \vec{\alpha}_{k}^{n-1}}{\partial c_{ek}}\biggr)\biggr)\vec{\eta}_{k}^{n} +\\
\displaystyle & \nu_k \bigl(V_{5}(\vec{\alpha}_{k}^{n-1})-V_{6}(\vec{\alpha}_{k}^{n-1})\bigr)\frac{\partial \vec{\eta}_{k}^{n}}{\partial c_{ek}} +\biggl(\nu_k\delta_k (V_7(\vec{\alpha}_{k}^{n-1})-V_1)+\biggl(\frac{1}{\Delta t} - K_T\biggr)V_{1}\biggr)\frac{\partial \vec{\alpha}_{k}^{n-1}}{\partial c_{ek}}+\\
\displaystyle &\nu_k \delta_k V'_{7}\biggl(\frac{\partial \vec{\alpha}_{k}^{n-1}}{\partial c_{ek}}\biggr)\vec{\alpha}_{k}^{n-1},\\ \\
\displaystyle U_{1}\frac{\partial \vec{\beta}_{k}^{n}}{\partial c_{ek}}&=\gamma_k^2U_{6}\frac{ \partial \vec{\alpha}_{k}^{n}}{\partial c_{ek}}+E_kc_{ek}U_{22}(\vec{\alpha}_{k}^{n})\frac{ \partial \vec{\alpha}_{k}^{n}}{\partial c_{ek}}-E_kU'_{3}\biggl(\frac{ \partial \vec{\alpha}_{k}^{n-1}}{\partial c_{ek}}\biggr)\vec{\alpha}_{k}^{n-1}-E_kU_3(\vec{\alpha}_{k}^{n-1})\frac{ \partial \vec{\alpha}_{k}^{n-1}}{\partial c_{ek}}-\\
\displaystyle &E_kc_{ek}U_{4}\frac{ \partial \vec{\alpha}_{k}^{n-1}}{\partial c_{ek}}+E_kU_{21}(P_{2}^TU_{\psi'})^{-1}\psi_1'(P_{2}^T\Phi\alpha_{k}^{n})-
E_kU_{4}\vec{\alpha}_{k}^{n-1}- E_kU_{5},\\ \\
\displaystyle \frac{W_{1}}{\Delta t}\frac{\partial \vec{\eta}_{k}^{n}}{\partial c_{ek}}&=-W_{2}\frac{\partial \vec{\eta}_{k}^{n}}{\partial c_{ek}}+(S_{nk}-\delta_{nk})W'_3\biggl(\frac{\partial \vec{\alpha}_{k}^{n-1}}{\partial c_{ek}}\biggr)\vec{\eta}_{k}^{n}+\biggl((S_{nk}-\delta_{nk}) W_{3}(\vec{\alpha}_{k}^{n-1})-S_{nk}W_{1}\biggr)\frac{\partial \vec{\eta}_{k}^{n}}{\partial c_{ek}}\\
\displaystyle &-S_{nk}W_5\frac{\partial \vec{\alpha}_{k}^{n-1}}{\partial c_{ek}}+\frac{W_{1}}{\Delta t}\frac{\partial \vec{\eta}_{k}^{n-1}}{\partial c_{ek}};
\end{cases}
\end{equation}

\bibliographystyle{plain}
\bibliography{biblio_opt} 

\end{document}